\input amstex

\pageno=1

\def\Edg {{\text {Edg}}}
\def\Stab {{\text {Stab}}}

\def\Dim {{\text {Dim}}}
\def\Mat {{\text {Mat}}}
\def\Ker {{\text {Ker}}}
\def\Vert {{\text {Vert}}}

\def\1{\hbox{\rm\rlap {1}\hskip .03in{\rm I}}}
\def\Int{{\text {Int}}}

\def\id{{\text {id}}}

\def\card{{\text {card}}}

\def\Tr {{\text {Tr}}}
\def\Tr {{\text {Tr}}}
\def\Tr {{\text {Tr}}}
\def\Tr {{\text {Tr}}}
\def\Tr {{\text {Tr}}}

\hsize=12.5cm
\vsize=19.0cm
\parindent=0.5cm
\parskip=0pt
\baselineskip=12pt
\topskip=12pt
\def\skipaline{\vskip 12pt plus 1pt}
\def\hal{\vskip 6pt plus 1pt}

\def\H#1\par{\skipaline\noindent\bf
#1\rm\par\nobreak\skipaline\nobreak\noindent}

\def\HH#1\par{\skipaline\noindent\bf
#1\rm\par\nobreak\skipaline\nobreak\noindent}

\def\Hom{{\text {Hom}}}

\def\End {{\text {End}}}
\def\Aut {{\text {Aut}}}

\def\Tr {{\text {Tr}}}

\def\Hom {{\text {Hom}}}
\def\Map {{\text {Map}}}
 
\def\mod {{\text {mod}}\,}

\def\mid {{\text {bas}}}

\def\id {{\text {id}}}

\def\Ker {{\text {Ker}}}
\def\Im {{\text {Im}}}

 \def\hal{\vskip 6pt plus 1pt}

 \centerline {  \bf   Homotopy    field  theory in dimension 2 and group-algebras}

\skipaline
\centerline {Vladimir  Turaev}

 \skipaline
 \centerline {  \bf   Abstract}
\skipaline
 
We   apply the   idea  of a
topological quantum field theory  (TQFT) to maps from manifolds into  
topological spaces. This leads to a notion of a $(d+1)$-dimensional homotopy quantum field
theory (HQFT) which may be   described as a TQFT   for  closed  $d$-dimensional
manifolds and $(d+1)$-dimensional cobordisms endowed with homotopy classes of maps into a
given space. 
For a group $\pi$,  we introduce cohomological   HQFT's with target $K(\pi,1)$
derived from cohomology classes of   $\pi$ and its subgroups of finite index. The main body of
the paper is concerned with  $(1+1)$-dimensional HQFT's. We
  classify them in terms of so called crossed group-algebras.
In particular, the cohomological $(1+1)$-dimensional HQFT's over a field of characteristic 0
are classified by  simple crossed group-algebras.  We introduce two state sum models for
$(1+1)$-dimensional HQFT's and prove that the resulting HQFT's are direct sums of
rescaled cohomological HQFT's.
 We also discuss a version of the Verlinde formula in this setting.

\skipaline
\skipaline

\centerline {\bf    Introduction}

\skipaline

Topological quantum field theories (TQFT's) produce topological invariants of
manifolds using ideas   from quantum field theory. For $d\geq 0$, a
$(d+1)$-dimensional TQFT over a   commutative ring  $K$
 assigns to every closed oriented $d$-dimensional manifold $M$ a $K$-module 
$A_M$    and assigns to  every compact oriented 
$(d+1)$-dimensional cobordism  $(W,M_0,M_1)$ a   $K$-homomorphism
$\tau(W):A_{M_0}\to A_{M_1}$. These modules and homomorphisms should satisfy
a few axioms, the main axiom being  the multiplicativity of $\tau$ with
respect to gluing of cobordisms. The study of  TQFT's has been especially successful
in low dimensions $d=1$ and $d=2$.   Deep
algebraic theories come up in both cases. The $(1+1)$-dimensional TQFT's  
bijectively correspond to finite-dimensional commutative Frobenius algebras
(see  [Di], [Du]). The $(2+1)$-dimensional
TQFT's    are closely related to   quantum groups and
braided categories (see [Tu]).

In this paper we apply the basic ideas of a
TQFT   to maps from manifolds into  
topological spaces. This suggests a notion of a $(d+1)$-dimensional homotopy quantum field
theory (HQFT) which may be briefly described as a TQFT   for  closed  $d$-dimensional
manifolds and $(d+1)$-dimensional cobordisms endowed with homotopy classes of maps into a
  pointed path-connected space $X$.
The $(0+1)$-dimensional HQFT's  over a field   correspond
bijectively to finite-dimensional representations of 
   $\pi_1(X)$ or, equivalently, to finite-dimensional flat vector bundles
over $X$. 
Thus,   HQFT's may be viewed as      higher-dimensional versions 
 of flat
vector bundles. The standard notion of a TQFT may be interpreted in this language as an HQFT
with contractible target.

In the case where   $X$ is an Eilenberg-MacLane space $K(\pi,1)$ corresponding to a
group $\pi$, the homotopy classes of maps to $X$ classify   prinicipal $\pi$-bundles.
Thus, a $(d+1)$-dimensional  HQFT  with target $X$ yields  
   invariants of  prinicipal $\pi$-bundles over   closed $d$-dimensional manifolds
and $(d+1)$-dimensional cobordisms.

In the first part of the paper we discuss a general setting of HQFT's. In particular we introduce
   cohomological HQFT's determined by   
 cohomology classes of groups. For $d=1$, we
define a wider class of semi-cohomological HQFT's.

In the second part of the paper  we focus   on     algebraic structures
underlying    $(1+1)$-dimensional   HQFT's. For a group $\pi$, we  
introduce a notion  of a    $\pi$-algebra.
Briefly speaking, this     is an associative algebra   $L$ endowed with a splitting     
$L=\bigoplus_{\alpha\in \pi} L_{\alpha}$ such that $L_{\alpha}   L_{\beta} \subset L_{\alpha
\beta}$ for any $\alpha, \beta\in \pi$. We  introduce various classes of such group-algebras
including
 so-called crossed, semisimple, biangular, and non-degenerate group-algebras.
 The  
crossed group-algebras play the key role in the study of $(1+1)$-dimensional HQFT's. Our main
result is that  each  $(1+1)$-dimensional HQFT with target $X$  has an underlying   crossed
$\pi_1(X)$-algebra and, moreover,  this establishes a  bijection between   
  (the isomorphism classes of)
  $(1+1)$-dimensional HQFT's with target $X$
 and   crossed $\pi_1(X)$-algebras.
For $\pi=1$, we obtain the well known equivalence between the 
$(1+1)$-dimensional TQFT's and commutative Frobenius algebras
(see  [Di], [Du]).

The   semi-cohomological $(1+1)$-dimensional HQFT's (arising from  
2-cohomology of groups) seem to be especially
important since they satisfy a  version of the Verlinde formula.   The
underlying  crossed group-algebra of a semi-cohomological $(1+1)$-dimensional HQFT is    
  semisimple;   we establish the converse     provided  the ground ring $K$ is a field of
characteristic 0. 

In the third part of the paper we discuss lattice    models for $(1+1)$-dimensional
HQFT's with target $K(\pi,1)$. We introduce two lattice models derived  from   
   biangular  (resp. non-degenerate) $\pi$-algebras.  The first
model  generalizes the well known lattice model for $(1+1)$-dimensional TQFT's, see [BP],
[FHK]. We first present a   map  from a surface to
$K(\pi,1)$ by a $\pi$-system, i.e., a system  of elements of $\pi$ associated with
1-cells of a    CW-decomposition of the surface. We fix  a  biangular   $\pi$-algebra and use it
to define a partition function (or a state sum) of each $\pi$-system.  This
partition function    is  homotopy invariant and   determines a $(1+1)$-dimensional
HQFT's with target $K(\pi,1)$.  The   model   based on a non-degenerate
$\pi$-algebra is more subtle: generally speaking, the corresponding   partition functions are not
homotopy invariant. To obtain an HQFT  we   sum up  the partition functions  over all
$\pi$-systems in a given homotopy class. 

We prove that   the lattice  $(1+1)$-dimensional HQFT's over an algebraically closed field of
characteristic 0    are semi-cohomological. This
 implies     that these HQFT's  satisfy the
Verlinde formula.

 This paper is a preliminary step towards a study of similar phenomena for $d=2$ and $d=3$. 
  
Throughout the paper, the symbol $K$ denotes a  commutative ring with unit.
The symbol $\pi$ denotes a group.

\skipaline \centerline  {\bf Part I.  Homotopy  quantum field  theories}

\skipaline \centerline  {\bf 1.  Basic definitions and examples}

\skipaline \noindent {\bf   1.1.  Preliminaries.} 
We shall use the language of   pointed homotopy theory. A
topological space is {\it pointed} if all its connected components are provided
with base points.  
A map  between pointed spaces   is a continuous map  sending base points
into base points. Homotopies  of such maps are always supposed to be constant on
the base points.  We shall work in the topological category although all our definitions apply in
the smooth and piecewise-linear categories. Thus, by   manifolds    we shall mean   topological
manifolds.

Let $X$ be   a    path-connected topological space  with base point $x\in X$.
We call an $X$-manifold 
any pair (a   pointed closed  oriented  manifold $M$, a map   
$g_M:M\to X$). The   map $g_M$ is called the
{\it characteristic map}. It 
   sends  the base points of  all    components of $M$ into $x$.
 It is clear that a disjoint union of $X$-manifolds is an $X$-manifold.
 An empty set $\emptyset$ is
considered as an $X$-manifold of any given dimension.
An $X$-homeomorphism
of $X$-manifolds $f:M\to M'$ is an   orientation
preserving   homeomorphism sending the base
points of $M$ onto those of $M'$ and such that $g_M=g_{M'}f$ where
$g_M, g_{M'}$ are the characteristic maps of $M,M'$, respectively.

By a
cobordism we shall mean  a triple $(W,M_0,M_1)$ where $W$ is a
compact  oriented   manifold  whose boundary is 
a disjoint union of    pointed closed  oriented  manifolds $M_0,M_1$ such that the
orientation of 
  $ M_1$ (resp. $M_0$) is induced by the one of $W$ (resp. is opposite to the one induced
from $W$).   The manifold $W$ itself is not supposed to be pointed.  

An
$X$-cobordism is a cobordism   $(W,M_0,M_1)$ endowed with a  
map  $W\to X$  sending the base points of the boundary components into $x$.
Both  bases $M_0$ and $M_1$ are considered as  
  $X$-manifolds with  characteristic maps   obtained by restricting the
 given map $W\to X$.     If  $(W,M_0,M_1)$  is an  $X$-cobordism, then 
 $(-W,M_1,M_0)$   is  also an
$X$-cobordism where  $-W$ denotes $W$ with opposite orientation.  
An $X$-homeomorphism
of $X$-cobordisms $f:(W,M_0,M_1)\to (W',M'_0,M'_1)$ is an   orientation
preserving  homeomorphism inducing $X$-homeomorphisms $M_0\to M'_0, M_1\to
M'_1$  and such that $g_W=g_{W'}f$ where $g_W, g_{W'}$ are the characteristic
maps of $W,W'$, respectively.

We can glue $X$-cobordisms along the  bases. If 
$(W_0,M_0,N), (W_1,N',M_1)$ are   $X$-cobordisms 
and $f:N\to N'$ is an $X$-homeomorphism then  the gluing of
$W_0$ to $W_1$ along $f$ yields a new $X$-cobordism with bases $M_0$ and $M_1$.
Here it is essential that $g_N=g_{N'}f$.

\skipaline \noindent {\bf   1.2.  Definition of HQFT's.}  Fix   an integer $d\geq 0$ and  a
path-connected topological space $X$ with base point $x\in X$.   We define     a   
$(d+1)$-dimensional homotopy quantum field theory     $(A,\tau)$ with target 
$X$. It will take values in the category of   projective $K$-modules of  finite type 
( = direct summands of $K^n$ with $n=0,1,...$).  The reader  
may restrict himself/herself to the case where $K$ is a field so that 
projective $K$-modules of  finite type are just finite-dimensional vector spaces over $K$.

A  $(d+1)$-dimensional HQFT   $(A,\tau)$   with target $X$ assigns
a projective $K$-module  of  finite type $A_M$ to any $d$-dimensional $X$-manifold  $M$, a 
 $K$-isomorphism $f_{\#}:A_M\to A_{M'}$ to any $X$-homeomorphism of 
$d$-dimensional $X$-manifolds
$f:M\to M'$, and a   $K$-homomorphism  $\tau(W): A_{M_0}\to A_{M_1} $ to any  
$(d+1)$-dimensional $X$-cobordism     $(W,M_0,M_1)$. These modules and
homomorphisms should satisfy the following  eight axioms.

(1.2.1) For any $X$-homeomorphisms of 
$d$-dimensional $X$-manifolds $f:M\to M', f':M'\to M''$, we have
$(f'f)_{\#}=f'_{\#}f_{\#}$. The isomorphism $f_{\#}:A_M\to A_{M'}$  is invariant
under   isotopies of  $f$ in the class of $X$-homeomorphisms. 

(1.2.2) For any disjoint $d$-dimensional $X$-manifolds $M,N$, there
is a natural  isomorphism $A_{M\amalg N}=A_M\otimes A_N$ where
$\otimes$ is the tensor product over $K$.

(1.2.3)  $  A_{\emptyset}=K$.

(1.2.4) The  homomorphism  $\tau $ associated with $X$-cobordisms is natural
with respect to $X$-homeomorphisms.

(1.2.5) If a $(d+1)$-dimensional $X$-cobordism $W$    is a disjoint
union of two $X$-cobordisms $W_1,W_2$ then    $\tau(W)=\tau(W_1)  \otimes   \tau(W_2)$.

(1.2.6) If an $X$-cobordism     $(W,M_0,M_1)$   is obtained from
two   $(d+1)$-dimensional $X$-cobordisms $(W_0,M_0,N)$ and $(W_1,N',M_1)$ by gluing 
along an $X$-homeomorphism $f:N\to N'$ then
$$\tau(W)= \tau(W_1)\circ f_{\#} \circ \tau(W_0):A_{M_0}\to A_{M_1}.$$

(1.2.7) For any   $d$-dimensional
$X$-manifold
 $(M,g :M\to X)$  and any map $F:M\times [0,1] \to X$
such that $F\vert_{M\times 0}=F\vert_{M\times 1}=g $ and $F(m\times
[0,1])=x$ for all base points $m\in M$,
we have $\tau(M\times [0,1], F)=\id: {A_M}\to {A_M}$ where the cylinder $ M\times
[0,1] $ is viewed as an $X$-cobordism  with bases   $M\times  0=M, M\times
1=M$ and characteristic map $F$.

(1.2.8) For any   $(d+1)$-dimensional
$X$-cobordism 
 $W=(W,g:W\to X)$, the homomorphism $\tau(W)$ is preserved
under any homotopy of   $g$ relative to   
$\partial W$.

Axioms  (1.2.1) - (1.2.7) form a version  of the standard definition of a TQFT, cf.  [At],  
[Tu, Chapter III]. It is sometimes convenient to consider the homomorphism
$\tau(W)  $ associated to an
  $X$-cobordism     $(W,M_0,M_1)$
as a vector 
$$ \tau(W)\in \Hom_K (A_{M_0}, A_{M_1})= A_{M_0}^*\otimes A_{M_1}.$$
In this language,   axiom (1.2.6) says that 
$\tau(W)$ is obtained  from $ \tau(W_0)\otimes  
 \tau(W_1) $ by the tensor contraction induced by the pairing 
$a\otimes b \mapsto  b(f_{\#}(a)):  A_{N}  \otimes A_{N'}^*\to K$.

Any closed oriented $(d+1)$-dimesional
manifold $W$ endowed with a map $ g:W\to X$ can be considered as a cobordism
with empty bases. The corresponding $K$-linear endomorphism of
$A_{\emptyset}=K$ is multiplication by a certain    $\tau(W)\in K$. By  
(1.2.8), $\tau(W)$ is a homotopy invariant of   $g$. More generally,
 the modules and homomorphisms
provided by any HQFT   depend only on the homotopy classes of the
characteristic maps, see Section 2.1.

We define a few simple operations on $(d+1)$-dimensional HQFT's   with    target   $X$. The direct
sum of    HQFT's $(A,\tau) \oplus (A',\tau')$  is defined by $(A\oplus A')_M=A_M\oplus
A'_M$ and $(\tau\oplus \tau')(W)=\tau(W)\oplus \tau'(W)$. The tensor product is
defined similarly using $\otimes$ instead of $\oplus$. 
The dual  $(A^*,\tau^*)$ of an  HQFT $(A,\tau)$ is defined by 
$A_M^*=\Hom_K (A_M,K)$ for any $X$-manifold $M$ with action of
$X$-homeomorphisms obtained  by transposition from the one given by  
$(A,\tau)$.  
We define $\tau^*(W): A_{M_0}^*\to A_{M_1}^*$ as the transpose of
$\tau(-W):A_{M_1}\to  A_{M_0}$.  All the axioms of an HQFT are straightforward.

We   define a category denoted  $Q_{d+1}(X,x)$ (or shorter  $Q_{d+1}(X)$) whose objects are
$(d+1)$-dimensional HQFT's  with     target  
$X$. A morphism
 $(A,\tau)\to  (A',\tau')$ in this category 
      is a family of $K$-homomorphisms
$\{\rho_M:A_M\to A'_M\}_M$ where $M$ runs over   $d$-dimensional $X$-manifolds  such
that:   $\rho_{\emptyset}=\id_K $;
 for  disjoint   $X$-manifolds $M,N$, we have 
$\rho_{M\amalg  N} =\rho_M\otimes \rho_N$;     the natural square
diagrams associated with homeomorphisms  of $X$-manifolds and with $X$-cobordisms are
commutative.  It can be shown (we shall not use it) that all morphisms in the category  
$Q_{d+1}(X,x)$ are isomorphisms.

It is easy to deduce from definitions that  the isomorphism  classes of
$(0+1)$-dimensional HQFT's with target
$X$  correspond bijectively to the isomorphism  classes of linear actions of 
$
\pi_1(X,x)$ on projective $K$-modules of finite type.  
Under this correspondence, the   invariant
$\tau$ of a map $g:S^1\to X$ equals   the trace of the conjugacy class of linear endomorphisms 
determined by 
$g$.  

If the space $X$    consists
of only one point $x$ then  all references to maps
into $X$ are redundant and we obtain   the usual definition of a 
 topological quantum field theory (TQFT)
for pointed closed oriented $d$-dimensional manifolds  and their   cobordisms. 
Restricting any HQFT $(A,\tau)$ with target $(X,x)$ 
to those
$X$-manifolds and $X$-cobordisms whose
    characteristic map is a constant map into $x\in X$,
we obtain an underlying TQFT of $(A,\tau)$. 

 \skipaline \noindent {\bf   1.3.  Primitive cohomological HQFT's.} Let $X$ be a    
Eilenberg-MacLane space of type $K(\pi,1)$ where $\pi$ is a group. (Speaking about
Eilenberg-MacLane spaces we always assume that they are CW-complexes). 
Recall that the symbol  $K^* $ denotes the  multiplicative group consisting of the invertible
elements of
$K$.
For each $\theta\in    H^{d+1}(X;K^*)=H^{d+1}(\pi;K^*)$, we shall define 
a 
$(d+1)$-dimensional HQFT $(A,\tau)$ with target  $X$ called the
{\it primitive cohomological HQFT} associated with $\theta$ and denoted $(A^{\theta},
\tau^{\theta})$. 
This construction  is inspired by the work of Freed and
Quinn [FQ] on TQFT's associated with finite groups.

Choose  a singular
$(d+1)$-dimensional cocycle  on
$X$ with values  in  $K^* $ representing $\theta$. By abuse of notation we denote this cocycle
by the same symbol $\theta$.   Let $M$ be a   $d$-dimensional $X$-manifold. Then $A_M$ is a 
free $K$-module of rank 1 defined as follows.  
 A $d$-dimensional singular
cycle $a\in C_d(M)=C_d(M;\bold  Z)$ is said to be {\it fundamental} if  it represents   
the fundamental class   $[M]\in
H_d(M;\bold  Z)$ defined as  the sum of the fundamental  classes of the
components of $M$.
Every  fundamental
cycle $a\in C_d(M)$ 
determines a non-zero element $\langle a\rangle\in A_M$. 
If $a,b\in C_d(M)$ are two  fundamental cycles, then we impose the equality 
$\langle a\rangle=g^*(\theta) (c) \langle b\rangle$
where $g:M\to X$ is the characteristic map of $M$ and $c$ is a 
$(d+1)$-dimensional singular  chain in $M$ such that $\partial c=a-b$.
Note that $g^*(\theta) (c)\in K^*$ does not depend on the choice of $c$:
if $ \partial c=\partial c'$ with $c,c'\in C_{d+1}(M)$
then   $c-c'=\partial
e$ with $e\in C_{d+2}(M)$ and   $$g^*(\theta) (c)/g^*(\theta) (c')= g^*(\theta)
(\partial e)=   \theta  (\partial g_*(e))=
  \partial \theta  ( g_*(e))=1.$$
It is easy to check that $A_M$ is a well defined 
free $K$-module of rank 1.  
An $X$-homeomorphisms of $d$-dimensional $X$-manifolds
$f:M\to M'$  induces an isomorphism $A_M\to A_{M'}$ sending the generator $\langle
a\rangle\in A_M$ as above  into $
\langle
 f_*(a) \rangle\in A_{M'}$.  
 
Consider 
a $(d+1)$-dimensional  $X$-cobordism     $(W,M_0,M_1)$. Let
 $B\in C_{d+1}(W)=C_{d+1}(W;\bold  Z)$ be a fundamental chain, i.e., a singular chain
representing the fundamental class $[W]\in H_{d+1}(W, \partial W;\bold  Z)$  defined as
the sum of the fundamental classes of the components of $W$. It is clear that
$\partial B=-b_0+b_1$ where $b_0,b_1$ are fundamental cycles of $M_0, M_1$,
respectively. We define $\tau(W): A_{M_0}\to A_{M_1} $ by
 $$\tau(W) (\langle b_0\rangle)=(g^*(\theta) (B))^{-1} \, \langle b_1\rangle$$
where $g:W\to X$ is the characteristic map of $W$. Let us verify that 
$\tau(W)$ does not depend on the choice of $B$. 
Let
 $B' \in C_{d+1}(W)$ be  another  fundamental chain     with $\partial B'=-b'_0+b'_1$
where $b'_i$ is a fundamental cycle  of $M_i$ with $i=0,1$.
Choose $c_i\in C_{d+1} (M_i)$ 
 such that $\partial c_i=b_i-b'_i$ for $i=0,1$.
By definition, 
$\langle b_i\rangle=g^*(\theta) (c_i) \langle b'_i\rangle$
where $g_i:M_i\to X$ is the characteristic map of $M_i$.
To see that $\tau(W)$ is well defined it suffices to check the equality in $K^*$
$$g^*(\theta) (c_0) (g^*(\theta) (B'))^{-1}=g^*(\theta) (c_1) (g^*(\theta)
(B))^{-1}.\leqno (1.3.a)$$ Clearly,  
$B+c_0-B'-c_1  \in C_{d+1}(W)$  is a cycle
representing  $0$  in $H_{d+1}(W, \partial W;\bold  Z)$. 
Note that  the inclusion $H_{d+1}(W;\bold  Z) \to H_{d+1}(W, \partial W;\bold  Z)$
is injective. Therefore the cycle 
$B+c_0-B'-c_1$ is   a boundary in $W$. This implies 
(1.3.a).

It remains to verify  the axioms of an HQFT.  Axioms (1.2.1) - (1.2.5) are
straightforward.
Let us check (1.2.6). Let     $(W,M_0,M_1)$ be   a $(d+1)$-dimensional 
$X$-cobordism     obtained from two  $X$-cobordisms $(W_0,M_0,N)$ and
$(W_1,N',M_1)$ by gluing  along an $X$-homeomorphism $f:N\to N'$.
Let
 $B_0\in C_{d+1}(W_0)$ be a  fundamental chain  
with  $\partial B_0=-b_0+b $
where $b_0,b $ are fundamental cycles of $M_0, N$, respectively.
Clearly, $f_*(b)$ is a fundamental cycle of $N'$. Choose
a fundamental chain $B_1\in C_{d+1}(W_1)$ such that   $\partial B=-f_*(b)+b_1 $
where $b_1 $ is a fundamental cycle  of $M_1$.
By definition, the composition 
$\tau(W_1)\circ f_{\#} \circ \tau(W_0)$ sends the generator
$\langle b_0 \rangle$ of $A_{M_0}$
into $(g_0^*(\theta) (B_0) g_1^*(\theta) (B_1))^{-1} \langle b_1\rangle$
where $g_j$ is the characteristic map of $W_j$ for $j=0,1$. 
Observe that under the gluing of $W_0$ to $W_1$ along
$f:N\to N'$, the chain $B_0+B_1$ is mapped into
a   fundamental chain $B\in C_{d+1}(W)$ with   $\partial B=-b_0+b_1 $.
Therefore, $g_0^*(\theta) (B_0) g_1^*(\theta)
(B_1)=g^*(\theta) (B) $ where $g:W\to X$ is the characteristic map of $W$. 
By definition,
$$\tau(W) (\langle b_0 \rangle)= (g^*(\theta) (B))^{-1} \langle b_1 \rangle
= (\tau(W_1)\circ f_{\#} \circ \tau(W_0)) (\langle b_0 \rangle).$$

Let us check
(1.2.7). It is here that we   use the fact that we work in the
pointed category and that   $X= K(\pi,1)$.  Consider a $d$-dimensional  $X$-manifold
 $(M,g:M\to X)$  and a  map $F:M\times [0,1] \to X$
such that $F\vert_{M\times 0}=F\vert_{M\times 1}=g$ and $F(m\times
[0,1])=x$ for all   base points $m$ of the
components of $M$. We  choose
  a fundamental chain $B\in C_{d+1}(M\times [0,1])$ so that $\partial
B=-(b\times 0)+(b\times  1)$ where $b $ is a fundamental cycle  of $M $.
By definition, the homomorphism
$\tau(M\times [0,1], F)$ sends the generator $\langle b\rangle\in A_M$ into 
$(F^*(\theta) (B))^{-1} \langle b \rangle$. It remains to prove that 
$F^*(\theta) (B)=1$.

Consider the   map  $p:M\times [0,1] \to M\times S^1$ obtained by the gluing
$M\times 0=M=M\times 1$.  It is clear that $p_*(B)$ is a fundamental cycle in
$M\times S^1$. The map $F$ induces a map $\tilde F:M\times S^1\to X $ such
that $F=\tilde F \circ p$. Therefore 
$F^*(\theta) (B)=\theta(F_*(B))=\theta (\tilde F_*([M\times S^1]))$
where $\tilde F_*([M\times S^1])\in H_{d+1}(X;\bold  Z)$. We claim that the latter
homology class is trivial.  This would imply the equality $F^*(\theta) (B)=1$.
To prove our claim it suffices to observe that the map  $\tilde F$ extends to a
map $M\times D^2 \to X$ where $D^2$ is a 2-disc bounded by $S^1$.
This follows easily from the assumptions that $X=K(\pi,1)$ and each 
component of $M$ contains a point $m$ such that $\tilde F(m\times S^1)=x$.

Let us verify (1.2.8). 
Let     $(W,g:W\to X)$ be   an $X$-cobordism of dimension $d+1$. We should verify that
for a  fundamental chain 
 $B\in C_{d+1}(W )$, the element $g^*(\theta) (B)\in K^*$ is preserved under 
   any homotopy of $g$ relative to
$\partial W$.  Consider the manifold $\tilde W$
obtained from $W\times [0,1]$ by contracting each interval $w\times [0,1]$
with $w\in \partial W$ to a point.  The homotopy of $g$ induces a map $\tilde g:\tilde
W\to X$.  The manifold $\partial \tilde W$ is  obtained by gluing 
$W\times 0$ to $W\times 1$ along  
 $\partial W\times 0= \partial W=\partial W\times 1$.
The chain  $(B\times 0) - (B\times 1)$ in $\partial \tilde W$ represents the
fundamental class of $\partial \tilde W$ and therefore bounds a singular chain in
$\tilde W$.  This implies that ${\tilde g}^*(\theta)(B\times 0)=
{\tilde g}^*(\theta)(B\times
1)$.  This is exactly the equality we need.
 
It is easy to compute the element $\tau(W,g)\in K$ associated by this HQFT with 
    a  (closed oriented) $(d+1)$-dimensional
$X$-manifold $(W,g:W\to X)$:
$$\tau(W,g)=g^*(\theta) ([W])=\theta(g_*([W]))\in K^*.$$

Note finally that the primitive cohomological HQFT's arising as above from   different  singular
cocycles  representing  $\theta \in  H^{d+1}(X;K^*)=H^{d+1}(\pi;K^*)$ are isomorphic. 
Therefore   the isomorphism
class of these HQFT's    depends only on    $\theta$. It is
easy to describe explicitly the 
 primitive cohomological HQFT  associated with $0\in H^{d+1}(\pi;K^*)$: this HQFT
assigns $K$ to all $d$-dimensional $X$-manifolds and assigns $\id_K:K\to K$ 
to all $X$-homeomorphisms  and to all $(d+1)$-dimensional $X$-cobordisms.

\skipaline \noindent {\bf   1.4.  Rescaling of  HQFT's.} Further examples of
HQFT's can be obtained from additive invariants of $X$-cobordisms. An integer valued function 
$\rho$ of $(d+1)$-dimensional $X$-cobordisms is   an {\it  additive invariant} if it is preserved
under 
  $X$-homeomorphisms and   homotopies of the characteristic map and is additive under the
gluing of
$X$-cobordisms described in Section 1.1.  Fix $k\in K^*$.  Using  an additive invariant $\rho$,
we can transform
  any $(d+1)$-dimensional HQFT $(A,\tau)$ into a $k^{\rho }$-rescaled HQFT which coincides
with $(A,\tau)$ except that the homomorphism associated with a $(d+1)$-dimensional 
$X$-cobordism     $ (W,M_0,M_1)$ is equal to $k^{\rho(W)} \tau(W)$.

For any $d\geq 0$, an   additive invariant 
of a  $(d+1)$-dimensional 
$X$-cobordism     $ (W,M_0,M_1)$ is  given by     $\rho (W) =\chi (W,M_0)$
where $\chi$ is  the Euler characteristic. 
This example  can be refined using the semi-characteristic. 
Consider for concreteness the case $d=1$. Set 
$$\rho_0 (W) = (\chi (W) + b_0(M_0)- b_0(M_1))/2 \in \bold Z \leqno (1.4.a)$$
 where $b_0(M)$ is
the number of components of $M$.  It is obvious that  $\rho_0$ is an    additive invariant of
$2$-dimensional $X$-cobordisms.

For
$d=0\, (\mod 3)$, examples of additive invariants are provided by the signature  of $W$   and
the $G$-signatures of $W$ determined by the 
 homomorphism   $ 
\pi_1(W)\to \pi_1(X)$ induced by the characteristic map $W\to X$ and a fixed homomorphism
from
$\pi_1(X) $ into a finite group $G$.

 \skipaline \noindent {\bf   1.5.  Transfer.} 
 Let $X$ be a  connected   CW-complex with base point $x\in X$.
Let $p:E\to X$ be a  connected finite-sheeted covering   with base point  
 $e\in
p^{-1}(x)$.  From any $(d+1)$-dimensional HQFT $(A,\tau)$ with target 
$E$  we shall derive a $(d+1)$-dimensional HQFT  $(\tilde A, \tilde \tau)$  with target  
$X$.   It is called the {\it transfer} of $(A,\tau)$. 

 We first fix a map $q:E\to E$ which
is homotopic to the identity and sends the   set  
$p^{-1}(x)$ into the point $e$.  (To construct such a map $q$ one may choose for each  $y\in
p^{-1}(x)$ a path  
from  $y$ to $e$  in $E$ and then push all  $y\in p^{-1}(x)$ into $e$ 
along these paths. This extends to a homotopy of the identity map into   $q$).
Let $(M,g:M\to X)$ be a $d$-dimensional  $X$-manifold. There is a finite number of
lifts of $g$ to $E$, i.e., of maps $\tilde g:M\to E$ such that
$p \tilde g=g$. Note that   each pair $(M, q\tilde g)$ is an $E$-manifold.
Consider  the $K$-module 
$$\tilde A_M= \bigoplus_{\tilde g, p \tilde g=g} A_{(M, q\tilde g)}.$$
 Any $X$-homeomorphism of 
$d$-dimensional $X$-manifolds
$f:(M,g)\to (M',g')$ induces a $K$-isomorphism $f_{\#}:\tilde A_{(M,   g)}\to \tilde
A_{(M',  g')}$ as follows.  Since $g f=g'$, any
lift $\tilde g:M\to E$ of $g$ induces a lift $ \tilde g f:M'\to E$ of   $g'$. The
HQFT $(A, \tau)$ yields an isomorphism $f_{\#}: A_{(M, q\tilde g)}\to  A_{(M' ,
q\tilde g f)}$. The direct sum of these isomorphisms yields the desired 
isomorphism $f_{\#}:\tilde A_{(M,   g)}\to \tilde A_{(M',  g')}$. 

Now, consider a $(d+1)$-dimensional    $X$-cobordism     $(W,M_0,M_1)$   with
characteristic map $g:W\to X$. As above, there is a finite number of lifts $\tilde
g:W\to E$   such that $p \tilde g=g$. Each such lift $\tilde g$ can be
restricted to the bases of $W$ and induces in this way certain lifts,
$  \tilde g_0, \tilde g_1$,  of the characteristic maps $ {M_0}\to X, M_1\to X$. 
Then the HQFT $(A, \tau)$ yields a homomorphism
$$\tau(W,q\tilde g): A_{(M_0, q\tilde g_0)}\to A_{(M_1, q\tilde g_1)}.$$
The   sum of these homomorphisms corresponding to all
  $\tilde g$ defines a homomorphism $\tilde \tau(W,g): \tilde A_{M_0}\to \tilde
A_{M_1} $.  (Warning:  a
lift of   the characteristic map  $ {M_0}\to X$  to
 $E$ may extend to different lifts of
$g $ so
 that the sum in question is in general not a direct sum.)
It is easy   to verify that these definitions yield a
$(d+1)$-dimensional HQFT $(\tilde A, \tilde \tau)$
with target $X$. 
In particular, for   a   (closed oriented)  $(d+1)$-dimensional
$X$-manifold $(W,g:W\to X)$,
$$\tilde \tau(W,g)=\sum_{\tilde g:W\to E, p\tilde g=g}  
 \tau (W, q\tilde g)\in  K.$$
We leave it as an exercise for the reader to check that the isomorphism class of the HQFT 
 $(\tilde A, \tilde \tau)$ does not depend on the choice of the map $q$ (cf. Section 2.3).
 
\skipaline \noindent {\bf   1.6.  Cohomological and semi-cohomological HQFT's.}  Let $\pi$ be a
group and  $X$ be  an Eilenberg-MacLane
space of type $K(\pi,1)$ with base point $x$. Let
$G\subset \pi$ be a subgroup of finite index $n$ and 
$p:E\to X$ be the connected $n$-sheeted covering  of $X$ corresponding to $G$ with base point
$e\in p^{-1}(x)$.  
For each $\theta\in     H^{d+1}(G;K^*)=H^{d+1}(E;K^*)$, the 
transfer transforms    the primitive cohomological
 HQFT with target $E=K(G,1)$ 
associated with 
  $\theta $   into  a $(d+1)$-dimensional HQFT 
 with target $X$  called the
{\it   cohomological HQFT} associated with $\theta$ and denoted by $(A^{\pi,G, \theta},
\tau^{\pi,G, \theta}) $.
  It follows from definitions that  for a  connected $d$-dimensional
$X$-manifold $(M, g:M\to X)$ with base point $m\in M$, the $K$-module $ A^{\pi,G, \theta}_M$ is free
of rank $$ 
\card\, \{ i\in G\backslash\pi
\,\vert\,
\Im (g_{\#}: \pi_1(M,m)\to \pi_1(X,x)) \subset i^{-1}Gi \}  $$
(this rank does not depend on $\theta$).  
To compute the invariant $ \tau^{\pi,G, \theta}(W)$ for  a  (closed oriented) connected 
$(d+1)$-dimensional $X$-manifold $(W, g:W\to X)$, we first deform $g$ so that  it sends   a  
point $w\in W$ into $x$. For each right coset class $i\in G\backslash \pi$, choose a
representative    $\omega_i\in    i$ so that $i=G\omega_i   \subset \pi$.
 Then
$$ \tau^{\pi,G, \theta}(W)= \sum_{ {{i\in   G\backslash \pi}},\,  \,
  g_{\#}(\pi_1(W,w)) \subset i^{-1}Gi  }
g_i^*(\theta) ([W])\,\,\in K$$
where $g_i:(W, w) \to  (E,e)$ is a map inducing the homomorphism
$\omega_i g_{\#}  \omega_i^{-1}$ of the fundamental groups. (The value   $g_i^*(\theta) ([W])
\in K^* $ does not depend on the choice of $\omega_i$.) 

We define   semi-cohomological HQFT's as direct sums of rescaled
cohomological HQFT's.
 Specifically, 
a $(1+1)$-dimensional HQFT with target $X=K(\pi,1)$ is said to be {\it semi-cohomological} if
it splits as a finite direct  sum $\oplus_{i} (A_i,\tau_i)$   where  each $(A_i,\tau_i)$ is a
$(1+1)$-dimensional 
HQFT with target $X$ obtained by  $k_i^{\rho_0 }$-rescaling from a cohomological HQFT
$(A^{\pi,G_i, \theta_i}, \tau^{\pi,G_i, \theta_i})$ where $k_i\in K^*$, $\rho_0$ is the additive
invariant of 2-dimensional cobordisms defined by  (1.4.a), $G_i\subset \pi$ is a subgroup of
finite index, and $\theta_i\in     H^{2}(G_i;K^*)$.

\skipaline \noindent {\bf   1.7.  Hermitian and unitary HQFT's.}
   Assume that  the ground ring  $K$ has a ring  involution $k\mapsto
\overline k:K \to K$.  
Let  $d\geq 0$ and  $(A,\tau)$ be a $(d+1)$-dimensional HQFT     with
target      a pointed path-connected  
space  $X$.    
A Hermitian structure on   $(A,\tau)$ 
 assigns to 
each   $d$-dimensional $X$-manifold  $M$ a non-degenerate Hermitian pairing $\langle
.,. \rangle_M:A_M\otimes_K A_M \to K$ satisfying the following two conditions.

\hal  (1.7.1)  The pairing $\langle .,. \rangle_M$ is natural with respect to
$X$-homeomorphisms and multiplicative with respect to disjoint union; for
$M=\emptyset$ the pairing  $\langle .,. \rangle_M$ is determined by the unit $1\times
1$-matrix.

\hal (1.7.2)  For any $X$-cobordism $(W,M_0,M_1)$ and any $a\in A_{M_0}, b\in 
A_{M_1}$, we have 
$$\langle \tau(W) (a), b\rangle_{M_1}=\langle a,\tau(-W) (b) \rangle_{M_0}.\tag
1.7.a$$  

An HQFT   with a Hermitian structure  
 is called a Hermitian
HQFT.  If $K=\bold  C$ with   usual complex conjugation and
 the   Hermitian form $\langle .,. \rangle_M$ is positive  definite for every $M$, we say
that the Hermitian HQFT is unitary.

Note two properties of a Hermitian HQFT $(A,\tau)$. First, 
if $W$ is a $(d+1)$-dimensional closed $X$-manifold  then
$\tau(-W)=\overline {\tau(W)}$.   Secondly, the
group of   $X$-self-homeomorphisms of any 
$d$-dimensional closed $X$-manifold $M$ acts in $A_M$ preserving the Hermitian form
$\langle .,. \rangle_M$. For a unitary HQFT  we obtain a unitary
action.  
 This implies in particular the following   estimate for the  value of
 $\tau$ on the mapping torus $W_g$ of an   $X$-homeomorphism
$g:M\to M$:
if  $(A,\tau)$ is  a unitary HQFT  then  $\vert \tau(W_g)\vert \leq
\dim_{\bold  C} A_M$ (cf.  [Tu],
Chapter III).

It is   easy  to check that direct sums, tensor products,     and 
transfers of  Hermitian (resp. unitary) HQFT's are again   Hermitian (resp.
unitary) HQFT's.  We can define   a category, $HQ_{d+1}(X)$ (resp. 
$UQ_{d+1}(X)$),  whose
objects are
$(d+1)$-dimensional Hermitian (resp. unitary)  HQFT's  with     target  
$X$. The morphisms   in this category are defined as in Section 1.2 with
  additional requirement that the homomorphisms
 $\{\rho_M \}_M$ preserve  the Hermitian pairing. 

The construction of primitive cohomological
 HQFT's can be refined to yield Hermitian and unitary
HQFT's.   Consider the
multiplicative group $S= \{k\in K^*\,\vert\,
k \overline k =1\} \subset K^*$.  For each $\theta \in  H^{d+1}(\pi;S)$, the
$(d+1)$-dimensional HQFT $(A^\theta,\tau^\theta)$ defined  in Section 1.3 can be 
provided with a  Hermitian structure  as follows. For a  $d$-dimensional
$X$-manifold $M$ (where $X=K(\pi,1)$) and a 
fundamental
cycle $a\in C_d(M)$,  set 
$\langle \langle a\rangle , \langle a\rangle \rangle_M=1\in K$ where
$\langle a\rangle\in A_M$ is  the   vector  represented by $a$.
This yields a well defined Hermitian form on   $A_M=K  \langle a\rangle$  satisfying   (1.7.1)
and (1.7.2). The isomorphism
class of the resulting  Hermitian HQFT  depends only on    $ \theta$.  If
$K=\bold  C$ then $S=S^1=\{z\in \bold  C\,\vert\, \vert z \vert =1\}$ and the Hermitian HQFT    
determined by    any $ \theta \in  H^{d+1}(\pi;S^1)$ is unitary. 

 Note    that   the     isomorphism classes of $(0+1)$-dimensional   unitary  HQFT's  
with target   $X$  correspond bijectively to  the  isomorphism classes of 
 finite-dimensional
flat   unitary bundles over $X$.

\skipaline \centerline  {\bf 2.  Homotopy properties of HQFT's}

 \skipaline \noindent {\bf   2.1.  Homotopy invariance of HQFT's.} 
It was already mentioned above that the modules and homomorphisms
provided by an HQFT   depend only on the homotopy classes of the
characteristic maps. Here we discuss this   in detail.

Fix a
path-connected topological space $X$ with base point $x$. Fix 
a $(d+1)$-dimensional HQFT $(A,\tau)$ with target $X$.  For a
pointed space $Y$, we  shall  denote the set of homotopy classes of maps (in the
category of pointed spaces) $YÊ\to X$ by $\Map(Y,X)$.

We  introduce    
homotopy $X$-manifolds    exactly as
$X$-manifolds  with the difference that instead of maps to $X$  we
  speak of homotopy classes of maps.  Thus, a homotopy $X$-manifold  is 
a  pair (a   pointed closed  oriented  manifold $M$, an element 
$G\in \Map(M,X)$). For a homotopy $X$-manifold $(M,G)$ we define a $K$-module 
$A_{(M,G)}$ as follows.
Observe that any        homotopy $F:M\times [0,1]
\to X$ between  two maps  $ g,g':M\to X $   belonging to the class $G$ provides the
cylinder   $M\times [0,1] $ with the structure of an $X$-cobordism. This cobordism
gives a homomorphism $\tau_F: A_{(M,g)}\to A_{(M,g')}$.  It follows from (1.2.6) that  
$\tau_{FF'}=\tau_{F'} \tau_F$ for any  composable homotopies $F,F'$. Axiom (1.2.7)
implies that  for any homotopy $F_0$ relating a map to itself, $\tau_{F_0} =\id$. This
implies that  for any homotopy $F$ between $g$ and $g'$ as above, the
homomorphism $\tau_F$ is an   isomorphism     depending only
on $g,g'$ and independent of the choice of   $F$. We
identify  the modules $\{A_{(M,g)}\vert\, g\in G\}$    along the  isomorphisms
$\{\tau_F\}_F$.
This gives  a $K$-module  $A_{(M,G)}$ depending  on $(M,G)$ and canonically
isomorphic to   each $A_{(M,g)}$  with   $g\in G$.

 An $X$-homeomorphism
of homotopy $X$-manifolds $f:(M,G)\to (M',G')$ is an   orientation
preserving   homeomorphism $M\to M'$ sending the base
points of $M$ onto those of $M'$ and such that $G=G'f$ where
  the equality is understood as an equality of
homotopy classes. We define the action $f_{\#}:A_{(M,G)}\to
A_{(M',G')}$ of $f$ 
as the composition 
 $$
\CD
 A_{(M,G)}   @>=>>  A_{(M,g'f)}     @>f_{\#}>>  A_{(M',g')}   @>=>>
A_{(M',G')}
\endCD $$
where $g':M'\to X$ is a  map representing $G'$ and the first and third homomorphisms
are the canonical identifications. We claim that
this composition does not depend on the choice of $g'\in G'$.
Let $g'':M'\to X$  be another representative of the class $G'$ and let $F'$
be a homotopy between $g'$ and $g''$.
It is clear that   $F'$  
induces a homotopy $F=F'\circ (f\times \id_{[0,1]})$ between $g'f$ and $g''f$. Axiom 
(1.2.4)
  implies that  the  diagram
$$
\CD
 A_{(M,g'f)}   @>f_{\#} >> A_{(M',g')}   \\ 
  @V\tau_FVV    @VV\tau_{F'}V        \\
  A_{(M,g''f)} @>f_{\#}>>  A_{(M',g'')} 
\endCD $$
 is commutative. This   implies our claim.

We define  a homotopy $X$-cobordism as a
 cobordism   $(W,M_0,M_1)$ endowed with a homotopy classes of maps $W\to X$. 
(The maps   should send  the base points of $\partial W$ into $x$
and  the homotopies should be constant on the base points; we denote the set of the
corresponding homotopy classes by   $\Map (W,X)$).   
  For a $(d+1)$-dimensional homotopy $X$-cobordism
 $(W,M_0,M_1,G\in \Map(W,X))$ we shall define  a homomorphism 
$\tau(W,G): A_{(M_0,G_0)} \to A_{(M_1,G_1)}$ where $G_j \in \Map (M_j,X)$ is the
restriction   of $G$ to $M_j$, for $j=0,1$. 
Let      $g:W\to X$   be a representative  of the homotopy class  $G$.
Let  $g_j \in G_j$ be the restriction of $g$ to $M_j$.  We define  
$\tau(W,G)$ as the composition 
 $$
\CD
 A_{(M_0,G_0)}   @>=>>  A_{(M_0,g_0)}     @>\tau(W,g)>>  A_{(M_1,g_1)}   @>=>>
A_{(M_1,G_1)}
\endCD $$
where the first and third homomorphisms are the canonical identifications and the
second homomorphism is determined by the $X$-cobordism $(W,g)$.  We claim that
$\tau(W,G)$ does not depend on the choice of $g$ in the class $G$. Let 
$g':W\to X$ be another representative of   $G$ 
and let $F$ be a homotopy between   $g$ and $g'$.
We should prove the commutativity of the diagram
$$
\CD
 A_{(M_0,g_0)}     @>\tau(W,g)>>  A_{(M_1,g_1)}  \\ 
  @V\tau_{F_0}VV    @VV\tau_{F_1}V        \\
A_{(M_0,g'_0)}     @>\tau(W,g')>>  A_{(M_1,g'_1)}
\endCD $$
 where $F_j$ denotes  the restriction of $F$ to $M_j\times [0,1]$ for $j=0,1$.
Let $F_1^{-1}$  be the homotopy $M_1\times  [0,1]\to X$ inverse to $F_1$
(i.e., $F_1^{-1} (a,t)= F_1  (a,1-t) $ for any $a\in M_1, t\in  [0,1]$). 
Consider the $X$-cobordism $W'$ obtained by gluing   
 the $X$-cobordisms $(M_0\times  [0,1], F_0), (W,g') $, and 
$(M_1\times  [0,1], F_1^{-1})$ along $M_0\times  1 =M_0\subset \partial  W$
and $M_1\times  0 =M_1\subset \partial  W$. By axiom (1.2.6),
$\tau(W')=(\tau_{F_1})^{-1} \tau(W,g') \,\tau_{F_0}$. On the other hand, it is clear that 
$W'$ is just the same cobordism $W$ with another characteristic map to $X$. Moreover,
the homotopy $F$ induces a homotopy of this characteristic map into $g$
relative to $\partial W$. By axiom (1.2.8),
$\tau(W')=  \tau(W,g)  $. Hence the diagram above is commutative.

The constructions of this subsection show  that  from the very beginning we can 
formulate the definition of an HQFT in terms of the homotopy classes of characteristic
maps.  In the sequel we shall make no difference
between characteristic
maps and their homotopy classes.

We end this subsection with a simple but useful lemma. 
We
  say that two maps
$g,g'$ from a cobordism $W$ to $X$ are equivalent  if they are
homotopic (rel
$\partial W$) in the complement of a small ball in $\Int W$.  It is easy to see that
this is indeed an equivalence relation. Note that equivalent maps $W\to X$ coincide 
on $\partial W$.

\skipaline \noindent {\bf    2.1.1.   Lemma. }   {\sl   Let $d\geq 1$ and  $(W,M_0,M_1)$
be a $(d+1)$-dimensional cobordism. If   $g,g':W\to X$ are two equivalent maps
  then 
$$\tau (W,g)=\tau (W,g'):A_{(M_0, g\vert_{M_0})}\to A_{(M_1, g\vert_{M_1})}.$$}

\skipaline  {\sl  Proof.}  By (1.2.5), it suffices to prove the lemma for  connected $W$. If
$\partial W=\emptyset$ then we split
$W$ along a  small embedded 
$d$-dimensional sphere into a union of two cobordisms.  By (1.2.6),  
the claim for these smaller cobordisms would imply the claim for $W$. Thus,  it suffices to
consider the case where $W$ is connected and 
$\partial W\neq
\emptyset$. Suppose for concreteness that   $M_0 \neq
\emptyset$.   By assumptions, $g':W\to X$ is homotopic (rel  $\partial W$)
to a map
$g'':W\to X$ which coincides with $g$ outside a  $(d+1)$-dimensional ball 
$B\subset W$. Choose  a regular neighborhood $M_0\times [0,1]\subset W$
of
$M_0$  such that $B\subset M_0\times (0,1)$  
and for any base point $m$ (of a component) of $M_0$  the arc $m\times [0,1]$ is disjoint from
$B$. Deforming if necessary $g,g''$ on $\Int(W) \backslash B$, we can assume that    
$g(m\times [0,1])=g''(m\times [0,1])=x$ for any base point $m$   of $M_0$.  By the argument
given at the beginning of Section 2.1, $\tau  (M_0\times [0,1],g)=\tau  (M_0\times [0,1],g'')$. 
  The cobordism $W$ is obtained by
gluing   $M_0\times [0,1]$ to   $V=  {W\backslash  (M_0\times [0,1))}$ along $M_0\times 1$.
 By   (1.2.6), (1.2.8),   $$\tau  (W,g')=
\tau  ( W,g'')=\tau  ( V,g'') \,\tau  (M_0\times [0,1],g'')
$$
$$=\tau  ( V,g) \,\tau  (M_0\times [0,1],g)=\tau  ( W,g).$$

\skipaline \noindent {\bf   2.2.  Functoriality.}
The HQFT's can be pulled back along the maps between the target spaces.
Having   a map $f:(X',x')\to (X,x)$  of path-connected pointed
spaces we can transform any HQFT $(A,\tau)$ with target $X$ into an HQFT with target $X'$. It
suffices to compose the characteristic maps    with $f$ and to apply $(A,\tau)$.
This induces a functor $f^*:Q_{d+1} (X,x) \to Q_{d+1} (X',x')$ (cf. Section 1.2). 

\skipaline \noindent {\bf    2.2.1.   Theorem. }   {\sl   Let $d\geq 1$ and 
$X,X'$ be path-connected spaces with base points $x\in X, x'\in X'$. 
If a map $f:(X',x')\to (X,x)$ induces an isomorphism $\pi_i(X',x') \to \pi_i(X,x)$ for all 
$i\leq d$, then the functor $f^*:Q_{d+1} (X,x) \to Q_{d+1} (X',x')$
is an equivalence of categories.}

\skipaline  {\sl  Proof.}  The standard obstruction theory shows  that for any $d$-dimensional
pointed manifold $M$, the composition with $f$ defines a bijection 
$\Map (M, X') \to \Map (M,X)$. For  a   $(d+1)$-dimensional cobordism $(W,M_0,M_1)$,
the composition with $f$ defines a bijection 
$\Map (W, X')/\sim \, \to \Map (W,X)/\sim  $ where $\sim $ is the equivalence
relation introduced before the statement of Lemma 2.1.1.  This lemma implies that from the
viewpoint of HQFT's there is no difference between $X$-manifolds and $X$-cobordisms on  
one hand and  $X'$-manifolds and $X'$-cobordisms on the other
hand.
In a formal language this means that $f^*:Q_{d+1} (X,x) \to Q_{d+1} (X',x')$ is
an equivalence of categories.

\skipaline \noindent {\bf    2.2.2.  Corollary. }   {\sl   A homotopy equivalence  
$(X',x')\to (X,x)$ of path-connected pointed
spaces induces an equivalence of categories $Q_{d+1} (X,x) \to Q_{d+1} (X',x')$ for all $d$.}

\skipaline \noindent {\bf    2.2.3.  Corollary. }   {\sl   For any 
   connected CW-complex $X$ with base point  $x\in X $, the categories  
$Q_{2} (X,x)$ and  $Q_{2} (K(\pi_1(X,x),1))$ are  equivalent.}

\skipaline The  equivalence is induced by the natural map $X\to K(\pi_1(X,x),1)$
inducing the identity   of the fundamental groups.

\skipaline \noindent {\bf   2.3.  Independence of the base point.}  We show in this section that
the notion of an HQFT with path-connected target $X$ is essentially independent of the choice of a
base point  
$x\in X$. To stress the role of the base point, we  use here the terms   $(X,x)$-manifolds  
and $(X,x)$-cobordisms for $X$-manifolds and $X$-cobordisms.

Let $\alpha:[0,1]\to X$ be a path in $X$ connecting the points $x=\alpha(0), y=\alpha(1)$.
For every  $(X,x)$-manifold $M$ with characteristic map $g:M\to X$, consider
a map $F_{\alpha}:M\times [0,1] \to X$ such that
 
$(\ast)$ $F_{\alpha}\vert_{M\times 0}= g $ and $F_{\alpha}(m\times t)=\alpha(t)$ for all the base points
$m$ of the components of $ M$ and all $t\in [0,1]$. 

The existence of $F_{\alpha}$ follows from the fact that 
$(M\times 0)\cup (\cup_m m\times [0,1])$ is a strong deformation retract  
of $M\times [0,1]$. Any two   maps   satisfying  
$(\ast)$ are homotopic in the class of maps satisfying 
$(\ast)$. Therefore, restricting $F_{\alpha}$ to the base $M\times 1$ we obtain a well-defined  $(X,y)$-manifold.
It is  denoted
$M^{\alpha}$. The same construction applies to  $(X,x)$-cobordisms and   transforms an  
$(X,x)$-cobordism $W$ into an $(X,y)$-cobordism $W^{\alpha}$. Now, every HQFT $(A,\tau)$
with target $(X,y)$ gives rise to an HQFT $({}^{\alpha}\! A, {}^{\alpha}\!\tau)$
with target $(X,x)$   by $({}^{\alpha}\! A)_M=A_{M^{\alpha}}$   and $({}^{\alpha}\! \tau)
(W)=\tau({W^{\alpha}})$.  The action of homeomorphisms is defined by a similar formula.
The notation is chosen so that $({}^{\alpha\beta}\! A, {}^{\alpha\beta}\!\tau)=({}^{\alpha}\! ({}^{\beta}\!A), {}^{\alpha}\!
({}^{\beta}\!\tau))$ for   paths 
$\alpha,\beta:[0,1]\to X$ with $\alpha(1)=\beta (0)$. Note also that if two paths $\alpha,\alpha'$ are
homotopic rel $\{0,1\}$, then   $({}^{\alpha}\! A, {}^{\alpha}\!\tau)=({}^{\alpha'}\! A,
{}^{\alpha'}\!\tau)$. Thus,  we can transport HQFT's along  any path in the target space.
This implies the claim at the beginning of this subsection. 

These constructions may be applied to  $y=x$. This gives a left action of  
$\pi_1(X,x)$ on   HQFT's with target $(X,x)$.  Observe that   for any HQFT $(A,\tau)$
with target $(X,x)$ and any $\alpha\in \pi_1(X,x)$,  the HQFT 
$({}^{\alpha}\! A, {}^{\alpha}\!\tau)$ is   isomorphic to $(A,\tau)$. The isomorphism
is given by the   $K$-isomorphisms
$\{\tau(M\times [0,1], F_\alpha): A_M \to  ({}^{\alpha}\! A)_M\}_M$
where $M$ runs over $(X,x)$-manifolds.

\skipaline \centerline  {\bf Part II.  Crossed group-algebras and $(1+1)$-dimensional   HQFT's}

\skipaline \centerline  {\bf 3.  Crossed group-algebras}

 \skipaline \noindent {\bf   3.1.    Group-algebras.}   Let $\pi$ be a group. 
A {\it $\pi$-graded algebra}  or, briefly, a {\it  $\pi$-algebra} over the ring  $K$ is an
associative algebra   $L$ over   $K$
endowed with a splitting     
$L=\bigoplus_{\alpha\in \pi} L_{\alpha}$ such that 
each $L_{\alpha}$ is a projective $K$-module  of finite type, $L_{\alpha}   L_{\beta} \subset
L_{\alpha \beta}$ for any $\alpha, \beta\in \pi$, and
 $L$ has a (right and left) unit $1_L\in L_1$ where     $1$ is the  neutral element of
$\pi$.    

An example of a $\pi$-algebra is provided by the group ring $L=K[\pi]$ with
$L_\alpha=K\alpha$ for all $\alpha\in \pi$.  More generally, for any associative unital
$K$-algebra $A$ we have  a $\pi$-algebra   $L=A[\pi]$ with
$L_\alpha=A\alpha$ for   $\alpha\in \pi$.  Multiplication in $A[\pi]$ is given by
$(a\alpha) (b\beta)= (ab) (\alpha\beta)$ where $a,b\in A$ and $ \alpha,\beta\in \pi$.

We describe here a few simple operations on $\pi$-algebras.
The direct sum $L\oplus L'$ of two  $\pi$-algebras 
$L,L'$ is   a   $\pi$-algebra
defined by $(L\oplus L')_{\alpha}= L_{\alpha}\oplus  L'_{\alpha}$ for $\alpha\in \pi$.
The tensor product $L\otimes L'$ of    $\pi$-algebras 
$L,L' $ is   a $\pi$-algebra
defined by $(L\otimes L')_{\alpha}= L_{\alpha}\otimes  L'_{\alpha}$.
Multiplication in  
$L\oplus L'$ and  $L\otimes L'$ is induced by multiplication in $L,L'$ in the obvious way.
The dual  $L^*$ of a   $\pi$-algebra $L$ is defined 
as the same module  $L$ with   opposite multiplication $a\circ b=ba$ and the splitting 
$L^* =\bigoplus_{\alpha\in \pi} L^*_{\alpha}$ given by
$  L^*_{\alpha}  =  L_{\alpha^{-1}}$. 

 The group-algebras
can be pulled back and pushed forward along group homomorphisms. Given   a group
homomorphism  $q:\pi'\to \pi$ we can transform any   $\pi$-algebra
$L$ into  a  
$\pi'$-algebra
$ q^*(L)$ defined by $(q^*(L))_{\alpha}= L_{q(\alpha)}$ for any $\alpha\in \pi'$. Multiplication in
$q^*(L)$ is induced by multiplication in $L$ in the obvious way. If the kernel of $q$ is finite then 
we can transform any   $\pi'$-algebra
$L'$ into  a  
$\pi$-algebra
$ q_*(L')$. For    $\alpha\in \pi$, set 
$$(q_*(L'))_\alpha=\bigoplus_{u\in q^{-1} (\alpha)}   L'_u.$$
Multiplication in
$q_*(L')$ is induced by multiplication in $L'$.    Note that 
$q_*(L')=L'$ as algebras:
$$q_*(L')=\bigoplus_{\alpha\in \pi} (q_*(L'))_\alpha=\bigoplus_{u\in \pi'}  L'_u=L'.$$

A {\it Frobenius $\pi$-algebra} is a 
$\pi$-algebra $L$ endowed with  a symmetric $K$-bilinear   form  (inner product) $\eta : L
\otimes L \to K $ such that 

 (3.1.1) $\eta(L_{\alpha}\otimes L_{\beta})=0$ if $\alpha \beta\neq 1$ and     the
restriction of  $\eta$ to $L_{\alpha}\otimes L_{\alpha^{-1}}$ is     non-degenerate for
all $\alpha \in \pi$; 

(3.1.2) $\eta(ab,c)=\eta(a,bc)$ for any $a,b,c\in L$.

Recall that  for   $K$-modules $P,Q$, a  bilinear form $P\otimes Q\to K$ is  
{\it non-degenerate} if both adjoint homomorphisms $P\to Q^*=\Hom(Q,K)$ and $Q\to P^*$ are
isomorphisms.   The group ring $L=K[\pi]$ is a   Frobenius $\pi$-algebra with inner product 
determined by $\eta(\alpha,\beta)=1$   if $\alpha \beta= 1$
and $\eta(\alpha,\beta)=0$   if $\alpha \beta\neq 1$ where  $\alpha,\beta\in \pi$.
 
It is clear that the   direct sum and the tensor product  of Frobenius  $\pi$-algebras 
$L,L'$ are  Frobenius $\pi$-algebras; the inner products in $L,L'$ extend to $L\oplus L'$ 
(resp. to $L\otimes L'$)  by linearity (resp. by multiplicativity).  
The inner product in $L$ induces an inner product in the dual $\pi$-algebra  $L^*$
via the equality of modules $L^*=L$; this makes $L^*$ a Frobenius $\pi$-algebra.
 The pull-backs and 
push-forwards of Frobenius  group-algebras are  Frobenius 
group-algebras in a natural way.

 \skipaline \noindent {\bf   3.2.  Crossed $\pi$-algebras.}   
A {\it crossed $\pi$-algebra} over    $K$ 
is a Frobenius $\pi$-algebra over    $K$  endowed with
  a group homomorphism $ \varphi : \pi \to \Aut  (L)$ satisfying the
following four axioms:

(3.2.1) for all $\beta\in \pi$, $\varphi_{\beta}=\varphi (\beta)$ is an algebra automorphism of
$L$ preserving  $\eta$ and  such that $\varphi_{\beta} (L_{\alpha})\subset
L_{\beta\alpha\beta^{-1}}$ for all $\alpha \in \pi$;

(3.2.2) $\varphi_{\beta}\vert_{L_{\beta}}=\id$, for all $\beta\in \pi$;

(3.2.3) for any $a\in L_{\alpha}, b\in L_{\beta}$, we have  
$\varphi_{\beta}(a)b=ba$;

(3.2.4) for any $\alpha, \beta \in \pi$ and any $c\in   L_{\alpha\beta \alpha^{-1}
\beta^{-1}}$ we have
$$\Tr\, (c\, \varphi_{\beta}:L_{\alpha}\to L_{\alpha})=
\Tr\, ( \varphi_{\alpha^{-1}} c:L_{\beta}\to L_{\beta}).\leqno (3.2.a)$$
Here $\Tr$ is the    $K$-valued trace of endomorphisms of projective $K$-modules of
finite type, see for instance [Tu, Appendix I].  If $K$ is a field then $\Tr$ is the standard trace
of matrices.
   The homomorphism  on the left-hand side of (3.2.a) sends any
$a\in L_{\alpha}$ into $c\, \varphi_{\beta}(a)\in L_{\alpha}$
and the homomorphism  on the right-hand side sends
any
$b\in L_{\beta}$ into $\varphi_{\alpha^{-1}} (cb) \in L_{\beta}$.

 Note a few corollaries of the definition. The module $L_1\subset L$ is
a commutative associative  $K$-algebra with unit. (The commutativity follows from (3.2.3)
 since $\varphi_{1}=\id$.) The restriction of $\eta$ to $L_1$ is
non-degenerate so that the pair $(L_1,\eta)$ is a commutative Frobenius algebra over $K$. 
  The
group $\pi$ acts on $L_1$ by algebra automorphisms preserving $\eta$.

The algebra $L_1$ acts on each $ L_{\alpha}$ by multiplications on the left and on the
right. Axiom (3.2.3)  with $\beta=1$ implies that left and right multiplications by
elements of $L_1$  coincide.

Axiom (3.2.4)  with $\beta=1$  and 
$c=1_L\in L_1$  implies that  for any $\alpha  \in \pi$,
$$ \Dim\,  L_{\alpha} =  \Tr\, (  \id:L_{\alpha}\to L_{\alpha})
=\Tr\, (  \varphi_{1}:L_{\alpha}\to L_{\alpha})
=\Tr\, (  \varphi_{\alpha^{-1}} :L_{1}\to L_{1})\in K.$$
Note that if $K$ is a field then $\Dim\,  L_{\alpha} =(\dim L_\alpha) 1_K$
where  
$\dim $ is the standard integer-valued dimension of vector spaces over $K$
and $1_K\in K $ is the unit of $K$.
In particular if $K$ is a field of   characteristic  0 then   the dimensions of all
$\{L_\alpha\}$ are determined by the character of the representation 
$ \varphi\vert_{L_1}: \pi \to \Aut  (L_1)$.
 
We   define a category, $Q_{2}(\pi)=Q_{2}(\pi;K)$,   whose objects are crossed $\pi$-algebras
over $K$.
  A morphism $L\to L'$  in this category  is an algebra
homomorphism  $L\to L'$ mapping each $ L_{\alpha}$  to  $L'_{\alpha}$,
preserving the unit and the inner product and commuting with the action of $\pi$. 
  It is easy to show  (we shall not use it) that all morphisms in the category  
$Q_{2}(\pi)$ are isomorphisms.

The operations on   group-algebras discussed in Section 3.1 apply also to 
crossed group-algebras. 
The   direct sum and the tensor product  of crossed  $\pi$-algebras 
$L,L'$ are  crossed $\pi$-algebras: the   action of $\pi$ on $L,L'$ extends  to
$L\oplus L'$  (resp. to $L\otimes L'$)  by linearity (resp. by multiplicativity).  
The    action of $\pi$ on $L$ induces an action of $\pi$ on the  dual 
$\pi$-algebra  $L^*$ via the equality   $L^*=L$; this makes $L^*$ a crossed
$\pi$-algebra.

  Given   a group
homomorphism  $q:\pi'\to \pi$ and a  crossed $\pi$-algebra
$L$ we can provide the 
$\pi'$-algebra
$L'= q^*(L)$ with the structure of a crossed $\pi'$-algebra.
The inner
product in $L$ induces an inner product in $L'$ in the obvious way. The action of $\pi$
on $L$ induces an action of $\pi'$ on $L'$ by  
$$\varphi_{\beta}=\varphi_{q(\beta)}: L'_{\alpha}=L_{q(\alpha)}\to 
 L_{q(\beta\alpha\beta^{-1})} =
L'_{\beta\alpha\beta^{-1}} $$ for $\alpha, \beta\in \pi'$.
A similar push-forward construction for crossed group-algebras will be discussed  in
Section 10.3.

 A useful operation on a crossed $\pi$-algebra $(L,\eta, \varphi)$  consists in
rescaling the inner product: for  $k\in K^*$, the triple
$(L, k\eta, \varphi)$ is also a crossed $\pi$-algebra.

Crossed  algebras over   the trivial group   $\pi=1$ are nothing but  commutative  Frobenius
algebras over $K$ whose underlying $K$-modules are projective of finite type. Each
such  Frobenius algebra
$A$  determines a crossed
$\pi$-algebra
$A[\pi]$ over any group $\pi$: it suffices to take the pull-back of $A$ along the trivial
homomorphism $\pi \to \{1\}$. Clearly, $A[\pi]_\alpha=A\alpha$ for all $\alpha\in \pi$. 
The underlying $\pi$-algebra of $A[\pi]$ is the one described at the beginning of  Section 3.1.
The group $\pi$ acts on $A[\pi]$ by permutations of the copies of $A$. 

In the remaining part of Section 3 we   describe several   
constructions of crossed group-algebras.  Our principal result is a classification of
semisimple crossed group-algebras in terms of 2-dimensional group cohomology.

\skipaline \noindent {\bf   3.3.  Example: crossed $\pi$-algebras from 2-cocycles.} Let 
$\{\theta_{\alpha,\beta} \in K^*\}_{\alpha,\beta\in \pi}$
 be a normalized 2-cocycle of the group $\pi$ with values in the
multiplicative group $K^*$. Thus
$$ \theta_{\alpha,\beta}\,
\theta_{\alpha\beta,\gamma}=\theta_{\alpha,\beta\gamma}\, \theta_{ \beta,\gamma}
\leqno (3.3.a)$$
for any $\alpha,\beta, \gamma\in \pi$ and  $\theta_{1,1}=1$.
 We
define  a crossed $\pi$-algebra $L=L^{\theta}$ as follows.

For $\alpha \in \pi$, let $L_{\alpha}$   be the  free  $K$-module  of rank one 
  generated by a vector $l_{\alpha}$, i.e., $L_{\alpha}= K l_{\alpha}$.
Multiplication is defined by $ l_{\alpha} l_{\beta}= \theta_{\alpha,\beta}
l_{\alpha\beta}$. 
Note that multiplication induces an isomorphism  $L_{\alpha} \otimes_K
L_{\beta} \to L_{\alpha\beta}$.
The inner product 
$\eta$ on $L$ is determined by $\eta(l_{\alpha}, l_{\alpha^{-1}})=
\theta_{\alpha,\alpha^{-1}}$ for all $\alpha$ and  $\eta(l_{\alpha}, l_{\beta})=0$ for
$\beta \neq \alpha^{-1}$.  The value of the endomorphism $\varphi_{\beta}:L\to L$ on
each $l_{\alpha}$ is uniquely determined by the
condition   $\varphi_{\beta}(l_{\alpha}) l_{\beta}=l_{\beta} l_{\alpha}$. 

Let us verify the axioms of a crossed $\pi$-algebra.  The associativity of
multiplication   follows from (3.3.a).
Substituting $\beta=\gamma=1$ in (3.3.a)
we obtain that $\theta_{\alpha,1}=1$ for all $\alpha\in \pi$. Substituting 
$\alpha=\beta=1$ in (3.3.a), we obtain that $\theta_{1,\gamma}=1$ for all $\gamma
\in \pi$. This implies that $1_L=l_1\in L_1$ is both right and left unit   of  
 $L$. 

Substituting $\beta=\alpha^{-1}$ and $\gamma=\alpha$ in (3.3.a)
we obtain that $\theta_{\alpha,\alpha^{-1}}=\theta_{\alpha^{-1}, \alpha}$ for all
$\alpha\in \pi$. Hence the form $\eta$ is symmetric. Axiom (3.1.1)  follows
from definitions. A  direct computation shows that $\eta (l_{\alpha},
l_{\beta})=\eta (l_{\alpha}  l_{\beta}, 1_L)$ for all $\alpha,\beta$.
Therefore $\eta(a,b)=\eta(ab,1_L)$ for all $a,b\in L$. This implies   (3.1.2). 

Let us verify (3.2.1). For $\alpha, \alpha', \beta\in \pi$, we have
$$\varphi_{\beta}(l_{\alpha} l_{\alpha'}) l_{\beta}=
\theta_{\alpha,\alpha'}  \varphi_{\beta}(l_{\alpha \alpha'}) l_{\beta}
=\theta_{\alpha,\alpha'}   l_{\beta}  l_{\alpha \alpha'} =
l_{\beta} l_{\alpha} l_{\alpha'}$$
$$
=\varphi_{\beta}(l_{\alpha})  l_{\beta} l_{\alpha'}
=
\varphi_{\beta}(l_{\alpha})  
\varphi_{\beta}(l_{\alpha'}) l_{\beta}.$$
This implies that 
$\varphi_{\beta}(l_{\alpha} l_{\alpha'}) =
\varphi_{\beta}(l_{\alpha})  
\varphi_{\beta}(l_{\alpha'})$.
Substituting $\alpha=\alpha'=1$, we obtain
 that
$\varphi_{\beta}(1_L)=1_L$   and therefore
$\varphi_{\beta}\vert_{L_1}=\id$ for all $\beta\in \pi$. 
For  $a,b\in L$, we have  
$$\eta(\varphi_{\beta}(a), \varphi_{\beta}(b))
=\eta(\varphi_{\beta}(a) \varphi_{\beta}(b), 1_L)=
\eta(\varphi_{\beta}(a b), 1_L).$$
Note that $\oplus_{\alpha\neq 1} L_{\alpha}$ is
orthogonal to
$L_1$ and $\varphi_{\beta}\vert_{L_1}=\id$.
Therefore $\eta(\varphi_{\beta}(a b), 1_L)=\eta( a b, 1_L)=\eta(a,b)$ which proves 
the invariance of $\eta$ under $\varphi_{\beta}$. 

Axioms (3.2.2) and  (3.2.3)    follow  directly from the definition 
of $\varphi_{\beta}$.

Let us check the last axiom (3.2.4). The homomorphism 
$c\, \varphi_{\beta}:L_{\alpha}\to L_{\alpha}$ sends 
$l_{\alpha}$ into $k l_{\alpha}$ with  certain $k\in K^*$. 
The homomorphism 
$ \varphi_{\alpha^{-1}} c:L_{\beta}\to L_{\beta}$ 
sends $l_{\beta}$ into $k' l_{\beta}$ with $k'\in K^*$. 
Note that
$$ k l_{\alpha} l_{\beta} =c\, \varphi_{\beta} (l_{\alpha}) l_{\beta}=c\,
l_{\beta} l_{\alpha} = l_{\alpha} \varphi_{\alpha^{-1}} (c\, l_{\beta})
=  k' l_{\alpha}  l_{\beta}.$$
Therefore $k=k'$ and 
$$\Tr (c\, \varphi_{\beta}:L_{\alpha}\to L_{\alpha})=k=k'=
\Tr ( \varphi_{\alpha^{-1}} c:L_{\beta}\to L_{\beta}).$$

It is easy to see that the isomorphism class of the crossed
$\pi$-algebra $L^{\theta}$ depends only on the  cohomology class  $\theta \in H^2(\pi;
K^*)$ represented by the 2-cocycle $\{\theta_{\alpha,\beta}
\}_{\alpha,\beta\in \pi}$. It is obvious that $L^{\theta\theta'}=L^{\theta}\otimes
L^{\theta'}$
for any $\theta,\theta'\in H^2(\pi; K^*)$ (we use multiplicative notation for the group operation
in $H^2(\pi; K^*)$). The crossed
$\pi$-algebra  corresponding to the neutral element of $  H^2(\pi;
K^*)$ coincides
with the crossed $\pi$-algebra $K[\pi]$ defined at the end of Section 3.2 where we take $A=K$
and  set $\eta (a, b)=ab$ for $a,b\in K$.

\skipaline \noindent {\bf   3.4.  Transfer for crossed algebras.} 
Let   $G\subset \pi$ be a subgroup of $\pi$  of finite index $n=[\pi:G]$.
We show that each crossed $G$-algebra $(L,\eta,\varphi)$  gives rise to
a crossed $\pi$-algebra $(\tilde L, \tilde \eta,\tilde  \varphi)$. It is called the
{\it transfer} of $L$. 

For each right coset class $i\in G\backslash \pi$, fix a representative   
$\omega_i\in   \pi$ so that $i=G\omega_i   \subset \pi$. For  
$\alpha\in \pi$, set $$N(\alpha)= \{i\in G\backslash \pi\, \vert \, \omega_i \alpha \omega_i^{-1} \in
G\}$$ and $$\tilde L_{\alpha}=\bigoplus_{  i  \in N(\alpha)} L_{\omega_i \alpha \omega_i^{-1}}.$$ 
(In particular, if ${\alpha}$ is not conjugated to an element of $G$, then $\tilde L_{\alpha}=0$.)
We provide $\tilde L=\oplus_{\alpha}\tilde L_{\alpha}$ with multiplication as
follows. The multiplication $\tilde L_{\alpha} \otimes \tilde L_{\beta} \to \tilde
L_{\alpha\beta}$ with $\alpha,\beta \in \pi$ sends $L_{\omega_i \alpha \omega_i^{-1}}\otimes
L_{\omega_j \beta \omega_j{-1}}$  into $0$ if
$i\neq j$ and is induced by multiplication   in $L$
$$L_{\omega_i \alpha \omega_i^{-1}}\otimes L_{\omega_i \beta \omega_i^{-1}} \to 
L_{\omega_i \alpha \beta  \omega_i^{-1}}$$  if $i=j\in
N(\alpha)\cap N(\beta)$. Clearly, $\tilde L$ is an
associative algebra.  

By definition, $\tilde L_{1}=\oplus_{  i\in  G\backslash \pi}   L_{1}$
is a direct sum of  $n$ copies of $L_1$. 
The corresponding sum of  $n$   copies of   $1_L\in  L_{1}$ is the unit   $1_{\tilde L}\in
\tilde L_{1}$.  

The inner product $\tilde \eta:\tilde L \otimes \tilde L \to K$ is
determined by 
$$\tilde \eta\vert_{\tilde L_{\alpha}\otimes \tilde L_{\alpha^{-1}}}=
\bigoplus_{i\in N(\alpha)=N(\alpha^{-1})}
\eta\vert_{L_{\omega_i\alpha \omega_i^{-1}} \otimes L_{\omega_i \alpha^{-1}\omega_i^{-1}}} $$
where $\alpha$ runs over $\pi$.
 Clearly,
$\tilde \eta$ is a symmetric bilinear form verifying (3.1.1). It is easy to deduce from
definitions  that $\tilde \eta(a,b)=\tilde \eta(ab,1_{\tilde L})$ for  any $a,b\in 
{\tilde L}$. This
implies   (3.1.2) for $\tilde \eta$.

The action $\tilde \varphi $ of $\pi$ on $\tilde
L$ is defined as follows.  The group $\pi$ acts on
$G\backslash \pi$ by $  {\beta} (i) = i\beta^{-1}$ for  $\beta\in \pi, i \in 
G\backslash \pi$. Then $G \omega_{\beta(i)}=G\omega_i \beta^{-1}$ so that 
$\beta_i=\omega_{\beta(i)} \beta
\omega_i^{-1} \in G$. For any given $\alpha\in \pi$, the map $i\mapsto 
{{\beta(i)}}$
sends  bijectively  $N(\alpha)$ onto $ N(\beta \alpha
\beta^{-1})$. For  every $i\in N(\alpha)$,
 we have the homomorphism $$\varphi_{\beta_i}: L_{\omega_i
\alpha \omega_i^{-1}} \to L_{\omega_{\beta(i)} \beta \alpha \beta^{-1}
(\omega_{\beta(i)})^{-1}}.\leqno  (3.4.a)$$
The direct
sum of these homomorphisms over all $i\in N(\alpha)$ is a homomorphism
 $\tilde \varphi_{\beta}:\tilde L_{\alpha} \to
\tilde L_{\beta \alpha \beta^{-1}}$.
It  extends by additivity to an endomorphism,
$\tilde \varphi_{\beta}$,  of $\tilde L$. 
An easy computation
shows that  
$(\beta \beta')_i={\beta}_{{\beta'}(i)}  \beta'_i$
for any $ \beta, \beta'\in \pi$. This implies that
$\tilde \varphi_{\beta\beta'}=\tilde \varphi_{\beta}\tilde \varphi_{\beta'}$.
It follows from definitions that $\tilde \varphi_1=\id$. Hence 
 $\tilde \varphi$ is an action of $\pi$  on $\tilde L$. 

Note that
the homomorphism
(3.4.a)  preserves multiplication and   inner product in $L$
and therefore $\tilde \varphi_{\beta}$ 
 preserves multiplication and   inner product in $\tilde L$. 
This implies   (3.2.1). 

It is clear that 
for any $i\in N(\beta)$, we have ${\beta(i)}=i$ and 
$ {\beta}_i=\omega_i \beta \omega_i^{-1}$. Therefore, by (3.2.2), the homomorphism 
(3.4.a) is the identity in the case $\alpha=\beta$. This implies   (3.2.2)
for the action $\tilde \varphi_{\beta}$. 

Let us verify axiom (3.2.3) for $\tilde L$. Let 
$a\in L_{\omega_i \alpha \omega_i^{-1}}\subset \tilde L_{\alpha}, b\in L_{\omega_j \beta \omega_j^{-1}}
\subset \tilde L_{\beta}$
 with
$i\in N(\alpha), j\in N(\beta)$.  By definition,
$$\tilde \varphi_{\beta} (a)=\varphi_{{\beta}_i}(a)
\in 
L_{\omega_{{\beta}(i)}  { \beta} \alpha \beta^{-1}
 (\omega_{{ \beta}(i)})^{-1}}.$$
The inclusion $j\in N(\beta)$ implies that 
${ \beta}(j)=j$. 
Therefore, if   $i\neq j$
then ${ \beta}(i)\neq { \beta}(j)=j$. Hence in the case $i\neq j$, we have
$\tilde \varphi_{\beta} (a) b=0=ba$.
Assume that $i=j$. Then ${ \beta}(i)={ \beta}(j)=j$ and  
${ \beta}_i=\omega_i { \beta} \omega_i^{-1}$.
By axioms  (3.2.3)  and (3.2.2) for $L$,
$$\tilde \varphi_{\beta} (a) b= \varphi_{\beta_i} (a) b
=\varphi_{\omega_i { \beta} \omega_i^{-1}} (a) b=\varphi_{\omega_j { \beta} \omega_j^{-1}} (a) b=ba.$$

Let us verify that
for any $\alpha, \beta \in \pi$ and any $c\in   \tilde L_{\alpha\beta \alpha^{-1}
\beta^{-1}}$ we have
$$\Tr (c\, \tilde\varphi_{\beta}:\tilde L_{\alpha}\to \tilde L_{\alpha})=
\Tr ( \tilde\varphi_{\alpha^{-1}} c:\tilde L_{\beta}\to \tilde L_{\beta}). \leqno
(3.4.b)$$   Since both sides of this formula are linear with respect to $c$,
it suffices to consider the case where $$c\in 
L_{\omega_i\alpha\beta \alpha^{-1} \beta^{-1}\omega_i^{-1}} \subset
\tilde L_{\alpha\beta \alpha^{-1}
\beta^{-1}}$$
with $i\in N( {\alpha\beta \alpha^{-1}
\beta^{-1}})$. A direct  application of the definitions shows that both sides of 
  (3.4.b) are equal to $0$ unless $i\in N(\alpha)\cap N(\beta)$. 
If $i\in N(\alpha)\cap N(\beta)$ then the left-hand side of
(3.4.b)  equals the trace of the endomorphism
 $c\, \varphi_{\omega_i\beta \omega_i^{-1}}$ of $\tilde L_{\omega_i\alpha \omega_i^{-1}}$ and 
the right-hand side of
(3.4.b)  equals the trace of the endomorphism
$\varphi_{\omega_i \alpha^{-1}\omega_i^{-1}} c$ of 
$\tilde L_{\omega_i\beta \omega_i^{-1}}$. The equality of these two traces follows from axiom
(3.2.4) for   $L$. 

We leave it as an exercise to the reader to verify that the isomorphism class of the
crossed
$\pi$-algebra $\tilde L$ does not depend on the choice of  the representatives
$\{\omega_i\}$.

\skipaline \noindent {\bf   3.5.  Semisimple crossed $\pi$-algebras. }  
A  crossed $\pi$-algebra  $L=\bigoplus_{\alpha\in \pi} L_{\alpha}$
  is said to be {\it semisimple} if
 the commutative $K$-algebra $L_1$ is semisimple, i.e., if
$L_1$ is a direct sum of
  several copies of the  ring $K$. Note
that direct sums, tensor products,  pull-backs, duals,  and 
transfers of  semisimple crossed  algebras are     semisimple. 
In this subsection we study   the structure 
of semisimple crossed $\pi$-algebras.  

We first briefly discuss semisimple   commutative algebras  of    finite type over $K$. Each
such algebra, $R$,  is a direct sum of  a finite number of
 copies of   $K$, say $\{K_u\}_u$. Let $i_u\in K_u$ be the unit
element of $K_u$. We have $1_R=\sum_u  i_u$ and $i_u i_v =\delta_u^v i_u$ for any $u,v$
where $\delta$ is the Kronecker symbol. 
Each  element $r\in
R$ can be uniquely written  in the form  $r=\sum_u r_u i_u$ with $r_u\in K$. It is
clear that $r$ is an idempotent  (i.e., $r^2=r$) if and only if $r_u \in \{0,1\}$ for all
$u$. We call the set $\{u\,\vert\, r_u \neq 0\}$ the support of $r$. 
Now,  the product of two idempotents   is zero if and only if their supports are disjoint.
This characterizes the set  $\{i_u\}_u  $ as the unique basis of   $R$
consisting of
   mutually annihilating   idempotents. We denote this set 
  by $\mid (R)$; its elements are called {\it  basic idempotents} of $R$.  In
particular, any   algebra automorphism  of $R$ acts by permutations
on   $\mid (R)$.

Let us come back to   crossed algebras.   By a {\it basic idempotent}  of
a semisimple 
crossed $\pi$-algebra $L$ we shall mean a   basic idempotent of    $L_1$. Set  
$\mid(L)=\mid(L_1)\subset L_1$. The equalities $ii'=\delta_{i'}^i i$ for $i, i'\in \mid (L)$
and
$1_L=\sum_{i\in \mid(L)} i$ imply that  
$L$ is a direct  sum of its  subalgebras $iL$  with $i\in \mid(L)$, i.e.,
$L=\oplus_{i\in \mid(L)} iL$ where 
 $$iL=\bigoplus_{\alpha\in
\pi} iL_\alpha. \leqno (3.5.a)$$
The  subalgebras $\{iL\}_{i\in \mid(L)}$ of $L$ are mutually annihilating and 
orthogonal with respect to the inner product.
The action $\varphi$ of $\pi$ permutes these  subalgebras.

We say that  
$L$ is {\it  simple} if it is semisimple and the action 
 of $\pi$ on   $L_1$ is transitive on $\mid(L)$.  Observe that every
semisimple crossed $\pi$-algebra 
$L $
splits uniquely as a direct sum of simple 
crossed $\pi$-algebras. Indeed, the algebra $L_1$ splits as a direct   sum
of its
subalgebras $\{mL_1\}_m $ generated by the orbits $m\subset \mid(L)$ of the $\pi$-action on
$\mid(L)$. 
Then  $L$ is a direct (orthogonal) sum of the  subalgebras $m L =\oplus_{i\in m} iL$
preserved by the action of $\pi$. Clearly, these subalgebras are simple crossed $\pi$-algebras.

For  a simple crossed $\pi$-algebra $L$, any two
elements of $\mid(L)$ can be obtained from each other by the action of $\pi$.
Therefore the value $\eta(i,i)\in K$ with $i\in \mid (L)$   does not depend on the
choice
of $i$. The non-degeneracy of $\eta$ implies that  
$\eta(i,i)\in K^*$. We call $L$   normalized if $\eta(i,i)=1$ for any $i\in \mid (L)$.
It is clear that any simple crossed $\pi$-algebra is obtained from a normalized 
simple crossed $\pi$-algebra by rescaling, see Section 3.2.

\skipaline \noindent {\bf   3.5.1.  Examples. }   1. The
crossed  $\pi$-algebra $L$ associated with a 2-dimensional cohomology class of $\pi$ as in
Section 3.3 is simple and normalized since $L_1=K1_L$ and  $\eta(1_L,1_L)=1$. 

2. Let $G$ be a subgroup of $\pi$ of finite index $n$ and $L$ be a crossed $G$-algebra
associated with    $\theta\in H^2(G;K^*)$. Then the crossed $\pi$-algebra
$\tilde L=L^{\pi, G,\theta}$ obtained as the transfer  of $L$ is   semisimple because the
algebra $\tilde L_{1} $ is a direct sum of  $n$ copies of   $  L_{1}=K$. The action of 
$  \pi$ on 
$\tilde L_{1}=\oplus_{i\in G\backslash \pi}   L_{1} $ permutes the   copies of
$L_1$ via the natural   action of
$\pi$ on $G\backslash \pi$. Hence the action 
 of $\pi$ on      $\mid(\tilde L)=G\backslash \pi$  is   transitive and the crossed $\pi$-algebra
$\tilde L$   is  simple.  
It is easy to check that  $\tilde L$ is normalized.
We have  a distinguished
element  $i_0(\theta)\in \mid(\tilde L)$: this is the unit element of the copy of 
$L_1$ corresponding to $G\backslash G \in G\backslash  \pi$. 

3. The crossed
$\pi$-algebra $A[\pi]$ considered at the end of Section 3.2 is  semisimple if and only if
  $A$ is  semisimple.

 \skipaline

The following theorem yields a   classification of  simple crossed
$\pi$-algebras over fields of   characteristic 0. We introduce the following notation. Denote
by
$C(\pi)$ the set of pairs (a    subgroup   $G\subset \pi$ of finite
index, a  cohomology class  $\theta\in H^2(G;K^*)$).
The group $\pi$ acts on $C(\pi)$ as follows: $\alpha (G,\theta)=
(\alpha G\alpha^{-1}, \alpha_*(\theta))$ where $\alpha\in \pi$
and $\alpha_*:H^2(G;K^*)\to H^2(\alpha G\alpha^{-1};K^*)$ is the isomorphism induced
by the conjugation by $\alpha$. Denote by $D(\pi)$ the set of  
isomorphism classes of  pairs (a normalized 
simple crossed $\pi$-algebra   $(L,\varphi,\eta)$, a    basic idempotent
$i_0\in L_1$).   The group $\pi$ acts on $D(\pi)$ by $\alpha (L,i_0)=
(L, \varphi_{\alpha} (i_0))$ where $\alpha\in \pi$. The set of orbits
$D(\pi)/\pi$ is just the set of  
isomorphism classes of  normalized simple  crossed $\pi$-algebras. 

\skipaline \noindent {\bf    3.6.   Theorem. }   {\sl 
Let the ground ring $K$ be a field of characteristic 0.
 Then the  formula 
 $(G,\theta)\mapsto  (L^{\pi, G,\theta}, 
i_0(\theta))  $ defines a  
$\pi$-equivariant bijection $C(\pi)\to D(\pi)$.
Hence $$D(\pi)/\pi=C(\pi)/\pi.$$}

\skipaline  {\sl  Proof.} We define the inverse mapping 
$ D(\pi)\to  C(\pi) $ as follows.
Let $(L,\varphi,\eta)$ be a normalized 
simple crossed $\pi$-algebra   with distinguished    basic idempotent $i_0\in
L_1$. We compute the dimension of $iL_\alpha \subset L_\alpha$ for $i\in \mid (L),
\alpha\in \pi$  by applying (3.2.4) to $\beta=1$ and $c=i\in L_1$. This gives 
$$\dim (iL_\alpha)=\Tr(a\mapsto ia:L_\alpha\to L_\alpha) \leqno
(3.6.a)$$
$$=\Tr(a\mapsto \varphi_{\alpha^{-1}}(i a):L_1\to L_1)=
\cases 1,~ {  {if}}\,\,\, \varphi_\alpha(i)=i, \\
  0,~    otherwise.\endcases
$$ 
Let  
$G= \{\alpha \in \pi\, \vert\, \varphi_\alpha(i_0)=i_0\}$
be the stabilizer of  $i_0$ with respect to the action of $\pi$ on $\mid (L)$.  We have 
$\dim (i_0L_\alpha)=1$ for all $\alpha \in G$.  
For any  $\alpha \in G\backslash \{1\}$, choose a  generator
$s_\alpha\in i_0L_\alpha$. For $\alpha=1$  set 
$s_\alpha=i_0\in i_0 L_1$. For any $\alpha,\beta\in G$, we have  $s_\alpha
s_\beta =\theta_{\alpha,\beta} s_{\alpha  \beta}$ with 
$\theta_{\alpha,\beta} \in K$. We claim that  $\theta_{\alpha,\beta}\in K^*$. Indeed, by
the non-degeneracy of $\eta$ we have $\eta(s_\beta s_{\beta^{-1}},1_L)=\eta( 
s_\beta,s_{\beta^{-1}})\in K^*$. Therefore $s_\beta s_{\beta^{-1}}=ki_0$ with $k\in
K^*$ and   $$\theta_{\alpha,\beta} \theta_{\alpha\beta,{\beta^{-1}}} s_{\alpha } =
\theta_{\alpha,\beta}
s_{\alpha  \beta} s_{\beta^{-1}}= s_\alpha s_\beta s_{\beta^{-1}}=ki_0 s_{\alpha}= k
s_{\alpha}.$$
 Hence, $\theta_{\alpha,\beta} \in K^*$. The associativity of multiplication in $L$
and the choice $s_1=i_0$ imply that $ \{\theta_{\alpha,\beta}\}_{\alpha,\beta}$ is a
 normalized 2-cocycle of $G$ representing a certain cohomology class
$\theta \in H^2(G;K^*)$. 
Under a different choice of the generators $\{s_\alpha\}_{\alpha \in G}$ we  
obtain a cohomological 2-cocycle. Thus, the formula $(L,i_0)\mapsto (G, \theta)$
yields   a well defined mapping  $ D(\pi)\to  C(\pi) $. It follows   from
definitions that this mapping is $\pi$-equivariant. (The key point is that the action of
$\pi$ on $L$ via $\varphi$ preserves multiplication.) 

We can apply the construction of the previous paragraph to the crossed $\pi$-algebra
$\tilde L$ derived as in Sections 3.3,  3.4  from a subgroup of finite index $G\subset
\pi$ and a cohomology class   $ \theta \in H^2(G;K^*)$. 
The distinguished basic
idempotent of $\tilde L$ is   $i_0=G\backslash G\in G\backslash \pi=\mid(\tilde L)$. 
The stabilizer of $i_0$ with respect to the natural action of $\pi$   is   $G$.  As the
representative $\omega_{i_0}\in G$ of $i_0$ (used in Section 3.3) we   take $1\in G$. For
$\alpha\in G$, we take as the generator $s_\alpha$  of  $i_0\tilde L_\alpha=L_\alpha$
the element $l_\alpha$ used in Section 3.3. Now, it is obvious that  the construction of
the previous paragraph applied to the pair $(\tilde L,i_0)$ gives $(G,  \theta )$. 
Thus, the mapping $ D(\pi)\to  C(\pi) $ contructed above is the left inverse of 
the mapping $ C(\pi)\to  D(\pi) $ in the statement of the theorem. 

To accomplish the proof, it suffices  to show that 
the mapping $ D(\pi)\to  C(\pi) $  is injective.
We need to prove that
any   normalized 
simple crossed $\pi$-algebra 
 $(L,\varphi,\eta)$  with distinguished    basic idempotent $i_0\in
L_1$ can be uniquely reconstructed from its subalgebra $i_0L=\oplus_{\alpha\in G} K
s_\alpha$ where $G\subset \pi$ is the stabilizer of $i_0$
and $s_\alpha$ is a generator of $i_0 L_\alpha$. Note first that the form
$\eta$ on $L$ is completely determined by the formulas  $\eta(a,b)=\eta(ab, 1_L)$ and
$\eta (\sum_{i\in \mid (L)} k_i i,1_L) =\sum_{i\in \mid (L)} k_i$ for any $k_i\in K$.  
Since the action of $\pi$ on $\mid (L)$ is transitive, for each $i\in \mid (L)$ there is
an element $\omega_i\in \pi$ such that $\varphi_{\omega_i}(i_0)=i$. We take
$\omega_{i_0}=1$. The homomorphism $\varphi_{\omega_i}:L\to L$ maps $i_0 L$ isomorphically
onto $iL$.  Therefore the elements $\varphi_{\omega_i} (s_{\alpha})$ with $i\in \mid(L)$
and $\alpha \in G$ form an additive basis of $L$. The product of
two basis elements is computed by
$$\varphi_{\omega_i} (s_{\alpha}) \,\varphi_{\omega_j} (s_{\beta})=
\cases \varphi_{\omega_i} (s_{\alpha} s_{\beta}),~ {  {if}}\,\,\,i=j, \\
  0,~ {  {if}}\,\,\,i \neq j.\endcases$$
It remains to recover 
the action $\varphi$ of $\pi$ on $L$.
For $\alpha\in G$, the homomorphism
$\varphi_{\alpha}:i_0L\to i_0L$  is uniquely determined by the
condition   $\varphi_{\alpha}(s_{\beta}) s_{\alpha}=s_{\alpha} s_{\beta}$
for all $\beta\in G$. 
Each     $\beta\in \pi$ splits as a product  $\omega_i\alpha$ with $\alpha\in G$. 
This gives $\varphi_{\beta}=\varphi_{\omega_i}\varphi_{\alpha}$ and computes the
restriction  of $\varphi_{\beta}$ to $i_0L$.  Knowing these restrictions for all
$\beta\in \pi$, we can uniquely recover the whole action of $\pi$ on $L$
because each   basis
element  of $L$  given above has the form $\varphi_{\alpha}(s)$ with $s\in i_0L$.

\skipaline \noindent {\bf    3.7.  Corollary. }   {\sl   
Let $K$ be a field of characteristic 0.
Then there is a bijection from $H^2(\pi;K^*)$ onto the set of isomorphism classes of 
crossed $\pi$-algebras $L$ such that $\dim L_1=1$ and $\eta(1_L,1_L)=1_K$.}

\skipaline \centerline  {\bf 4.  Two-dimensional HQFT's} 

\skipaline
We shall relate  $(1+1)$-dimensional HQFT's to crossed
$\pi$-algebras where $\pi$ is the fundamental group of the target space. In view of Corollary
2.2.3, we can restrict ourselves to HQFT's with target 
$K(\pi,1)$.

\skipaline \noindent {\bf   4.1.  Theorem. }   {\sl   Let $\pi$ be a group.
Every $(1+1)$-dimensional HQFT  with target   $K(\pi,1)$  determines an \lq\lq underlying"
crossed
$\pi$-algebra. This establishes  
 an equivalence   between   the category $Q_2(K(\pi,1))$ of 
$(1+1)$-dimensional HQFT's with target   $K(\pi,1)$ and   the category $Q_2( \pi )$ of  crossed
$\pi$-algebras. The  direct sums, tensor products,  pull-backs, duals and 
transfers for   $(1+1)$-dimensional  HQFT's correspond 
to  direct sums, tensor products,  pull-backs, duals and 
transfers for  crossed  algebras. }
\skipaline    

It is understood that the pull-backs of HQFT's are effected along maps 
between Eilenberg-MacLane spaces of type $  K(\pi,1)$. By transfers of HQFT's we mean the
transfers described in Section 1.6.  

\skipaline \noindent {\bf   4.2.  Corollary. }   {\sl    
There is a bijective correspondence between   the isomorphism classes of 
$(1+1)$-dimensional HQFT's with target   $K(\pi,1)$ and  the isomorphism classes of crossed
$\pi$-algebras.   }
\skipaline    

Note also another corollary of Theorem 4.1: the group of endomorphisms of 
a $(1+1)$-dimensional HQFT  with target   $K(\pi,1)$ is isomorphic to the group of endomorphisms
of the underlying crossed
$\pi$-algebra.

We complement Theorem  4.1 with
a description of crossed $\pi$-algebras corresponding to   semi-cohomological
HQFT's (cf. Section  1.6).

\skipaline \noindent {\bf    4.3.  Theorem. }   {\sl    
The  underlying 
crossed
$\pi$-algebra of  a cohomological (resp. semi-cohomological) $(1+1)$-dimensional   HQFT  
  with target   $K(\pi,1)$  
   is normalized and simple (resp. semisimple).  The  converse is also true provided   the ground
ring  $K$ is a field of
characteristic 0.  } \skipaline   

Since the splitting of a semisimple crossed group-algebra as a direct sum of simple
crossed group-algebras is unique (up to permutation of   summands), we obtain the following
corollary. 

\skipaline \noindent {\bf    4.4.  Corollary. }   {\sl   
Let the ground ring $K$ be a field of characteristic 0.   
   The 
  splittting of a semi-cohomological $(1+1)$-dimensional  HQFT over $K$ into a direct sum of
rescaled cohomological HQFT's is unique.} \skipaline

 In the remaining part of Section 4 we
construct  the underlying  
crossed algebra of a   $(1+1)$-dimensional  HQFT and show the functoriality of this construction.  The proof
of  Theorems 4.1 and 4.3 will be  given in Section  5.

Throughout this section we fix a group $\pi$ and a  $(1+1)$-dimensional HQFT
$(A,\tau)$  with target   $X=K(\pi,1)$ and base point $x\in X$.

\skipaline \noindent {\bf   4.5.  Preliminaries on $(1+1)$-dimensional HQFT's.}  Here we  
reformulate the data provided by  a $(1+1)$-dimensional HQFT  $(A,\tau)$
in a form convenient for the sequel.  Observe first that a
1-dimensional connected  $X$-manifold $M$ is just a  pointed oriented circle  endowed
with a map into $X$ sending the base point  into $x$. This is nothing but a
  loop in $X$ with endpoints in $x$.  The $K$-module  $A_M$   
depends only on the  class of this loop in $\pi_1(X,x)=\pi$, see Section 2.1. In this way,
for each $\alpha\in \pi$ we obtain  a $K$-module  denoted
$L_{\alpha}$.  
A non-connected  $X$-manifold $M$ is a finite
non-ordered family 
$\{\alpha_i\}_i$ of  loops in $X$  and   $A_M=\otimes_i L_{\alpha_i}$.
If $M=\emptyset$, then $A_M=K$.

Let   $W$ be a compact oriented surface 
 with pointed oriented boundary   endowed with a map
$g:W\to X$ sending the base points of all the components of $\partial W$ into $x$.
We   write $C_{+}$ for  a  component $C \subset \partial W$
with orientation induced from 
$ W$ and we write $C_-$ for $C$ with opposite orientation.
Let $\{(C_-^p,\alpha_p \in \pi)\}_p$ be the components of $\partial W$  whose orientations are 
opposite to the one induced
from $ W$. Here $\alpha_p$ is  represented by 
the restriction of $g$ to $C_-^p$. 
Let $\{(C_+^q,\beta_q \in \pi)\}_q$ be the components of $\partial W$  whose orientations are 
  induced
from $ W$. Here $\beta_q$ is   represented by 
the restriction of $g$ to $C_+^q$.  
 We   view $W$ 
as an $X$-cobordism between $\cup_p (C_-^p,\alpha_p)$ and
 $\cup_q (C_+^q,\beta_q)$.  By the remarks of  Section 2.1, the HQFT $(A,\tau)$ gives
rise to a homomorphism  $$\tau(W)\in \Hom_K (\bigotimes_{p} L_{\alpha_p}, \bigotimes_{q}
L_{\beta_q})= (\bigotimes_{p} L_{\alpha_p}^*) \otimes (\bigotimes_{q}
L_{\beta_q}).$$
This homomorphism is preserved  under any homotopy of the map $g:W\to X$ relative to
the base points on $\partial W$.
The axioms (1.2.6) and  (1.2.7) tell  us that the homomorphism $\tau(W)$ is 
multiplicative under the gluing of cobordisms and that $\tau (W)=\id$
if $W$ is a cylinder as in (1.2.7).

\skipaline \noindent {\bf   4.6.   Annuli  and discs  with   holes.}
In this subsection we discuss the structures
of $X$-cobordisms on   annuli  and discs with 2 holes.

Let $C$ denote the annulus  $S^1\times [0,1]$. We fix an orientation of $C$ once for
all. Set   $C^0=S^1\times 0 \subset \partial C$ and $C^1=S^1\times 1\subset \partial
C$.   Let us provide $C^0,C^1$ with base points $c^0=s \times 0, c^1=s \times 1$,
respectively, where $s\in S^1$.
 For any signs   $\varepsilon, \mu=\pm$ we  denote by $C_{\varepsilon,\mu}$
the triple $(C, C^0_{\varepsilon},C^1_{\mu})$.  This is an  
annulus  with oriented pointed boundary.  By definition,
$$\partial C_{\varepsilon,\mu}=(\varepsilon C^0_{\varepsilon}) \cup 
( \mu C^1_{\mu}).$$
 The
homotopy class of a map  $g:C_{\varepsilon,\mu} \to X$ is  
determined by the homotopy classes $\alpha,\beta \in \pi$ represented by  the loops 
$g\vert_{C^0_{\varepsilon}}$ and  
$g\vert_{s \times [0,1]}$, respectively. Here   the interval $[0,1]$ is oriented from 0
to 1. Note that  the loop $g\vert_{C^1_\mu}$  
represents   $( \beta^{-1} \alpha^{-\varepsilon} \beta)^{\mu}$. 
We denote by 
$C_{\varepsilon,\mu} (\alpha; \beta )$  the annulus  $C_{\varepsilon,\mu}$
endowed with the map to $X$ corresponding to the pair $\alpha,\beta \in \pi$.

 Let $D$ be an oriented  2-disc with two holes. Denote the boundary components of $D$
by $Y,Z,T$  and provide them  with base points $y,z,t$, respectively. 
For any signs $\varepsilon , \mu,\nu=\pm$ we denote by $D_{\varepsilon,\mu,\nu}$
the tuple 
$(D, Y_{\varepsilon}, Z_{\mu},T_{\nu})$. This is a 
2-disc with two holes  with oriented pointed boundary.
By definition,
$$\partial D_{\varepsilon,\mu,\nu}=(\varepsilon Y_{\varepsilon}) \cup ( \mu Z_{\mu})\cup (\nu T_{\nu}).$$
To analyse the homotopy classes of maps 
$D_{\varepsilon,\mu,\nu}\to X$, we fix two proper embedded  arcs
$ty$ and $tz$ in $D$ leading from  $t$ to $y,z$  and  mutually disjont  
  except in the endpoint $t$. 
To every map $g:D_{\varepsilon,\mu,\nu}\to X$ we assign  the homotopy classes of the
 loops
$g\vert_{Y_{\varepsilon}}, g\vert_{Z_{\mu}},g\vert_{ty}, g\vert_{tz}$. This
establishes a bijective correspondence between the set of homotopy classes of maps
$D_{\varepsilon,\mu,\nu}\to X$ and $\pi^4$. 
For any $\alpha,\beta,\rho, \delta\in \pi $ denote by 
$D_{\varepsilon,\mu,\nu} (\alpha, \beta;\rho, 
\delta)$  the disc $D_{\varepsilon,\mu,\nu}$ endowed with the
map to $X$ corresponding to the tuple $\alpha, \beta,\rho, 
\delta$. 
Note that 
the loops $g\vert_{Y_{\varepsilon}}, g\vert_{Z_{\mu}}, g\vert_{T_{\nu}}$  
represent  the classes   $\alpha, \beta,(\rho \alpha^{-\varepsilon} \rho^{-1} 
\delta \beta^{-\mu} \delta^{-1})^{\nu}$, respectively.

\skipaline \noindent {\bf   4.7.  Algebra $L $.}
We provide the direct sum 
$$L=\bigoplus_{\alpha\in \pi} L_{\alpha}$$
with the structure of an associative algebra as follows.
  The disc with two holes  
$D_{--+} (\alpha, \beta;1, 1)$ is an $X$-cobordism between  
$(Y_-,\alpha) \cup (Z_-,\beta)$ and $(T_+,\alpha  \beta)$.  Note that  the  map $D\to X$
in question (considered up to homotopy) sends the intervals $ty,tz$ into the base point
$x\in X$.  The corresponding  homomorphism  $$\tau (D_{--+} (\alpha, \beta;1, 1)): 
L_{\alpha}\otimes L_{\beta}\to  L_{\alpha \beta} $$
   defines a $K$-bilinear multiplication in $L$ by
$$ab=\tau
(D_{--+} (\alpha, \beta;1, 1)) (a\otimes b) \in L_{\alpha \beta} \leqno (4.7.a)$$
for $a\in  L_{\alpha}, b\in  L_{\beta}$.
By a standard argument, this multiplication is associative. The key point is that under the gluing of 
$D_{--+}(\alpha, \beta;1, 1)$ to $D_{--+}(\alpha \beta, \gamma;1, 1)$ along  an
$X$-homeomorphism $(T_+,\alpha\beta)=(Y_-,\alpha\beta)$ we obtain the same
$X$-cobordism   as  under the gluing of 
$D_{--+}(  \beta, \gamma;1, 1)$ to $D_{--+}(\alpha, \beta \gamma;1, 1)$ along  a
homeomorphism $(T_+,\beta \gamma)=(Z_-,\beta \gamma)$.

The unit of  $L$ is constructed as follows. Let $B_+$ be an oriented 2-disc
whose boundary is pointed and   endowed with the orientation induced from 
$B_+$.
There is only one homotopy class of maps  $B_+\to X$.  The 
corresponding  homomorphism $\tau(B_+):K\to L_1$ sends the unit $1\in K$ into an
element of $L_1$, denoted $1_L$.
 This element is a
right unit   in $L$ because the gluing of $B_+$ to $D_{--+} (\alpha, 1;1, 1)$ along  an
$X$-homeomorphism  $\partial B_+=Z_-$ yields  the annulus   
$C_{-+} (\alpha; 1)$  and
axioms (1.2.6, 1.2.7) apply. Similarly, $1_L$ is a left unit of $L$.

\skipaline \noindent {\bf   4.8.  Action of $\pi$ and the form $\eta$.}
The group $\pi$ acts on 
$L$
  as follows. For $\alpha,\beta \in \pi$, the annulus  
$C_{-+} (\alpha; \beta^{-1})$ is an $X$-cobordism  
between   $(C^0_-,\alpha) $ and $(C^1_+,\beta\alpha \beta^{-1})$.   Set
  $$\varphi_{\beta}=\tau (C_{-+} (\alpha;
\beta^{-1})):  L_{\alpha} \to  L_{\beta\alpha \beta^{-1}}. $$
Observe that  the gluing of  
$C_{-+}(\alpha; \beta^{-1})$ to $C_{-+}(\beta\alpha \beta^{-1}; \gamma^{-1})$ with
$\gamma\in \pi$ yields  $C_{-+}(\alpha; (\gamma \beta)^{-1})$. Axiom  (1.2.6)  implies
that $\varphi_{\gamma \beta }= \varphi_{\gamma} \varphi_{\beta} $.  By axiom
(1.2.7),    $\varphi_1=\id$.  

The annulus 
$C_{--} (\alpha; 1)$ is an $X$-cobordism between  
   $(C^0_-,\alpha) \cup (C^1_-,{\alpha^{-1}})$ and $\emptyset$.   Set
  $$\eta \vert_{L_{\alpha} \otimes   L_{\alpha^{-1}}}= \tau (C_{--} (\alpha;
1)):  L_{\alpha} \otimes   L_{\alpha^{-1}} \to K. $$
The direct sum of these pairings over all $\alpha\in \pi$ yields the form
$\eta:L\otimes L \to K$.

\skipaline \noindent {\bf    4.9.   Lemma. }   {\sl   
The algebra $L$ with action  $\varphi$ and form $\eta$
   is a crossed $\pi$-algebra.}

\skipaline  {\sl  Proof.} The existence of a unit was verified above.   Let us verify the other
axioms.

 Verification of (3.1.1).  Note   that by   axiom
(1.2.7),
$\tau (C_{-+} (\alpha; 1))= \id_{L_{\alpha}}$. The annulus 
$C_{-+} (\alpha; 1)$ can be obtained by gluing $C_{--}(\alpha; 1)$ to
$C_{++}(\alpha^{-1};1)$ along a homeomorphism
$(C^1_-,\alpha^{-1})=(C^0_+,\alpha^{-1})$. The multiplicativity of $\tau$  implies 
that  $$(\tau (C_{--}(\alpha; 1)) \otimes \id_{L_{\alpha}}) \circ
(\id_{L_{\alpha}}\otimes \tau (C_{++}(\alpha^{-1};1))) =\tau (C_{-+} (\alpha; 1)) =
\id_{L_{\alpha}}.$$  
If $K$ is a field then presenting the  homomorphisms 
$\tau (C_{--}(\alpha; 1))=
\eta \vert_{L_{\alpha} \otimes   L_{\alpha^{-1}}}$
and
$\tau (C_{++}(\alpha^{-1};1)): K\to  L_{\alpha^{-1}} \otimes  
L_{\alpha}$ by matrices (with  respect
to certain bases of $L_{\alpha},  L_{\alpha^{-1}}$) we easily obtain
  that $\dim L_{\alpha}\leq \dim L_{\alpha^{-1}}$ and the matrix of the
pairing $\eta \vert_{L_{\alpha} \otimes   L_{\alpha^{-1}}}$  admits a right inverse. By
symmetry, $\dim L_{\alpha}= \dim L_{\alpha^{-1}}$  and the
latter pairing is non-degenerate. In the    case of an arbitrary $K$ one should use the argument
given in [Tu, Section III.2].

 Verification of (3.1.2). Let us prove that $\eta(ab,c)=\eta(a,bc)$ for any $a \in
L_{\alpha}, b\in L_{\beta},c \in L_{\gamma}$ where $\alpha, \beta, \gamma\in \pi$.
 If $\alpha \beta \gamma\neq1$, then $\eta(ab,c)=0=\eta(a,bc)$.
Assume that $\alpha \beta \gamma=1$.  Gluing
the annulus  $C_{--} (\alpha\beta;1)$ to   $D_{--+}(\alpha,\beta;1,1)$ 
along an $X$-homeomorphism 
$(C^0_-,\alpha\beta)=(T_+,\alpha\beta)$ we obtain   $D_{---}(\alpha,\beta;1,1)$.
Therefore 
$$\eta(ab,c)= \tau (D_{---}(\alpha,\beta;1 ,1 )) (a\otimes b\otimes c)$$
where 
$$\tau (D_{---}(\alpha,\beta;1 ,1))\in \Hom_K ( L_{\alpha} \otimes
L_{\beta}  \otimes  L_{\gamma}, K).$$ 
Gluing   
  $C_{--} (\alpha;1)$ and    $D_{--+}(\beta,\gamma;1,1)$ 
along an $X$-homeomorphism 
$(C^1_-,\alpha^{-1})=(T_+,\beta\gamma)$ we obtain   $D_{---}(\beta,\gamma;1,1)$.
Hence
$$\eta(a,bc)= \tau (D_{---}(\beta,\gamma;1,1)) ( b\otimes c\otimes a)$$
where 
$$\tau (D_{---}(\beta,\gamma;1,1)) \in \Hom_K (L_{\beta} \otimes
L_{\gamma}  \otimes  L_{\alpha},K).$$ 
It remains to observe that there is an $X$-homeomorphism  
$D_{---}(\alpha,\beta;1,1)\to D_{---}(\beta,\gamma;1,1)$  mapping the boundary
components $Y,Z,T$ of the first disc with holes onto the  boundary   components
$T,Y,Z$  of the second disc with holes, respectively. By axiom (1.2.4),  
$\eta(ab,c)=\eta(a,bc)$.

 Verification of (3.2.1). Let us prove that $\varphi_{\gamma}(ab)=
\varphi_{\gamma}(a) \varphi_{\gamma} (b)$ for
any $a \in L_{\alpha}, b\in L_{\beta}$ where $\alpha, \beta,
\gamma\in \pi$.
   Gluing  
 $C_{-+} ( \alpha \beta ; \gamma^{-1})$ 
  to   $D_{--+}( \alpha, \beta;1,1)$  along an $X$-homeomorphism  
$(C^0_-, \alpha\beta)=(T_+,\alpha\beta)$, we obtain   $D_{--+}(\alpha,\beta; \gamma ,\gamma )$.
Hence
$$\varphi_{\gamma}(ab)= \tau ( D_{--+}(\alpha,\beta; \gamma ,\gamma )) (
a\otimes b)$$ where 
$$\tau ( D_{--+}(\alpha,\beta; \gamma ,\gamma ))
\in \Hom_K (L_{\alpha} \otimes
L_{\beta},   L_{ \gamma\alpha \beta\gamma^{-1}}).$$ 
Similarly,  gluing
   $C_{-+} (\alpha; \gamma^{-1})\cup C_{-+} (\beta; \gamma^{-1})$ to  
$D_{--+}(\gamma\alpha \gamma^{-1}, \gamma\beta\gamma^{-1};1,1)$  along  
$X$-homeomorphisms  $$(C^1_+,\gamma\alpha \gamma^{-1})=(Y_-,\gamma\alpha
\gamma^{-1})\,\,\,\, {\text {and}} \,\,\,\,
 (C^1_+,\gamma\beta\gamma^{-1})=(Z_-,\gamma\beta\gamma^{-1}),$$ respectively, we
obtain   $D_{--+}(\alpha,\beta; \gamma ,\gamma )$. Therefore  $$\varphi_{\gamma}(a)
\,\varphi_{\gamma} (b) = \tau (D_{--+}(\alpha,\beta;\gamma  ,\gamma  )) (a\otimes b)
=\varphi_{\gamma}(ab).$$

The proof of the identity $\eta (\varphi_{\gamma}(a),\varphi_{\gamma}(b))=
\eta (a,b)$ is similar.

  Verification of (3.2.2).   The Dehn twist along the circle  $S^1\times
(1/2)\subset  C_{-+} (\alpha; 1)$  yields
an  $X$-homeomorphism   $C_{-+} (\alpha; 1)\to C_{-+}(\alpha; \alpha)$.  Axiom (1.2.4)
implies that 
 $\tau(C_{-+}(\alpha; \alpha))=\tau(C_{-+}(\alpha;1))=\id$.
 Therefore $\varphi_{\alpha}\vert_{L_{\alpha}}=\id$.

  Verification of  (3.2.3).    Note first that 
 for any  $\rho,\delta\in \pi$ the 
homomorphism
$$\tau (D_{--+}(\alpha,\beta;\rho,\delta)): L_{\alpha}\otimes
L_{\beta}\to  L_{\rho\alpha\rho^{-1} \delta\beta\delta^{-1}}$$ can be   computed in
terms of $\varphi$ and multiplication in $L$: $$\tau
(D_{--+}(\alpha,\beta;\rho,\delta)) (a \otimes b)=\varphi_\rho (a)\,
\varphi_\delta(b).\leqno (4.9.a)$$ This follows from the   multiplicativity of $\tau$
using the  splitting of $D_{--+}(\alpha,\beta;\rho,\delta)$ as a union of
$D_{--+}(\alpha,\beta;1,1)$ with 
 $C_{-+}(\alpha; \rho^{-1})$ and $ C_{-+}(\beta; \delta^{-1})$.

Consider a self-homeomorphism $f$ of the disc with two holes $D $ which is the
identity on $T$ and  which permutes   $(Y,y)$ and $(Z,z)$.
 We choose $f$
so that    $f(tz)= ty$ and   $f(ty)$ is   an embedded
arc leading from $t$ to $z$ and homotopic to the product of four  arcs 
$ty, \partial Y, (ty)^{-1}, tz$.   An easy 
  computation   shows that $f$ is an $X$-homeomorphism
  $D_{--+}( \alpha, \beta;1,1) \to D_{--+}( \beta, \alpha;1, \beta^{-1})$.
Axiom (1.2.4) implies that 
the homomorphisms 
$$\tau (D_{--+}( \alpha, \beta;1,1)): L_{\alpha}\otimes L_{\beta}\to  L_{ \alpha
\beta} $$ and $$\tau (D_{--+}( \beta, \alpha;1,\beta^{-1})) : L_{\beta}\otimes
L_{\alpha}\to  L_{ \alpha \beta} $$ are obtained from  each other  by the permutation
of the   two tensor factors. Therefore for any $a\in L_\alpha, b\in L_\beta$
$$ab=
\tau (D_{--+}( \alpha, \beta;1,1))(a\otimes b)=\tau (D_{--+}( \beta, \alpha;1,\beta^{-1}))
(b\otimes a)= b \,\varphi_{\beta^{-1}}(a).$$
This is equivalent to (3.2.2). 

  Verification of  (3.2.4).  Fix an orientation of $S^1$ and consider the 
2-torus $ S^1\times S^1$ with product orientation. 
Let $B\subset S^1\times S^1$ be a closed 
embedded 2-disc   disjoint from the loops $S^1\times s$ and
$s\times S^1$ where $s\in S^1$. Consider the punctured torus
$H=(S^1\times S^1)\backslash \Int B$ with orientation induced from $S^1\times S^1$.  
Let us provide the boundary circle $\partial H=\partial B$   with orientation
 opposite to the one induced from $H$. 
We choose a base point   on $ \partial H$ and an arc $r\subset H$ joining this point 
  to  $s \times s\in H$. We can assume that $r$ meets
the loops  $S^1\times s$ and
$s\times S^1$ only in its endpoint $s \times s$.

Consider a map $g:H\to X=K(\pi,1)$ such that $g(r)=x\in
X$ and  the   restrictions of $g$ to  $S^1\times s, s\times S^1$  represent  
$\alpha,\beta \in \pi$, respectively.  (The orientations of $S^1\times s, s\times S^1$
are  induced by the one of $S^1$.) Then the loop $g\vert_{\partial H}$
represents $\alpha\beta \alpha^{-1} \beta^{-1}$. Now, the  pair $(H,g)$ is an
$X$-cobordism  between $ (\partial H_-,g\vert_{\partial H})$ and $\emptyset$.  
This gives a homomorphism $\tau (H,g): L_{\alpha\beta \alpha^{-1} \beta^{-1}}\to K$.
  We claim that both sides of formula (3.2.a)
are equal to  
 $\tau (H,g)(c)$.
This   would imply (3.2.4).

 The pair $(H,g)$ can be  obtained from 
  $D_{--+}  (\alpha\beta \alpha^{-1} \beta^{-1},\alpha;1,\beta)$ by gluing the boundary
components $(Z_-,\alpha)$ and $(T_+,\alpha)$ along an $X$-homeomorphism. 
(The circles  $Z_-$ and $T_+$  give the loop $S^1\times s\subset T$.)
A standard
argument in the theory of TQFT's shows that the homomorphism $\tau (H,g):
L_{\alpha\beta \alpha^{-1} \beta^{-1}}\to K$ is the partial trace of the
homomorphism
$$\tau (D_{--+}  (\alpha\beta \alpha^{-1} \beta^{-1},\alpha;1,\beta)):
 L_{\alpha\beta \alpha^{-1} \beta^{-1}}\otimes L_{\alpha}\to  L_{\alpha}.$$
For   $d\in  L_{\alpha}$, we have
$$\tau (D_{--+}  (\alpha\beta \alpha^{-1} \beta^{-1},\alpha;1,\beta))
(c\otimes d)=c\, \varphi_\beta (d).$$
Therefore  $\tau (H,g) (c)=\Tr (c\, \varphi_{\beta}:L_{\alpha}\to L_{\alpha})$.
Similarly, 
the pair $(H,g)$ is obtained from
  $D_{--+}  (\alpha\beta \alpha^{-1} \beta^{-1},\beta;\alpha^{-1},\alpha^{-1})$ by gluing the
boundary components $(Z_-,\beta)$ and $(T_+,\beta)$ along an $X$-homeomorphism.
(The circles  $Z_-$ and $T_+$  give the loop $s\times S^1\subset T$.)
Thus,
$\tau (H,g): L_{\alpha\beta \alpha^{-1} \beta^{-1}}\to K$ is the partial trace of the
homomorphism
$$\tau (D_{--+}   (\alpha\beta \alpha^{-1} \beta^{-1},\beta;\alpha^{-1},\alpha^{-1})):
 L_{\alpha\beta \alpha^{-1} \beta^{-1}}\otimes L_{\beta}\to  L_{\beta}.$$
For   $d\in  L_{\beta}$, we have
$$\tau (D_{--+}   (\alpha\beta \alpha^{-1} \beta^{-1},\beta;\alpha^{-1},\alpha^{-1}))
(c\otimes d)= \varphi_{\alpha^{-1}}(c) \,\varphi_{\alpha^{-1}}
(d)=\varphi_{\alpha^{-1}}(cd).$$ Therefore  $\tau (H,g) (c)=
\Tr ( \varphi_{\alpha^{-1}} c:L_{\beta}\to L_{\beta})$.

 \skipaline \noindent {\bf    4.10.   Functoriality. } The construction of the underlying crossed algebra
is functorial with respect to morphisms of HQFT's as defined in Section 1.2.
Such a morphism
 $\rho:(A,\tau)\to  (A',\tau')$ yields for each $\alpha\in \pi$ a $K$-linear homomorphism 
$\rho_\alpha:L_\alpha\to L'_\alpha$ where $L$ (resp. $L'$) is the underlying crossed
algebra of $(A,\tau)$ (resp. of $(A',\tau')$).  
The commutativity of  
   the natural square
diagrams associated with the cobordisms   
$D_{--+} (\alpha, \beta;1, 1)$, $B_+$, $C_{--} (\alpha;1)$, and 
$C_{-+} (\alpha;\beta^{-1}) $ 
  imply that  $\oplus_{\alpha} \rho_\alpha:L\to L'$
is an algebra homomorphism preserving the unit and  the inner product and commuting with
the action of
$\pi$. This establishes the functoriality of the underlying crossed algebra.

 \skipaline \noindent {\bf    4.11.  The Verlinde formula. } The arguments given at the 
 end of Section 4.9 allow us to compute explicitly the values of $\tau$  on closed oriented
$X$-surfaces. Note first that for any $\alpha,\beta\in \pi$ there is a
unique element $h_{\alpha,\beta} \in L_{\beta \alpha \beta^{-1} \alpha^{-1}}$
such that
$$\eta(c,h_{\alpha,\beta}) =\Tr (c\, \varphi_{\beta}:L_{\alpha}\to L_{\alpha})=
\Tr ( \varphi_{\alpha^{-1}} c:L_{\beta}\to L_{\beta}) $$
for any $c\in L_{ \alpha \beta \alpha^{-1} \beta^{-1}}$.
The element $h_{\alpha,\beta}$ can be geometrically described as follows. Let  $(H',g)$ be the
same punctured
$X$-torus $(H,g)$ as in Section 4.9 with opposite orientation of  $\partial H $  (so that it is 
induced from $H$). The pair $(H',g)$
 is an $X$-cobordism  between  $\emptyset$ and $ (\partial
H',g\vert_{\partial H'})$.   Then $h_{\alpha,\beta}= \tau (H',g)(1_K)$ where  $\tau (H',g):  
 K\to L_{\beta \alpha
\beta^{-1} \alpha^{-1}}$ is the 
homomorphism
associated with $(H',g)$.

To compute $\tau(W)$ for a  closed oriented $X$-surface $(W,g:W\to X)$ of genus $n$, we choose
a   point $w\in W$ and  a   system of generators
$  a_1, a_2,  ..., a_{2n}\in \pi_1(W,w)$ subject to the only relation
$\prod_{r=1}^n  a_{2r-1} a_{2r} a_{2r-1}^{-1} a_{2r}^{-1} =1$. 
We   deform $g$ so that $g(w)=x$. Then $g$ induces a group homomorphism 
$g_{\#}:\pi_1(W,w) \to  \pi=\pi_1(X,x)$. Using a decomposition of $W$ into
$n$ copies of   $H'$, a disc with $n$ holes, and a disc  $B_-$ (cf. Section
5.1 below) one can easily compute that 
$$\tau(W,g)=\eta( \prod_{r=1}^n h_{ g_{\#}(a_{2r}), g_{\#}(a_{2r-1})}, 1_L) \in K. \leqno
(4.11.a)$$ Note that the relation $\prod_{r=1}^n  a_{2r-1} a_{2r} a_{2r-1}^{-1} a_{2r}^{-1} =1$ yields 
$\prod_{r=1}^n h_{ g_{\#}(a_{2r}), g_{\#}(a_{2r-1})} \in L_1$. 

 If the crossed $\pi$-algebra $L$ 
underlying $(A,\tau)$ is
semisimple then we can split each 
$h_{\alpha,\beta} \in L_{\beta \alpha \beta^{-1} \alpha^{-1}}$
 as follows:
$$h_{\alpha,\beta}=\sum_{i\in \mid(L)}
h_{\alpha,\beta, i}$$
where $h_{\alpha,\beta, i}=i h_{\alpha,\beta}\in i L_{\beta \alpha \beta^{-1} \alpha^{-1}}$.
Now we can rewrite
(4.11.a) in the   form 
$$\tau(W,g)=\sum_{i\in \mid(L)} \eta( \prod_{r=1}^n h_{ g_{\#}(a_{2r}), g_{\#}(a_{2r-1}), i}, i).
\leqno (4.11.b)$$
This formula   generalizes
the well known Verlinde formula corresponding to $\pi=1$:  in this case
$h_{1,1, i}=k_i i$ with $k_i\in K$ and   formula (4.11.b) has the standard form 
$\tau(W)=\sum_{i\in \mid(L)} \eta( i,i)  (k_i)^n$.  

If $K$ is a field of characteristic 0, then  by   (3.6.a) the  
  $i$-th summand on the righ-hand side of (4.11.b) is 0 unless
$g_{\#}(a_{2r-1} a_{2r} a_{2r-1}^{-1} a_{2r}^{-1})(i)=i$ for all $r=1,...,n$.

Formula  (4.11.a)  extends to $X$-surfaces with boundary.
Consider for simplicity 
  an   $X$-cobordism  $W$ with  bottom base $\cup_{p=1}^{m} (C_-^p,\alpha_p \in \pi)$
as in Section 4.5 and
  empty top base. The homomorphism
 $\tau(W): \bigotimes_{p=1}^{m} L_{\alpha_p}\to K$  is computed as follows.
Choose a   point $w\in
\Int W$ and deform $g$ relative to $\partial W$ so that $g(w)=x$.  
 Now, 
connect $w$ by disjoint (except in $w$) embedded arcs,   $d_p$,  to the base points of the
circles 
$\{C_-^p\}$.  Denote by $c_p$ the element of $\pi_1 (W,w)$ represented by the loop  
obtained as the product of    $d_p$, the loop parametrizing $C_-^p$ and the path
inverse to $d_p$. 
The group $\pi_1 (W,w)$ can be presented by $2n+m$ generators
$  a_1, a_2, ..., a_{2n}, c_1,...,c_m$
and one relation
$$\prod_{r=1}^n a_{2r-1} a_{2r} a_{2r-1}^{-1} a_{2r}^{-1}\, \prod_{p=1}^{m} c_p=1.$$
The  paths  $d_p$  are mapped by $g$ into loops in $(X,x)$  representing
certain   $\{\gamma_p\in \pi\}_p  $.
Then for any $\{u_p \in  L_{\alpha_p}\}_p$
we have
$$ \tau(W)( \otimes_{p=1}^{m}  {u_p})=\eta(\prod_{r=1}^n h_{ g_{\#}(a_{2r}), g_{\#}(a_{2r-1})}
\prod_{p=1}^{m} \varphi_{ \gamma_p} (u_p),1_L).$$
 
 \skipaline \noindent {\bf    4.12.   Remark. } We may consider   generalized 
 $(1+1)$-dimensional
  HQFT's     such that the   homomorphism  $\tau$  is associated  only with
surfaces of genus 0. These  HQFT's correspond  to generalized crossed algebras   defined as  the
crossed algebras
above but without axiom (3.2.4).

 \skipaline \centerline  {\bf 5.  Proof of Theorems 4.1 and 4.3}

\skipaline

Throughout this section we shall use the term $X$-surface
for any 2-dimensional  $X$-cobordism. When the underlying surface is a disc, an
annulus or a disc with holes we  call the $X$-surface an $X$-disc, an
$X$-annulus, or an $X$-disc with holes, respectively.

\skipaline \noindent {\bf   5.1.  Proof of Theorem  4.1: the bijectivity for morphisms.}  
We show here that for any two 
$(1+1)$-dimensional HQFT's  $(A,\tau), (A',\tau')$ with target   $X=K(\pi,1)$
with undelying crossed algebras $L,L' $, the homomorphism
$$\Hom ((A,\tau), (A',\tau')) \to  \Hom (L,L')\leqno
(5.1.a)$$
constructed in Section  4.10 is bijective. The injectivity of this homomorphism is obvious since 
all 1-dimensional $X$-manifolds are disjoint unions of loops   and therefore any  two morphisms
$(A,\tau)\to  (A',\tau')$ coinciding on loops coincide on all  1-dimensional $X$-manifolds.

  To establish the surjectivity of (5.1.a) we first    show how to reconstruct 
   $(A,\tau)$  
(at least up to isomorphism) 
from the underlying crossed $\pi$-algebra $L=(L,\eta,\varphi)$.

It   follows from the topological classification of  
surfaces that every compact oriented surface  can be split along a finite
set of disjoint simple loops into a   union of  
 discs with $\leq 2$ holes, i.e.,  discs, annuli and discs with
  two holes.   This implies that every
$X$-surface $W$ can be obtained by gluing from a finite
collection of $X$-surfaces whose underlying surfaces are  discs with
$\leq 2$ holes. Axioms of an HQFT imply that the vector $\tau(W)$  is
determined by the values of $\tau$ on the discs  with $\leq 2$  holes. It
remains to show that these values are completely determined by $(L,\eta,\varphi)$.

We begin by computing   $\tau $ for annuli.
Each  $X$-annulus  is
$X$-homeomorphic either to  $C_{-+} (\alpha; \beta)$,
or to $C_{--} (\alpha; \beta)$, or to $C_{++} (\alpha; \beta)$
with  $\alpha,\beta \in \pi$. (Note that   $C_{+-} $  is
homeomorphic to $C_{-+} $.)
By definition, 
$$\tau (C_{-+} (\alpha; \beta))(a)=\varphi_{\beta^{-1}}(a)\leqno
(5.1.b) $$
for any $a\in  {L_\alpha}$. 
 
The annulus  $C_{--}(\alpha; \beta)$  can be obtained by the gluing of
two annuli $C_{-+}(\alpha; \beta)$ and  $C_{--}(\beta^{-1}\alpha\beta; 1)$ along an
$X$-homeomorphism $(C^1_+,\beta^{-1}\alpha\beta)=(C^0_-,\beta^{-1}\alpha\beta)$.   
Axiom (1.2.6) and the   definition of $\eta$  imply that 
   $$\tau (C_{--}(\alpha; \beta))\in \Hom_K(L_\alpha \otimes 
L_{\beta^{-1}\alpha^{-1}\beta}, K)$$
is computed by 
 $$\tau (C_{--}(\alpha; \beta))(a\otimes b)=
\eta (\varphi_{\beta^{-1}} (a),b) \leqno
(5.1.c) $$  
where  $a\in  L_\alpha $, $b\in L_{\beta^{-1}\alpha^{-1}\beta}$.

To compute  the vector $$\tau (C_{++}(\alpha; \beta))\in
 \Hom_K (K, L_\alpha\otimes L_{\beta^{-1}\alpha^{-1}\beta})=L_\alpha\otimes
L_{\beta^{-1}\alpha^{-1}\beta }$$  we present this vector as a finite sum
$\sum_i a_i\otimes b_i$ where $a_i\in  L_\alpha$ and $b_i\in 
L_{\beta^{-1}\alpha^{-1}\beta }$.
Observe that the
 gluing of $C_{--}(\alpha^{-1};1)$ to 
$C_{++}(\alpha; \beta)$ along an $X$-homeomorphism $(C^1_-,\alpha)=(C^0_+,\alpha)$
yields
  $C_{-+}(\alpha^{-1};\beta)$. The axioms of an HQFT and the computations above yield
 $$\sum_i\eta (a, a_i) \,b_i=\varphi_{\beta^{-1}}(a) \leqno (5.1.d)$$
for any $a\in L_{ \alpha^{-1}}$. Since the  restriction of $\eta$ to 
$ L_{ \alpha^{-1}}\times    L_{ \alpha}$ is non-degenerate, equality (5.1.d)
determines uniquely the vector  $\tau (C_{++}(\alpha; \beta))=\sum_i a_i\otimes b_i$.

There are two $X$-discs:   $B_+$ 
where the orientation of the boundary is induced by the one in the disc
and $B_-$  where the orientation of the boundary is opposite to the one induced from
the disc.    The vector  $\tau (B_+)=1_L\in L_1$
  is   determined uniquely as the unit element of $L$.
The $X$-disc $B_-$ may be obtained by the gluing of $B_+$ and  
$C_{--} (1;1)$  along $\partial B_+=C^0$. Therefore    $\tau (B_-)\in \Hom_K(L_1,K)$
is determined uniquely by $L$. 
One can   compute that   $\tau (B_-)(a)=\eta (a,1_L)=\eta (1_L,a)$ for any  $a\in
L_{1}$.

Each  disc with two holes   can be split
  along  3 disjoint simple loops parallel to its boundary components 
into a   union of 3 annuli 
and a smaller disc  with two holes. Choosing appropriate orientations of these 3 loops
we obtain that any $X$-disc with two holes  $D$
splits as a union of 3 annuli and 
an $X$-disc with two holes $X$-homeomorphic
to $D_{--+} (\alpha, \beta;1, 
1)$ for certain $\alpha, \beta\in \pi$. The homomorphism 
$\tau (D_{--+}(\alpha, \beta;1, 
1))$ is multiplication in    $L$. The values of $\tau $ on the annuli are also determined by 
$L$ by the arguments   above. Now, the axioms of an HQFT
allow us to recover   $\tau (D) $ uniquely.  

Now we can prove the surjectivity of (5.1.a). Every morphism of crossed $\pi$-algebras $\rho: L\to
L'$ defines a linear homomorphism $A_M\to A_{M'}$ for any connected 1-dimensional $X$-manifold
$M$. These homomorphisms extend to non-connected $X$-manifolds $M$ by multiplicativity.
We should   show that the resulting 
family of $K$-linear homomorphisms
$\{\rho_M:A_M\to A'_M\}_M$    makes 
the natural square diagrams associated with $X$-homeomorphisms   and  
$X$-surfaces   commutative. The part concerning the homeomorphisms is obvious. 
As it was explained above,  
every
$X$-surface   can be obtained by gluing from a finite
collection of    $X$-surfaces of type $B_+, C_{-+} (\alpha; \beta), C_{--} (\alpha; 1),
C_{++}(\alpha; \beta)$, and $D_{--+}
(\alpha, \beta;1,  1)$.  Therefore it suffices to check the commutativity of  the 
square diagrams associated with these $X$-surfaces. 
For  
$B_+$ and 
$D_{--+}
(\alpha, \beta;1,  1)$ this   follows from the assumption that 
$\rho: L\to
L'$  is an algebra homomorphism preserving the unit. For the annuli
$C_{-+} (\alpha; \beta), C_{--} (\alpha; 1),
C_{++}(\alpha; \beta)$, this follows from the formulas  (5.1.b) - (5.1.d)  and the assumption that 
$\rho: L\to
L'$  preserves the inner product and commutes with the action of $\pi$.

\skipaline \noindent {\bf   5.2.  Proof of Theorem  4.1: the surjectivity for objects.}  
Let $(L  =\bigoplus_{\alpha\in \pi} L_\alpha, \eta,\varphi)$ be a crossed
$\pi$-algebra. We shall  realize it as the underlying algebra of a $(1+1)$-dimensional HQFT
$(A,\tau)$ with target $X=K(\pi,1)$.

We  associate  with every 1-dimensional  $X$-manifold a $K$-module as in
Section 4.5. It remains to define the homomorphisms  $\tau$ for 2-dimensional
$X$-cobordisms. The construction of $\tau $ goes in nine steps.

Step 1. We define    the homomorphism  $\tau $ for     $X$-annuli
using formulas (5.1.b) - (5.1.d). To establish the topological invariance,
note that the $X$-homeomorphisms of  $X$-annuli  are generated by 
(i)  the Dehn twists 
$C_{\varepsilon,\mu}(\alpha;\beta) \to C_{\varepsilon,\mu}(\alpha;\beta\alpha)$
along the circle $S^1\times (1/2)$ where $\epsilon, \mu=\pm$
and   (ii) 
the homeomorphisms $C_{\varepsilon,\varepsilon}(\alpha;\beta)\to
C_{\varepsilon,\varepsilon}(\beta^{-1}\alpha^{-1} \beta;\beta^{-1})$
 permuting the boundary components of the annulus and preserving the
arc $s\times [0,1]$ (with $s\in S^1$) as a set.

The invariance of $\tau $ under the Dehn twists follows from   the equalities 
$\varphi_{\beta\alpha}\vert_{L_{\alpha}}=
\varphi_{\beta}\varphi_\alpha\vert_{L_{\alpha}}=\varphi_\beta\vert_{L_{\alpha}}$.

The homomorphism   $\tau $ defined
by (5.1.c) is invariant under the   homeomorphism (ii) with $
\epsilon=-$  because $$\tau (C_{--}(\alpha; \beta))(a\otimes b)=
\eta (\varphi_{\beta^{-1}} (a),b)
 = \eta (a,\varphi_{\beta} (b))=\eta (\varphi_{\beta} (b),a)$$
$$=\tau
(C_{--}(\beta^{-1}\alpha^{-1} \beta;\beta^{-1})) (b\otimes a).$$ 
To prove that 
the homomorphism  $\tau(C_{++}) $ defined by (5.1.d) is invariant under the
homeomorphism   $C_{++}(\alpha;\beta)\to
C_{++}(\beta^{-1}\alpha^{-1} \beta;\beta^{-1})$ described above, it suffices to deduce
from (5.1.d) that for any $b\in L_{\beta^{-1}\alpha \beta}$,
 $$\sum_i\eta (b, b_i)\, a_i= \varphi_{\beta}(b). \leqno (5.2.a)
$$
For any $a\in L_{\alpha^{-1}}$, we have
$$\eta(a,\varphi_{\beta}(b))=\eta( \varphi_{\beta^{-1}}(a),b)=\eta(b,
\varphi_{\beta^{-1}}(a))$$
$$=\eta(b,
\sum_i\eta (a, a_i) \,b_i) =\sum_i\eta(a,a_i)\, \eta(b,b_i)=\eta(a, \sum_i\eta (b, b_i)
\, a_i).$$ Now, the non-degeneracy of $\eta$ implies (5.2.a).

Step 2.   
We   check now   axiom  (1.2.6) in the case where 
$W,W_0$ and $W_1$  are  annuli. It suffices to consider
the case where the annulus $C_{\varepsilon,\mu}(\alpha;\beta)$ is glued
to $C_{-\mu,\nu}(\gamma;\delta)$ along  an $X$-homeomorphism
$(C^1_{\mu},(\beta^{-1} \alpha^{-\varepsilon}\beta)^\mu)=(C^0_{-\mu}, \gamma)$. Note
that the gluing is possible only if $\gamma= 
(\beta^{-1} \alpha^{-\varepsilon}\beta)^\mu$; the result of
the gluing is the annulus 
  $C_{\varepsilon, \nu}(\alpha;\beta \delta)$.   We    consider  8 cases depending on the
values of  $\varepsilon,\mu,\nu=\pm$  and indicate the key argument
implying (1.2.6);     the details are left to the reader.

Cases (1) and (2): $\varepsilon=-,\mu=+,\nu= \pm$.
Use that $\varphi_{\delta^{-1}} \varphi_{\beta^{-1}}=\varphi_{(\beta \delta)^{-1}}$.

Case (3): $\varepsilon= \mu= \nu= -$. This case follows from
Case (2) by permuting the annuli under gluing:  the gluing of $C_{--} $ to 
  $C_{+-} $ along an $X$-homeomorphism $C^1_-=C^0_+$ is the same operation as 
the gluing of $C_{-+} $ to 
  $C_{--} $ along an $X$-homeomorphism $C^1_+=C^0_-$.

Case (4): 
$\varepsilon= \mu=-,\nu= +$.
We should prove that
$$(\tau (C_{--}(\alpha;\beta))\otimes \id_{L_{\delta^{-1}\gamma^{-1}\delta}})
\, (\id_{L_\alpha} \otimes 
\tau (C_{++}(\gamma;\delta)))=\tau (C_{-+}(\alpha;\beta\delta))
.$$
By definition, $$\tau (C_{++}(\gamma;\delta)
)=\sum_i c_i\otimes d_i \in L_{\gamma}\otimes
L_{\delta^{-1}\gamma^{-1}\delta}$$ where  
 $$\sum_i\eta (e, c_i) \,d_i=\varphi_{\delta^{-1}}(e) \leqno (5.2.b)$$
for any $e\in L_{\gamma^{-1}}$.
Note that $\gamma=\beta^{-1}\alpha^{-1}\beta$.   For  $a\in
L_\alpha$, we have
$$(\tau (C_{--}(\alpha;\beta))\otimes \id_{L_{\delta^{-1}\gamma^{-1}\delta}})
\, (\id_{L_\alpha} \otimes 
\tau (C_{++}(\gamma;\delta))) (a)
$$
$$= \sum_i (\tau (C_{--}(\alpha;\beta))\otimes
\id_{L_{\delta^{-1}\gamma^{-1}\delta}}) (a\otimes c_i\otimes d_i)$$
$$
= \sum_i  \eta (\varphi_{\beta^{-1}}(a), c_i) d_i
=\varphi_{\delta^{-1}} \varphi_{\beta^{-1}}(a)=
\varphi_{(\beta \delta)^{-1}}  (a)=\tau (C_{-+}(\alpha;\beta\delta)) (a).$$

Case (5): $\varepsilon= \mu= \nu=+$.  Use  that formula
(5.1.d) implies the equality  $ \sum_i\eta (a, a_i)
\,\varphi_{\delta^{-1}} (b_i)=\varphi_{(\beta \delta)^{-1}}(a) $.

Cases (6) and (7):  $\varepsilon=+, \mu= \pm,  \nu=-$.  These
cases follow  from Cases (4) and (1)   by permuting the annuli under
gluing.
 
Case (8):  $\varepsilon=+, \mu=-, \nu=+$.
Use the same expression for $\tau (C_{++}(\gamma;\delta))$
 as in Case (4) and the   equalities
$$   \sum_i\eta (e, \varphi_{\beta} (c_i) )\, d_i = \sum_i\eta
(\varphi_{\beta^{-1}} (e),   c_i) \,d_i = \varphi_{(\beta \delta)^{-1}}(e)  $$
where $e\in L_{\beta \gamma^{-1}\beta^{-1}}$. 

Step 3.  At steps 3 - 5 we define $\tau$ for discs with two holes
  $D_{\varepsilon, \mu,\nu}(\alpha, \beta;\rho, \delta)$
where $\alpha, \beta,\rho, \delta\in \pi$.
 
For an $X$-disc with
two holes   $D=D_{--+} (\alpha, \beta;\rho, 
\delta)$, cf. Section 4.6, we define    $\tau(D)\in  \Hom_K (L_\alpha\otimes L_\beta,
L_{\rho \alpha\rho^{-1} \delta\beta\delta^{-1}})$ by $\tau(D)(a\otimes b)=
\varphi_{\rho}(a)\,\varphi_{\delta}(b)$  where $a\in
L_\alpha, b\in L_\beta$.  The topological
invariance of $\tau (D)$ follows from   axioms (3.2.2,   3.2.3)  and the fact that
any self-homeomorphism of 
  $D=D_{--+}$ is isotopic to a composition of   Dehn twists in   annuli neighborhoods 
of the  circles $Y,Z\subset \partial D $ and the homeomorphisms
$f^{\pm 1}:D \to D $ introduced in the proof of Lemma 4.9.
 
Axiom (1.2.6) holds for any gluing
of an annulus of type $C_{-+}$ to  $D_{--+} $:
if the gluing is performed along an $X$-homeomorphism $C^1_+=Y_-$ or $C^1_+=Z_-$
then this follows from the identity $\varphi_\alpha  \varphi_\beta=\varphi_{\alpha
\beta}$;  if the gluing is performed along $C^0_-=T_+$ then this follows from the
assumption  that each $\varphi_\alpha$
is an algebra homomorphism. 
 
Step 4.   Consider an $X$-disc with two holes  
$D=D_{---}(\alpha, \beta;\rho, \delta)$ and set $$\gamma=\delta\beta^{-1}\delta^{-1}
\rho \alpha^{-1}\rho^{-1}.$$ We define   
$\tau(D)\in\Hom_K (L_\alpha\otimes L_\beta \otimes L_{\gamma},K)$ by $\tau(D)(a\otimes b\otimes c)=\eta(
\varphi_{\rho}(a)\,\varphi_{\delta}(b),c)$   where $a\in L_\alpha, b\in L_\beta, c\in 
L_{\gamma} $. This definition results immediately
from   axiom (1.2.6) if we present 
$D $ as the result of a gluing  of $D_{--+}(\alpha, \beta;\rho, 
\delta)$
to  $C_{--} (\gamma^{-1};1)$ along   $
(T_+,\gamma^{-1})=(C^0_-,\gamma^{-1})$. 

Let us verify the topological invariance of $\tau(D)$. 
 Consider an $X$-homeomorphism $h:D_{---}\to D_{---} $  
  which maps   $(Y,y), (Z,z), (T,t)$ onto $(Z,z), (T,t), (Y,y)$,
respectively. We choose $h$ so that    the arc $tz$ is mapped onto $yt=(ty)^{-1}$ and
the arc $ty$ is mapped onto an embedded arc leading from $y$ to $z$ and homotopic to
the product of  the  arcs $yt$ and $ tz$.  An easy 
  computation   shows that $h$ is an  $X$-homeomorphism
  $D_{---}(\alpha, \beta;\rho, 
\delta) \to D_{---}(\gamma,\alpha ;\delta^{-1},\delta^{-1}\rho)$.
The   invariance of $\tau $ under $h$ follows from
the equalities
$$\eta(\varphi_{\delta^{-1}}(c) \,\varphi_{\delta^{-1}\rho} (a ) ,b)=
\eta(c \,\varphi_{ \rho} (a ) , \varphi_{\delta }(b))=
\eta(c ,\varphi_{ \rho} (a ) \, \varphi_{\delta }(b))
=\eta(\varphi_{ \rho} (a ) \, \varphi_{\delta }(b), c)
 .$$
It is clear that any self-$X$-homeomorphism  of
$D $ may be presented as a composition of   $h^{\pm 1}$ with 
a self-homeomorphism  of
$D  $ preserving all boundary components set-wise.
Therefore it remains only to check the topological invariance of $\tau(D) $ under 
  self-homeomorphisms of
$D  $ preserving the boundary components. Such homeomorphisms 
(considered up to isotopy) are
compositions of   Dehn twists in   annuli neighborhoods  of $Y,Z,T$. Invariance of
$\tau (D )$ under such Dehn twists follows from the already established
topological invariance of the homomorphisms    $\tau (D_{--+}
(\alpha, \beta;\rho, 
\delta))$ and 
$\tau (C_{--} (\gamma^{-1};1))$. 

Axiom (1.2.6) holds for any gluing
of an annulus of type $C_{-+}$ to  $D_{---} $ (such a gluing produces again  
$D_{---}$). If the gluing is performed along $Y$ or $Z$ then this follows from  the
identity $\varphi_{\alpha} \varphi_{\beta}=\varphi_{\alpha\beta}$. The existence of a
self-homeomorphism   of $D_{---}$ mapping $T$ onto $Y$ shows  that the claim holds
also for the  gluings along $T$.

Step 5. Now we define    $\tau $ for $D=D_{++-}(\alpha, \beta;\rho,
\delta) $. We can obtain  $D$ by gluing three annuli    $C_{++}(\alpha;1 ),
C_{++}(\beta;1 ), C_{--}(\gamma;1)$ with $\gamma=
  \rho \alpha^{-1}\rho^{-1}\delta\beta^{-1} \delta^{-1}$ to
    $D_{--+}(\alpha^{-1},
\beta^{-1};\rho, \delta)$ along $X$-homeomorphisms $$(C^1_+, \alpha^{-1})=(Y_-,
\alpha^{-1}), \,\,\,
(C^1_+,\beta^{-1})=(Z_-,\beta^{-1}), \,\,\, (C^0_-,\gamma)=(T_+,\gamma),$$
respectively.   These three annuli  form a regular neighborhood of $\partial D$ in
$D$.
Axiom (1.2.6) determines   
$\tau (D)$.  Its topological invariance   follows from 
the topological invariance of the values of $\tau $ for
$D_{--+}(\alpha^{-1},
\beta^{-1};\rho, \delta)$ and for the three annuli in question  and the  following
obvious fact:    any self-homeomorphism  of   $D$ is
isotopic to a homeomorphism mapping a given regular neighborhood of $\partial D$
onto  itself.

Similarly, we can obtain  $D_{+++}(\alpha, \beta;\rho,
\delta)$ by gluing 3 annuli  of type
$C_{++}(...; 1)$ to   
$D_{---}(\alpha^{-1},
\beta^{-1};\rho, \delta)$. Axiom (1.2.6) determines  
$\tau (D_{+++}(\alpha, \beta;\rho,
\delta))$ in a topologically invariant way.
 
Step 6. We   check   now  axiom (1.2.6)   for a  gluing
of an $X$-annulus   $C_{\varepsilon^0,\varepsilon^1} $ to an $X$-disc with two holes 
$D_{\varepsilon, \mu,\nu} $. By the topological invariance of
$\tau $, it is enough to consider the gluings performed along an $X$-homeomorphism 
$C^0_{\varepsilon^0 }=T_{ \nu}$ so that $\varepsilon^0=-\nu$. We have  16 cases
corresponding to different signs
$ \varepsilon^1,\varepsilon, \mu,\nu$. The cases
where $\varepsilon^0=-\varepsilon^1$ and the triple 
$\varepsilon, \mu,\nu$ contains at least two minuses were
considered   at Steps 3 and 4. The cases
where $\varepsilon=\mu $ are checked one by one using directly the
definitions and the properties of   $\tau $ established above, specifically, axiom
(1.2.6) for    annuli (Step 2).  The key argument in all these cases is that  the  
tensor contractions along different tensor factors commute.
 The case $\varepsilon=-, \mu=+$ reduces
to   $\varepsilon=+, \mu=-$ by the topological invariance.  Assume that
$\varepsilon=+, \mu=-$. If $\nu= +,\varepsilon^1=+$   then (1.2.6) follows again from
definitions. The remaining three cases 
$(\nu= +,\varepsilon^1=-), (\nu= -,\varepsilon^1=\pm)$
  can be   deduced from  the following 
 multiplicativity of $\tau$.

Let $\alpha, \beta,  \rho,
\delta,\Delta,\sigma \in \pi$. Set $\gamma= \rho\alpha\rho^{-1} \delta \beta
\delta^{-1}$.
Observe that the gluing of
$D_{--+}(\alpha, \beta;\rho,\delta)$ to $C_{--}(\gamma;\Delta)$ along
$(T_+,\gamma)=(C^0_-, \gamma)$ yields the $X$-disc with two holes
$D=D_{---}(\alpha, \beta;\Delta^{-1}\rho, \Delta^{-1}\delta)$.
The same $X$-disc with two holes $D$ is obtained by gluing
$\tilde D= D_{-+-}(\alpha, \sigma^{-1}\beta^{-1}\sigma;\Delta^{-1}\rho,
\Delta^{-1}\delta \sigma)$ to $C_{--}(\beta;\sigma)$ along $(Z_+,
\sigma^{-1}\beta^{-1}\sigma)=(C^1_-, \sigma^{-1}\beta^{-1}\sigma)$.
This allows us to compute 
$$
\tau(D)\in
\Hom_K(L_\alpha\otimes L_\beta \otimes L_{\Delta^{-1}\gamma^{-1}\Delta}, K)$$ 
applying (1.2.6) to these two splittings of $D$.   We claim that  these two
computations  give the same result. The first splitting implies that for any $a\in
L_\alpha, b\in  L_\beta, c \in  L_{\Delta^{-1}\gamma^{-1}\Delta}$,
$$\tau(D) (a\otimes b\otimes c)=\tau
(C_{--}(\gamma;\Delta)) (\tau(D_{--+}(\alpha, \beta;\rho,\delta)) (
a\otimes b)\otimes c)$$
$$=\eta(\varphi_{\Delta^{-1}}(
\varphi_{\rho}(a) \varphi_{\delta }(b))\,,\,c).$$
 To use the second splitting we first observe that 
$\tilde D$ is $X$-homeomorphic to
$$D_{--+}(\Delta^{-1}\gamma^{-1}\Delta,\alpha ; 
\sigma^{-1} \delta^{-1} \Delta, \sigma^{-1} \delta^{-1}\rho)$$
 via a homeomorphism mapping the
  boundary components $Y,Z,T$ onto  $Z,T,Y$, respectively.
Therefore applying   (1.2.6) to the second splitting  of $D$ we obtain
$$\tau(D) (a\otimes b\otimes c)=
\tau(C_{--}(\beta;\sigma))(b  \otimes \tau (\tilde D)(c \otimes a))$$
$$=\eta(\varphi_{\sigma^{-1}}(b), \varphi_{\sigma^{-1} \delta^{-1} \Delta}(c)
\,\varphi_{\sigma^{-1} \delta^{-1}\rho}(a))
=\eta(b, \varphi_{  \delta^{-1} \Delta}(c)
\,\varphi_{  \delta^{-1}\rho}(a))$$
$$=
\eta(  \varphi_{  \delta^{-1} \Delta}(c)
\,\varphi_{  \delta^{-1}\rho}(a),b) 
 =
\eta(  \varphi_{  \delta^{-1} \Delta}(c),
\varphi_{  \delta^{-1}\rho}(a) b)=
\eta(  
\varphi_{  \delta^{-1}\rho}(a) b,\varphi_{  \delta^{-1} \Delta}(c))$$
$$=
\eta(  
\varphi_{  \Delta^{-1}\rho}(a) \,\varphi_{  \Delta^{-1}\delta}(b),c)
=\eta(\varphi_{\Delta^{-1}}(
\varphi_{\rho}(a)\, \varphi_{\delta }(b)),c).$$
This completes the check of consistency and implies 
  (1.2.6)   for any  gluing
of an $X$-annulus    to an $X$-disc with two holes.

 Step 7.  Now we define   $\tau $ for the $X$-discs $B_+$
and $B_-$ (see Section 5.1 for notation).
Set  $\tau (B_+)=1_L\in L_1$.  
Axiom (1.2.6)   for a  gluing 
of $B_+$ to an $X$-annulus  of type $C_{-+}$ follows from the equality
$\varphi_\beta(1_L)=1_L$ for $\beta\in \pi$. Axiom (1.2.6)   for a  gluing 
of $B_+$ to  an
$X$-disc with two holes of type $D_{--+}$ follows from the equalities $1_La=a1_L=a$
for any $a\in L$. This and  the definition  of    
$\tau $ for $X$-discs with two holes of types $D_{---}$ or $D_{-++}$    imply  
(1.2.6)    for any  gluing  of $B_+$ to such discs with holes.

The disc $B_-$ can be obtained by the gluing of $B_+$
and   $C_{--} (1;1)$  along $\partial B_+=C^0_-$. This   determines    $\tau (B_-)\in
\Hom_K (L_1,K)$. Axiom (1.2.6)   for a  gluing 
of $B_-$ to   
$X$-discs with $\leq 2$ holes follows from the already established 
  properties    of  the  gluings 
of $C_{--} (1;1)$ and $B_+$.

Step 8. Now we   define    $\tau (W)$ for any connected $X$-surface
$(W,g:W\to X)$. By a splitting system of loops on $W$ we mean a finite set
of disjoint 
 embedded circles $\alpha_1,...,\alpha_N\subset W$ which split $W$ into a   union of  
  discs with  $\leq 2$  holes. 
 We provide each   $\alpha_i$ with an orientation and a
base point $x_i$.  Replacing if necessary $g$ by a homotopic map we can assume that
$g(x_i)=x\in X$ for all $i$.  
 The   discs with   holes obtained by the splitting of $W$ along
$\cup_i \alpha_i$ endowed with the restriction of $g$ are
$X$-surfaces.  
Axiom (1.2.6) determines   $\tau (W)$  from  the values of $\tau $ on these
discs with holes. 

We claim  
that  $\tau (W)$ does not depend on the choice of orientations and base
points on  $\alpha_1,...,\alpha_N$. The proof is as follows. Choose $i\in
\{1,...,N\}$ and set $\alpha=\alpha_i$.
Let $D_1, D_2$ be the discs with holes attached to $\alpha$ from two sides (they may
coincide). Let $C$ be a regular neighborhood of $\alpha$ in $D_1$. Clearly, $C$ is an
annulus in $D_1$ bounded by $\alpha$ and a parallel loop ${\tilde \alpha}$.   Consider the discs
with holes  $\tilde D_1 =\overline {D_1\backslash C}$ and 
$\tilde D_2 =\overline {D_2\cup C}$. We provide ${\tilde \alpha}$ with 
  a  base point $\tilde x$
and
the orientation opposite to the one of $\alpha$. 
We 
deform   $g$ in a small 
neighborhood of $\tilde x$
so that
   $g(\tilde x)=x$.   
In these way the  surfaces $C, \tilde D_1,\tilde D_2$
  acquire the structure of $X$-surfaces. It follows from the properties of
$\tau $ established above
that
$$\ast_{\alpha} (\tau (D_1) \otimes \tau (D_2))=
\ast_{\alpha} \ast_{{\tilde \alpha}} (\tau (\tilde D_1) \otimes \tau (C) \otimes  \tau (D_2))
=  \ast_{{\tilde \alpha}} (\tau (\tilde D_1) \otimes    \tau (\tilde D_2)) $$
where $\ast_{\alpha}$ denotes the tensor contraction corresponding to the gluing
along   $\alpha$.  Thus, replacing $\alpha$ with ${\tilde \alpha}$ in
our splitting system of loops we do not change  $\tau (W)$. This implies
that $\tau (W)$ does not depend on the choice of orientations and base
points on  $\alpha_1,...,\alpha_N$. A similar argument shows that $\tau (W)$ does not
depend on the choice of $g$ in its homotopy class (relative to the base points on
$\partial W$).

We claim that  the homomorphism  
$\tau (W)$  does not depend on the choice of a splitting system of loops on $W$.
The crucial argument  is provided by the fact (see [HT]) that any two splitting systems
of loops on $W$ are related by  the following transformations:

(i) isotopy  in $W$;

(ii) adding to   a splitting system of loops $\{\alpha_1,...,\alpha_N\}$
a  simple    loop  $\alpha \subset W \backslash \cup_i \alpha_i$; 

(iii) deleting a loop from   a splitting system of loops, provided the remaining loops
form a splitting system;  

(iv) replacing  one of the loops $\alpha_i$ of a splitting system  adjacent to two
different discs with two holes $D_1,D_2$  by a  simple loop
lying in $\Int D_1\cup \Int D_2 \cup \alpha_i$, meeting $\alpha_i$ transversally  in  
two points  and  splitting  both $D_1$ and $D_2$ into   annuli;

(v) replacing  one of the loops of a splitting system by a  simple loop
meeting it transversally in one  point   and disjoint from the other loops.

We should check the  invariance of $\tau (W)$ under these transformations. The
invariance of $\tau (W)$ under isotopy is obvious.  Consider the transformation (ii).
The loop $\alpha$ lies in a disc with $\leq 2$ holes obtained by the splitting of
$W$ along  $\cup_i \alpha_i$. The loop $\alpha$ splits this disc with $\leq 2$
holes into a union of a smaller disc with $\leq 2$ holes and   an annulus or a disc
(without holes). Therefore  the
invariance of $\tau (W)$ under   (ii) follows from the already
established multiplicativity of $\tau$ under   the gluings of 
a disc with $\leq 2$ holes to an annulus or a disc. The transformation (iii) is
inverse to (ii) and the same argument applies.

The associativity of multiplication in $L$ yields two equivalent
expressions for   the value of $\tau $ for a disc with 3 holes; these expressions are
obtained from two splittings of the disc with 3 holes as a union of two discs with 2
holes (cf. Section 4.7). Since we are free to choose orientations of the loops in a
splitting system, we can
  reduce the invariance of $\tau (W)$ under the transformation (iv) to this
model case. Note that in order to have the maps to $X$  as in the
model case we can use     transformations  (ii) to   add additional
annuli to the splitting. Similarly, the
invariance of $\tau (W)$ under   the transformation (v) follows from axiom 
(3.2.4) (cf. the end of the proof of Lemma 4.9). 

The topological invariance of $\tau (W)$ follows
from the topological invariance of $\tau $ for discs with $\leq 2$ holes and the fact that
any homeomorphism of connected surfaces maps a splitting system of loops onto a
splitting system of loops.

 Step 9.  We have defined   $\tau $ for connected $X$-surfaces.
We extend    $\tau $ to arbitrary $X$-surfaces by (1.2.5).
It follows directly
from definitions that $(A,\tau)$ is an HQFT.  
 (To prove (1.2.6) we
compute $\tau (W)$ using a splitting system of loops containing  
$N=N'\subset W$.)

\skipaline \noindent {\bf   5.3.  Last claim of Theorem  4.1.} The last claim of Theorem 4.1
is obvious for  direct sums, tensor products, and pull-backs. 
For the duality and transfer, this claim is a nice computational exercise.

\skipaline \noindent {\bf   5.4.  Proof of Theorem  4.3.}  Theorem  4.3  results from 
Theorems 3.6, 3.7,   4.1 and the fact that  the underlying crossed  algebra of a
$(1+1)$-dimensional primitive cohomological HQFT is the crossed  algebra defined in
Section 3.3.

\skipaline \centerline  {\bf 6. Hermitian and unitary crossed algebras}

\skipaline

In this section we assume that $K$ has a ring   involution $k\mapsto
\overline k:K \to K$.

 \skipaline \noindent {\bf   6.1.  Hermitian and unitary crossed $\pi$-algebras.}  Let
$\pi$ be a group.   A {\it Hermitian crossed $\pi$-algebra}   is a 
crossed $\pi$-algebra   $(L=\bigoplus_{\alpha\in \pi} L_{\alpha}, \eta, \varphi)$ over
$K$ endowed with an involutive antilinear antiautomorphism 
$a\mapsto
\overline a:L \to L$ which commutes with the action of $\pi$,
transforms $\eta$ into $\overline \eta$ and 
sends  each 
$ {L_{\alpha}}$   into $L_{\alpha^{-1}}$.
 This definition implies the identities $$
\overline{\overline a}=a,\,\,\,\overline{ka}= \overline k \overline a,\,\,\,
\overline{ab}= \overline b \overline a,\,\,\,
\varphi_{\beta} (\overline a)=\overline {\varphi_{\beta} (a)},\,\,\,
\eta( \overline a,   \overline b)=  \overline {\eta(a,b)} \leqno (6.1.a)$$ for any $a,b\in L, k\in K,
\beta\in \pi$. It follows from these conditions  that $\overline
{1_L}=1_L$.

Observe that for any $a\in L$,
$$  \overline {\eta(a,\overline a)} =\eta( \overline a,   a)= {\eta(a,\overline a)}.$$
In particular, 
if $K=\bold  C$ with   usual complex conjugation then ${\eta(a,\overline a)}\in \bold 
R$ for all $a$. If additionally,  
  $\eta (a,  \overline a)  >0 $ for all non-zero $a\in L$, then we say
that  $L$ is {\it unitary}. 

It is    easy   to check that direct sums, tensor products,  pull-backs, duals,  and 
transfers of  Hermitian (resp. unitary) crossed  algebras are again   Hermitian (resp.
unitary) crossed  algebras. We can define   a category, $HQ_{2}(\pi)$ (resp. 
$UQ_{2}(\pi)$),  whose
objects are
  Hermitian (resp. unitary)  crossed $\pi$-algebras. The morphisms     are
defined as in Section 3.2 with the additional condition that they
commute  with the involutions.

 \skipaline \noindent {\bf   6.2.  Examples.}  Let 
$\{\theta_{\alpha,\beta} \in S\}_{\alpha,\beta\in \pi}$
 be a normalized 2-cocycle of the group $\pi$ with values in the
multiplicative group $S= \{k\in K\,\vert\,
k \overline k =1_K\}$. We introduce a Hermitian structure on the crossed $\pi$-algebra
$L=\bigoplus_{\alpha\in \pi} K l_\alpha$ constructed in Section 3.3. 
 For all $\alpha\in \pi$, set 
$$\overline {l_\alpha}= \overline {\theta_{\alpha,\alpha^{-1}}} \,l_{\alpha^{-1}}=
({\theta_{\alpha,\alpha^{-1}}})^{-1} l_{\alpha^{-1}}.$$
This extends uniquely to an antilinear homomorphism 
$a\mapsto
\overline a:L \to L$. We claim that it satisfies   (6.1.a). 
The involutivity follows from the equalities
$$\theta_{\alpha,\alpha^{-1}} \overline {\theta_{\alpha^{-1},\alpha}} =
\theta_{\alpha,\alpha^{-1}} \overline {\theta_{\alpha,\alpha^{-1}}}=1_K $$
where we use the identity $\theta_{\alpha^{-1},\alpha}=\theta_{\alpha,\alpha^{-1}}$.
Similarly,
$$\eta (\overline {l_\alpha},\overline {l_{\alpha^{-1}}})=
\eta(  \overline {\theta_{\alpha,\alpha^{-1}}} \, l_{\alpha^{-1}},
\overline {\theta_{\alpha^{-1},\alpha}} \, l_{\alpha })
=  \overline {\theta_{\alpha,\alpha^{-1}}}\, \overline {\theta_{\alpha^{-1},\alpha}}
\, \eta(  l_{\alpha^{-1}},  
l_{\alpha })$$
$$=
\overline {\theta_{\alpha,\alpha^{-1}}} \, \overline
{\theta_{\alpha^{-1},\alpha}} \, {\theta_{\alpha^{-1},\alpha}}=
\overline {\theta_{\alpha,\alpha^{-1}}}=\overline {\eta ( {l_\alpha},
{l_{\alpha^{-1}}})}.$$
  To prove that this involution is an antiautomorphism we should check that
 $\overline {l_\alpha l_\beta}=\overline {  l_\beta} \,\overline {l_\alpha }$
for any $\alpha, \beta \in \pi$. This is equivalent to
$$(\theta_{\alpha,\beta}  \theta_{\alpha \beta,(\alpha\beta)^{-1}})^{-1}
=(\theta_{\alpha,\alpha^{-1}} \theta_{\beta,\beta^{-1}})^{-1}
\theta_{\beta^{-1}, \alpha^{-1}}. \leqno  (6.2.a)
$$
To prove this formula,   set $\gamma=\beta^{-1}$ in (3.3.a).
This  yields 
$ \theta_{\alpha\beta,\beta^{-1}} =\theta_{\beta,\beta^{-1}} (\theta_{\alpha,\beta})^{-1}$.
(We use that $\theta_{\alpha,1}=1$, see Section 3.3). 
Now, we  replace $\alpha, \beta,\gamma$ in (3.3.a)
with $\alpha  \beta,   \beta^{-1},\alpha^{-1}$, respectively.
Substituting  in the resulting formula the expression
$ \theta_{\alpha\beta,\beta^{-1}} =\theta_{\beta,\beta^{-1}} (\theta_{\alpha,\beta})^{-1}$,
we obtain a formula equivalent to  (6.2.a).
It remains to check the identity
$\varphi_{\beta} (\overline {l_\alpha})=\overline {\varphi_{\beta} ({l_\alpha})}$.
Note first that $ {l_\alpha} \overline {l_\alpha}=1$. Therefore
$\varphi_{\beta}( {l_\alpha})\, \varphi_{\beta}(\overline {l_\alpha})=1$.
This characterizes $\varphi_{\beta}(\overline {l_\alpha})$ as the (unique)  element of
$L_{\beta \alpha^{-1} \beta^{-1}}$  inverse to 
$\varphi_{\beta}( {l_\alpha})$. We claim that
$\overline {\varphi_{\beta} ({l_\alpha})}\in L_{\beta \alpha^{-1} \beta^{-1}}$   is also inverse
to  $\varphi_{\beta}( {l_\alpha})$. By definition,
$\varphi_{\beta}( {l_\alpha})=y  l_{\beta \alpha  \beta^{-1}}$ where
$y \in K$ is determined
from the equation $\varphi_{\beta}( {l_\alpha}) l_{\beta}=  l_{\beta} l_\alpha$.
Thus,
 $y = (\theta_{\beta\alpha\beta^{-1},\beta})^{-1} \theta_{\beta,\alpha}$.
This implies $ y \overline y=1$ and therefore 
$\varphi_{\beta}( {l_\alpha})  \overline {\varphi_{\beta} ({l_\alpha})}
= y \overline y  l_{\beta \alpha  \beta^{-1}} \overline {l_{\beta \alpha  \beta^{-1}}}=1$.
Hence $\varphi_{\beta} (\overline {l_\alpha})=\overline {\varphi_{\beta} ({l_\alpha})}$.

If   $K=\bold  C$ with   usual complex conjugation then $ \eta(l_{  \alpha  },\overline
{l_{  \alpha  }}) =1$ for all $\alpha$ so that $L$ is unitary.

\skipaline \noindent {\bf   6.3.  Theorem. }   {\sl   Let $\pi$ be a group.
The underlying 
crossed
$\pi$-algebra  of a Hermitian   $(1+1)$-dimensional HQFT  with target   $K(\pi,1)$ 
has   a Hermitian structure in a natural way. This establishes  
 an equivalence   between   the category $HQ_2(K(\pi,1))$ of 
$(1+1)$-dimensional Hermitian HQFT's with target   $K(\pi,1)$ and   the category $HQ_2( \pi )$
of  Hermitian crossed
$\pi$-algebras.  Similar results hold in the unitary setting. }

\skipaline  {\sl  Proof.} Consider a Hermitian HQFT $(A,\tau)$ with target $K(\pi,1)$
and its underlying 
crossed
$\pi$-algebra $L$. The Hermitian structure 
on   $(A,\tau)$ 
yields for each $\alpha\in \pi$  a non-degenerate Hermitian pairing $\langle
.,. \rangle_{\alpha}:L_{\alpha}\times L_{\alpha} \to K$. 
By the non-degeneracy of $\eta$ there is a unique antilinear isomorphism
$b\mapsto \overline {b}: L_{\alpha} \to L_{\alpha^{-1}}$ such that
$$ \langle
a,b\rangle_{\alpha} =\eta (a, \overline {b})\leqno (6.3.a)$$ for any $a,b\in L_\alpha$.
The identity $\langle
a,b\rangle_{\alpha}=\overline {\langle
b,a\rangle_{\alpha}}$   may be rewritten 
as  $$\eta(  a,   \overline b)=  \overline {\eta(\overline a,b)}  \leqno (6.3.b)$$ for any $a,b\in L$.
We shall check that the homomorphism
$b\mapsto \overline {b}$ defines  a Hermitian structure on $L$.
To this end we shall apply axiom (1.7.2) to 
the $X$-surfaces  of types
$C_{-+}(\alpha;\beta^{-1}), C_{--} (\alpha;1)$ and $ D_{--+} (\alpha,\beta;1,1)$ (cf.
Section 4.6).

For $W=C_{-+}(\alpha;\beta^{-1})$, we have $-W=C_{-+}(\beta\alpha\beta^{-1};\beta)$.
By definition, $\tau(W)= \varphi_{\beta}: L_\alpha\to L_{\beta\alpha\beta^{-1}}$
and
$\tau(-W)= \varphi_{\beta^{-1}}:L_{\beta\alpha\beta^{-1}}\to L_\alpha$.
Axiom (1.7.2) yields
$$\langle \varphi_{\beta} (a), b\rangle_{\beta\alpha\beta^{-1}}=\langle a, \varphi_{\beta^{-1}}
(b)
\rangle_{\alpha} $$
for any $a\in L_\alpha, b\in L_{\beta\alpha\beta^{-1}}$.
Substituting $c=\varphi_{\beta^{-1}}(b)\in L_\alpha$ we obtain an equivalent formula
$$\langle \varphi_{\beta} (a), \varphi_{\beta} (c)\rangle_{\beta\alpha\beta^{-1}}=\langle a,
c\rangle_{\alpha} $$ for any $a,c\in L_\alpha$. This   is equivalent to
 $\eta(\varphi_{\beta} (a), \overline {\varphi_{\beta} (c)})=\eta (  a,
\overline c)$. Since the form $\eta$ is invariant under  $\varphi_{\beta}$ and non-degenerate,
the latter formula is equivalent to
$\varphi_{\beta} (\overline c)=\overline {\varphi_{\beta}(c)}$.

For $W=C_{--}(\alpha;1)$, we have $-W=C_{++}( \alpha ;1)$.
By definition, 
$\tau(W)$ is the pairing $a \otimes b\mapsto \eta(a,b): L_\alpha\otimes 
L_{ \alpha ^{-1}} \to K$ where $a  \in  L_\alpha, b\in L_{ \alpha ^{-1}} $.
The homomorphism $\tau(-W):K\to L_\alpha\otimes 
L_{ \alpha ^{-1}}$ sends $1_K$ into 
a  sum
$\sum_i a_i\otimes b_i$ with $a_i\in  L_\alpha, b_i\in 
L_{\alpha^{-1} }$  such that $\sum_i\eta (a, a_i) \,b_i=a $
for any $a\in L_{ \alpha^{-1}}$ (cf. Section 5.1).
Axiom (1.7.2) for $W$ is equivalent to 
$$ \langle \eta(a,b), 1_K\rangle_{\emptyset}= \langle a \otimes b, \sum_i a_i\otimes    b_i
\rangle $$
where on the right-hand side we have the Hermitian form on 
$L_\alpha\otimes 
L_{ \alpha ^{-1}}$ induced by the Hermitian forms on the factors.
By (1.7.1), $ \langle \eta(a,b), 1_K\rangle_{\emptyset}=  \eta(a,b)$.
On the right-hand side we have
$$\langle a \otimes b, \sum_i a_i\otimes    b_i
\rangle =
 \sum_i \langle a ,   a_i \rangle \langle b, b_i
\rangle
= \sum_i  \eta (a, \overline {a_i})\,\eta (b, \overline {b_i})$$
$$
=\overline {\sum_i  \eta (\overline a,  {a_i})\,\eta (\overline b,  {b_i})}
 =\overline {\eta (\overline b, \sum_i  \eta (\overline a,  {a_i})  {b_i})} 
 =\overline {\eta (\overline b,  \overline a )}
 =\overline {\eta (\overline a,  \overline b )}
.$$
This gives   $\overline{\eta(a,b)}=  {\eta (\overline a,  \overline b )}$.
Using this formula and the non-degeneracy of $\eta$ we observe that (6.3.b) 
is equivalent to the involutivity of   $a\mapsto \overline a:L\to L$.

It remains to analyze the case $W=D_{--+} (\alpha,\beta;1,1)$.
By definition, 
$\tau(W)$ is the multiplication $a \otimes b\mapsto ab: L_\alpha\otimes 
L_{ \beta} \to L_{\alpha \beta}$ where $a  \in  L_\alpha, b\in L_{ \beta} $.
Clearly, $-W=D_{++-} (\beta, \alpha;1,1)$.
The homomorphism $\tau(-W): L_{\alpha \beta}\to L_{ \beta} \otimes L_\alpha$
sends any $c\in L_{\alpha \beta}$ to a sum 
$\sum_j t_j\otimes z_j$ with $t_j\in   L_{ \beta}, z_j \in  L_\alpha$
such that $cy= \sum_j  \eta(y,t_j) z_j$ for all $y\in L_{\beta^{-1}}$.
Axiom (1.7.2) for $W$ is equivalent to 
$$ \langle ab, c\rangle_{\alpha \beta}=\sum_j \langle a  , z_j \rangle_\alpha \,
\langle  b, t_j\rangle_{\beta}$$ for any   $a  \in  L_\alpha, b\in L_{ \beta}, c\in L_{\alpha
\beta}$. The left-hand side is equal to $\eta (ab, \overline c)=\eta (a, b \overline c)$
while the right-hand side is equal to
$$\sum_j \langle a  , z_j \rangle_\alpha 
\langle  b, t_j\rangle_{\beta}=\overline {\sum_j \eta( \overline a  , z_j) \,
\eta( \overline b  , t_j) } = \overline {\eta( \overline a  ,  \sum_j 
\eta( \overline b  , t_j) z_j) }= \overline {\eta( \overline a  ,c \overline b)}
=\eta(  a  , \overline{c \overline b}).$$
Thus   $\eta (a, b \overline c)= \eta(  a  , \overline{c \overline b})$
for all $a,b,c\in L$. By the non-degeneracy of $\eta$,  
  $b \overline c=\overline{c \overline b}$. Hence
$  \overline {b}  \overline c=\overline{c   b}$. 
We conclude  that   
$b\mapsto \overline {b}$ is  a Hermitian structure on $L$. 

Conversely, having a Hermitian structure on $L$  we define  a Hermitian structure on 
$(A,\tau)$ by (6.3.a) and (1.7.1). The arguments above  verify axiom (1.7.2) for the $X$-surfaces 
$C_{-+}(\alpha;\beta^{-1}), C_{--} (\alpha;1), D_{--+} (\alpha,\beta;1,1)$.
Since $C_{++}=-C_{--}$ and  
 (1.7.2) for $W$ is equivalent to  
 (1.7.2) for $-W$, this axiom holds also for  $W=C_{++}$.
Clearly $-B_+=B_-$. Therefore for $W=B_+$, axiom   (1.7.2) is equivalent to 
$$\langle 
  1_L, b \rangle_{1}= \langle \tau(B_+)(1_K), b \rangle_{1}= \langle 1_K,\tau(B_-) (b)
\rangle_{\emptyset}=  \overline {\tau(B_-) (b)}\leqno (6.3.c) $$
  for any  $  b\in L_1$. Recall (see Section 5.1) that $\tau(B_-) (b)=\eta(b,1_L)\in K$. 
Therefore (6.3.c) is equivalent to
$\langle 
 1_L, b \rangle_{1}= \overline {\eta(b,1_L)}$. Applying the conjugation in $K$, we obtain  
an equivalent formula
$\langle 
 b, 1_L \rangle_{1}=  {\eta(b,1_L)}$. Since $ \langle 
 b, 1_L \rangle_{1}=\eta (b, \overline {1_L})$, formula (6.3.c) follows from  $\overline
{1_L}=1_L$.

Now,  every $X$-surface $W$ can be obtained by gluing from the discs with holes of types
$B_+,C_{-+}(\alpha;\beta), C_{--} (\alpha;1), C_{++}  (\alpha;1), D_{--+} (\alpha,\beta;1,1)$ (cf.
Section 5.1). Therefore  (1.7.2) holds for $W$.
Thus, the functor $HQ_2(K(\pi,1))\to HQ_2( \pi )$ is surjective on   objects.
That this functor  is an equivalence of categories is  straightforward (cf. Section 5.1).

\skipaline \centerline {\bf Part III.  Lattice models for 
$(1+1)$-dimensional HQFT's}

\skipaline

\skipaline \centerline {\bf 7.  Biangular $\pi$-algebras and lattice HQFT's}

 \skipaline By a lattice topological quantum field theory one means a TQFT which can be
  computed
in terms of partition  functions (or state sums) on triangulations or cellular decompositions of
manifolds. Lattice TQFT's are well known in dimensions 2 and 3, see [BP], [FHK], [TV].  In 
this
section we   introduce a lattice HQFT in dimension 2.

\skipaline \noindent {\bf 7.1.  Biangular $\pi$-algebras.} 
Let $L=\oplus_{\alpha\in \pi} L_\alpha$ be a
 $\pi$-algebra.  Given $\ell\in 
L$, denote 
by 
$\mu_\ell$ the left multiplication by $\ell$ sending any $a\in L$ into 
$\ell 
a\in 
L$. 
Clearly, if $\ell\in L_1$  then $\mu_\ell(L_\alpha)\subset L_\alpha$ for  
all 
$\alpha\in \pi$.  We 
call 
$L$ a {\it 
biangular $\pi$-algebra} if it satisfies the following two conditions:

(i)  for any $\ell\in L_1$ and any
$\alpha \in \pi$,
$$\Tr(\mu_\ell\vert_{L_\alpha}:L_\alpha\to 
L_\alpha) =\Tr(\mu_\ell\vert_{L_1}:L_1\to 
L_1), \leqno (7.1.a)$$

(ii) for any 
$\alpha \in \pi$, the bilinear form $\eta:L_\alpha\otimes L_{\alpha^{-1}}\to K$
defined by 
$$\eta(a,b)=\Tr(\mu_{ab}\vert_{L_1}:L_1\to 
L_1)\leqno (7.1.b)$$
(where $a\in L_\alpha, b\in L_{\alpha^{-1}}$) is non-degenerate.

Condition (i)   implies   that 
 $$\Dim (L_\alpha)=\Tr(\id:L_\alpha\to 
L_\alpha)=\Tr(\mu_1 \vert_{L_\alpha}:L_\alpha\to 
L_\alpha)=\Tr(\mu_1\vert_{L_1}:L_1\to 
L_1)  $$ 
does not depend on $\alpha\in \pi$. This shows in particular that the crossed
$\pi$-algebras are not necessarily biangular.

The bilinear form $\eta$ on a biangular $\pi$-algebra 
$L$
  defined by   (7.1.b) is symmetric: if $a\in L_\alpha, b\in L_{\alpha^{-1}}$
  then
  $$\eta(b,a)=\Tr(\mu_{ba}\vert_{L_1})=\Tr(\mu_{b}\mu_a \vert_{L_1}:L_1\to 
L_1)=\Tr(\mu_a\mu_b\vert_{L_\alpha}:L_\alpha\to 
L_\alpha)$$
$$ = \Tr(\mu_{ab} \vert_{L_\alpha}:L_\alpha\to 
L_\alpha)=\Tr(\mu_{ab}\vert_{L_1}:L_1\to 
L_1)=\eta(a,b).$$
It is clear  that the inner product $\eta$ 
  makes $L$ a Frobenius $\pi$-algebra.

 The group ring $K[\pi
]$ with canonical $\pi$-algebra structure is  biangular. More generally,   the  
$\pi$-algebra 
constructed in Section 3.3 is  biangular.
 It follows from definitions and the standard properties of the trace that the 
direct sums and tensor products of 
biangular  $\pi$-algebras are biangular  $\pi$-algebras.
More examples of biangular  algebras can be obtained using the pull-back 
and 
 push-forward constructions, see Section 3.1. Namely,  the pull-back   of a 
 biangular $\pi$-algebra along any group homomorphism  
 $q:\pi' \to \pi$ is a biangular $\pi'$-algebra. If $q$ is 
a surjective homomorphism with finite kernel and the order of $\Ker \,q$ is 
invertible in $K$, then the push-forward along $q$ of a biangular $  \pi'$-algebra is a 
biangular $\pi$-algebra. Further examples of biangular  algebras will be 
given in Section 10.

Our interest in biangular $\pi$-algebras is due to the fact that   each such 
algebra $L$ gives rise to a
 $(1+1)$-dimensional HQFT with target $K(\pi,1)$.  The underlying crossed
 $\pi$-algebra of this HQFT is  called the {\it  $\pi$-center} of $L$. (This term is justified by
the examples given in Section 10.2).
 We give here a direct
   algebraic description of  the $\pi$-center of $L$.

    Observe first that for every $\alpha\in \pi$,
 the non-degenerate form $\eta:L_\alpha\otimes L_{\alpha^{-1}}\to K$     
yields   a 
canonical element $b_\alpha \in
L_\alpha\otimes L_{\alpha^{-1}}$  (cf.  Lemma 7.1.1 below).
We  expand  
$$b_\alpha= \sum_i p_i^{\alpha}\otimes q_i^{\alpha} 
\leqno (7.1.c)$$ 
where $i$ runs over a   finite set of 
indices and $p_i^{\alpha} \in L_{\alpha}, q_i^{\alpha} \in 
L_{\alpha^{-1}}$. Since   $\eta$ is symmetric,  
the vector $b_{\alpha^{-1}}$ is obtained from $b_\alpha$ by permutation of 
the 
tensor factors.
The element $b_\alpha$ is characterized by the following property: for any  
$p\in 
L_\alpha $,
$$p=\sum_{i }\,\eta(p, q_i^{\alpha} )\, p_i^{\alpha}.  $$
Note that the sum $ \sum_i p_i^{\alpha}  
q_i^{\alpha} 
\in L_1$ does not depend on the choice of  expansion (7.1.c).
Below we shal prove that for all $\alpha\in \pi$,  
$$\sum_i \,p_i^{\alpha}  \,
q_i^{\alpha}  = 1_L.\leqno (7.1.d)$$

Define   a 
$K$-homomorphism
 $\psi_\alpha:L\to L$   by $\psi_\alpha(a)= \sum_{i} 
 p_i^{\alpha}\,a\, q_i^{\alpha}$
 where $a\in L$.  This
 homomorphism
 does not depend on the choice of   expansion (7.1.c)   and sends 
each $ 
 L_\beta  \subset L$ into $L_{ \alpha\beta \alpha^{-1}}$. 
  It turns out that
$\psi_\alpha
 \psi_{\alpha'}=\psi_{\alpha \alpha'}$ for any $\alpha, \alpha' \in \pi$
and
 $\psi_\alpha$ is adjoint to $ \psi_{\alpha^{-1}}$ with respect to $\eta$.
 In particular, $\psi_1$ is an orthogonal projection onto a submodule
 $C=\oplus_{\alpha\in \pi} C_\alpha $ of $ L$ where $C_\alpha= \psi_1
 (L_\alpha) \subset L_\alpha$. 
  Formula (7.1.d) implies that $1_L\in C$.
   It turns out that $C$ is a subalgebra of $L$.
 Moreover, the pairing $\eta\vert_C:C\otimes C\to K$    and the homomorphism 
$\pi\to \Aut C$ sending each 
$\alpha\in 
\pi$
 to $\psi_\alpha\vert_C:C\to C$ make $C$ a crossed $\pi$-algebra.  
 The crossed $\pi$-algebra $C$ is   the $\pi$-center of
$L$.
All 
these
 claims follow directly by applying the  definitions of Section 4 to the 
HQFT
 associated with $L$ (cf. Section  8.6). 

It is easy to check that if $L=K[\pi]$ or more generally if $L$ is the 
$\pi$-algebra constructed in 
Section 3.3, then $\psi_1=\id_L$ so that $C^L=L$.

We finish this subsection with a proof of the existence of $b_\alpha$,  check 
(7.1.d),
and state a lemma which will be used in Section 8. In the remaining part of 
Section 7 we   give a 
construction 
of the lattice HQFT associated with a biangular 
$\pi$-algebra $L$.   A dual construction of the same HQFT is given in Section 
8.

   \skipaline \noindent {\bf 7.1.1.  Lemma.  } {\sl Let $P,Q$ be 
projective $K$-modules of finite type and let $\eta:  P\otimes Q\to K$ be a 
non-degenerate
bilinear form.  There is a unique element $b\in P\otimes Q$ such that 
for any expansion   $b=\sum_i p_i\otimes q_i$ into a finite sum with   
$p_i 
\in 
P, q_i\in 
Q$
and any $p\in P$,
$q\in Q$, we have 
$$p=\sum_i \eta(p, q_i) \,p_i,\,\,\,\,\, q= \sum_i \eta (p_i,q) \,q_i. $$
 Moreover, for any 
$K$-homomorphism $f:P\to P$, we have $\Tr(f)= \sum_i 
\eta(f(p_i),q_i)$.}

\skipaline {\sl Proof.}  Denote by $\rho$ the homomorphism $Q\to P^*$ 
defined 
by 
$\rho (q) 
(p)=\eta(p,q)$ for any $p\in P, q\in Q$.  Since
$\eta$ is non-degenerate, $\rho$ is an isomorphism.

To any element $b=\sum_i p_i\otimes q_i$   of $P\otimes Q$ we associate   
a 
homomorphism 
$\nu_b:
P^*\to Q$ sending each $x\in P^*$ into $\nu_b(x)=\sum_i x(p_i) q_i\in Q$.  
Since 
$P$
and $Q$ are projective modules, the formula $b\mapsto \nu_b$ defines an
isomorphism $P\otimes Q= \Hom (P^*,Q)$.

Observe that the identity $p=\sum_i \eta(p, q_i) p_i$ for all $p\in P$ is 
equivalent to 
$\rho\,
\nu_b=\id_{P^*}$.  Indeed,  $\rho\, \nu_b$ sends any 
$x\in 
P^*$
into the homomorphism $P\to K$ which maps each $p\in P$ into $\eta(p, 
\sum_i
x(p_i) q_i)= x(\sum_i \eta (p,q_i) p_i)$.  The equality $\rho\, \nu_b(x)=x$ 
holds
for all $x\in P^*$ if and only if $\sum_i \eta (p,q_i) p_i=p$ for all $p\in P$.  
A similar argument shows that the identity
$q= \sum_i \eta (p_i,q) q_i$ for all $q\in Q$ is equivalent to   $\nu_b 
\, \rho=\id_Q$.
  Thus, the only element $b$ satisfying conditions of the lemma is
determined from $\nu_b= \rho^{-1}\in \Hom (P^*,Q)$.

The second claim of the lemma is  standard (see for instance [Tu, Lemma II.4.3.1]).

\skipaline \noindent {\bf 7.1.2. Proof of (7.1.d).  }    Using the 
 expression for the trace given in Lemma 
7.1.1, 
we obtain for any $\ell \in L_1$, $\alpha\in \pi$,
 $$\eta(\ell,1_L)=  \Tr(\mu_\ell\vert_{L_1} )
=   \Tr(\mu_\ell\vert_{L_\alpha} :{L_\alpha} \to 
{L_\alpha})
=  \sum_i\eta(\ell \, p_i^{\alpha}, q_i^{\alpha} )
= \eta(\ell , \sum_i p_i^{\alpha} q_i^{\alpha} ).$$
Since the form $\eta : L_1\otimes L_{1}\to K$ is non-degenerate,  
$\sum_i p_i^{\alpha} q_i^{\alpha} =1_L$.

\skipaline \noindent {\bf 7.1.3.  Lemma.  } {\sl For any $\alpha\in \pi$, $a\in 
L_{\alpha^{-1}}, b\in 
L_\alpha$ and any 
expansion (7.1.c) of $b_\alpha$, we have
$\sum_i \, \eta( a \, p_i^{\alpha}, 1_L) \,\eta (q_i^{\alpha} b,1_L)= \eta 
(ab,1_L)$.}

\skipaline {\sl Proof.}
$$\sum_i \eta( a\, p_i^{\alpha}, 1_L) \,\eta (q_i^{\alpha}b,1_L)= 
\sum_i \eta( a,p_i^{\alpha}) \,\eta (q_i^{\alpha},b )= 
\sum_i \eta( a,p_i^{\alpha}) \, \eta (b,q_i^{\alpha} )
$$
$$= 
\eta( a, \sum_i \eta (b,q_i^{\alpha} ) \, p_i^{\alpha}) = \eta(a,b)=\eta 
(ab,1_L).$$

\skipaline \noindent {\bf 7.2.  Maps to $K(\pi,1)$ as combinatorial 
systems.}  
It
is well known that the homotopy classes of maps from a CW-complex to the
Eilenberg-MacLane space $X=K(\pi,1)$ admit a   combinatorial description 
in
terms of   elements of $\pi$ assigned to the 1-cells. These elements 
may 
be 
viewed as the holonomies of the corresponding principal $\pi$-bundles. Here 
we 
recall
this description   adapting it to our setting.

Let $T$ be a   CW-complex with underlying topological space $\tilde T$.  
By  vertices, edges and faces of $T$ we   mean  0-cells,   1-cells and   
2-cells of $T$, respectively. 
Denote the set of vertices of $T$ by $\Vert (T)$. 
Each oriented 
edge $e$ of $T$ leads from an initial vertex, $i_e \in \Vert (T)$, to a 
terminal vertex, 
$t_e\in \Vert (T)$ (they may coincide).
Denote by $\Edg(T)$ the set of oriented 
edges
of $T$.  Inverting the orientation we obtain an involution 
$e\mapsto e^{-1}$ on $\Edg(T)$.

Each face $\Delta$ of $T$ is obtained by adjoining a 2-disc to the 1-skeleton 
$T^1$ of $T$
along a map $f_\Delta:S^1\to T^1$. In general this map may be rather wild;    
we shall require the following property of regularity.
We say that the CW-complex $T$ is {\it 
regular}
if for any face $\Delta$ of $T$
the pre-image $f_{\Delta}^{-1} (\Vert (T))\subset S^1$ of $\Vert (T)$ is  a 
finite non-empty set which   
splits the 
circle $S^1$  into arcs   mapped by $f_{\Delta}$ 
homeomorphically onto
certain edges of $T$. These arcs in $S^1$   are called the {\it sides} of 
$\Delta$.  The image 
of each side of $\Delta$ under $f_\Delta$ is an edge of $T$. 
An orientation of $S^1$ determines an orientation and a cyclic order of  the 
sides of $\Delta$. The corresponding cyclically 
ordered oriented edges of $T$ form the boundary of $\Delta$.
Examples of  regular CW-complexes are provided by   triangulated spaces.

Assume   that $T$ is regular. 
Assume also that  a certain subset (possibly empty) of $\Vert (T)$ is 
distinguished;
the vertices belonging to this subset are called  the {\it base vertices}  of 
$T$.
By a {\it $\pi$-system} 
$\{g_e\}$ on $T$ we mean a function which
assigns to every   $e\in \Edg(T)$ an element $g_e\in \pi$ such 
that

(i) for any $e\in \Edg(T)$ we have $g_{e^{-1}}=(g_e)^{-1}$;

(ii) if ordered oriented edges $e_1,e_2,...,e_n$ of $T$ with $n\geq 1$ form the 
boundary 
of a
face   then $g_{e_1} g_{e_2} ...
g_{e_n}=1$.

Let $F_T$ be the set  of maps $ \Vert (T)\to \pi$ 
taking value $1\in \pi$ on    all   base vertices of $T$. 
The set $F_T$
is a group with pointwise multiplication.  
If $g=\{g_e\}$ is a $\pi$-system on $T$ and $\gamma\in F_T$ then the 
formula 
$(\gamma g)_e=\gamma (i_e) g_e (\gamma (t_e))^{-1}$ yields a new 
$\pi$-system 
$\gamma g$ on $T$. This defines a left action 
of 
$F_T$ on 
the set of $\pi$-systems on $T$. We say that two $\pi$-systems $g,g'$ on 
$T$ 
are 
{\it homotopic} if   they lie in the same orbit of this action, i.e., if there 
is 
$\gamma\in 
F_T$ such that $g'=\gamma g$.
 It is clear that homotopy is 
an 
equivalence 
relation.

Every $\pi$-system $g $ on $T$ gives rise to a map $\tilde g:\tilde 
T 
\to 
X=K(\pi,1)$
which sends all vertices of $T$ into the base point of $X $ and 
sends 
each oriented edge
$e\in \Edg(T)$ into a loop in $X$ representing  
$g_e\in \pi$.  
An
elementary obstruction theory shows that the formula $g \mapsto \tilde g$   
establishes a 
bijective
correspondence between the homotopy classes of   $\pi$-systems on $T$ and
the pointed homotopy classes  of maps $\tilde T\to X$. (These are classes 
of 
maps 
  sending  all the base vertices of $T$
in the base point of $X$   modulo homotopies constant on the base 
vertices).

By a CW-decomposition of an $X$-cobordism (or an  $X$-manifold) $W$
  we shall mean a regular 
CW-decomposition   such that $\partial W$ is a subcomplex and all the base 
points are among vertices.
  A $\pi$-system $g$ on    a CW-decomposition $T$  of $W$ is said to be {\it 
characteristic} if 
the map $\tilde g:W=\tilde T\to X$ is homotopic (in the pointed category) 
to the 
given
characteristic map $W\to X$. Note that the characteristic $\pi$-systems on 
$T$ 
form 
a single homotopy class of $\pi$-systems.

\skipaline \noindent {\bf 7.3.  State sums on closed $X$-surfaces.} Fix a 
biangular $\pi$-algebra $L$.  Let $W$ be   
 a closed   $X$-surface, i.e., a closed  oriented surface 
 endowed with a  map  $W\to X=K(\pi,1)$. (We do not provide 
closed 
surfaces 
with base points, cf. Section 1.1). Here we
define a state sum invariant  $\tau(W)=\tau_L(W) \in K$. 

Choose a 
(regular) CW-decomposition
$T$ of $W$. 
By a {\it flag} in $T$ we   mean  a  pair 
(a 
face  $\Delta$ of $T$, a side $e$ of $\Delta$). The flag
  $(\Delta,e)$   determines an orientation of $e$ such that
$\Delta$ lies on the right of $e$. This means that 
the pair (a  vector looking from a point of $e$ inside $\Delta$, the 
oriented side   $e$) is positively oriented with respect to the given 
orientation of $W$. By abuse of notation, we shall denote the   edge of 
$T$ corresponding to $e$ by the same letter $e$.

Let  $g $ be a characteristic $\pi$-system on $T$.
 With each flag 
$(\Delta,e) $ in $T$  we associate the $K$-module $L(\Delta,e,g)=L_{g_e}$
where the orientation of $e$ is determined by $\Delta$ as above.
Each unoriented edge $e$ of $T$ appears in two flags, 
$(\Delta,e) , (\Delta',e) $, and inherits from them opposite orientations. 
Since the corresponding values of $g$ are inverse to each other, we have the  
canonical vector $b_{g_e}\in 
L(\Delta,e,g)  \otimes L(\Delta',e,g)$ introduced in Section 7.1. The 
tensor 
product 
of these vectors over 
all unoriented
edges of $T$  is an element,   $B_g\in \otimes_{(\Delta,e) }
L(\Delta,e,g) $ where $(\Delta,e) $ runs over all flags in $T$.

Let $\Delta$ be   a   face of $T$ with $n\geq 1$ oriented ordered sides    
$e_1,...,e_n$ whose
orientation is determined by    $\Delta$ as above 
and   the  order is chosen  so that the terminal endpoint of $e_{r}$ is 
  the  initial 
endpoint 
of 
$e_{r+1}$  for $r=1,...,n \,(\mod \,n)$. Note that 
$g_{e_1}... g_{e_n}=1$. The $n$-linear form
$$ L(\Delta,e_1,g) \otimes L(\Delta,e_2,g)\otimes ... \otimes L(\Delta,e_n,g) 
\to K$$
defined by $(a_1,a_2,...,a_n)\mapsto \eta(a_1 a_2 ... a_n,1_L)$ 
with 
$a_r\in L(\Delta,e_r,g) $ for $r=1,...,n$ 
is invariant under cyclic permutations
and therefore is determined  by $\Delta$.
The tensor product of these  multilinear forms over all  faces of $T$
is a homomorphism, $D_g: \otimes_{(\Delta,e) }
L(\Delta,e,g) \to K$ where $(\Delta,e) $ runs over all flags in $T$. Set 
$\langle g \rangle=D_g(B_g)\in K$.

\skipaline \noindent {\bf 7.3.1. Claim.  } {\sl The state sum $\langle g 
\rangle$ does not depend on the choice of 
  $g$ in its homotopy class and does not depend on the choice of $T$.}

\skipaline

 For a    proof, see Section 8.
 
Claim 7.3.1 implies 
that 
$\tau(W)=\langle g \rangle=D_g(B_g)\in K $ is a well defined 
invariant of the 
$X$-surface $W$.
By the very definition, $\tau(W) $ depends only on the homotopy class
of the   characteristic map $W\to X$ and is multiplicative with respect 
to 
disjoint union 
of closed  $X$-surfaces.

For connected $W$, we can give an explicit formula for $\tau(W)$.
 We need the following notation. For  
$b\in L_{\alpha}\otimes L_{\alpha^{-1}}, b'\in L_{\beta}\otimes L_{\beta^{-1}}$ 
with $\alpha, \beta \in \pi$, we define  
$[b,b']\in L_{\alpha\beta \alpha^{-1} \beta^{-1}}$ as follows:   expand  $b,b'$ 
into finite sums 
$b=\sum_i p_i\otimes q_i, b'=\sum_j p'_j\otimes q'_j$ with $p_i \in  
L_{\alpha}, q_i\in  L_{\alpha^{-1}}, p'_j\in L_\beta, q'_j \in 
L_{\beta^{-1}}$ and    set
$$ [b,b']=\sum_{i,j} p_i p'_j q_i q'_j \in L_{\alpha\beta \alpha^{-1} 
\beta^{-1}}.$$
It is clear that $ [b,b']$   does not depend
on the choice of the expansions of $b,b'$.  Assume now that $W$ is a closed 
connected oriented $X$-surface of genus $n\geq 1$.  Consider the 
CW-decomposition of 
$W$
formed by one vertex, $2n$ oriented loops     $a_1, a_2,...,a_{2n}$ and one 
face
with boundary $\prod_{r=1}^n (a_{2r-1} a_{2r} a_{2r-1}^{-1} a_{2r}^{-1})$.
We choose this CW-decomposition so that the  intersection number $[a_1] 
\cdot  [a_2]$   of the 
integer homology classes represented by the loops $a_1, a_2$ equals 
$-1$.
We deform the characteristic map $g:W\to X$ so that is sends the only vertex of 
$W$ into the base point of $X$. Consider the induced homomorphism $g_{\#}$ 
of the fundamental groups and set $\alpha_s=g_{\#} (a_s)\in \pi$ for 
$s=1,...,2n$. Recall   the canonical element
$b_{\alpha_s}\in  L_{\alpha_s}\otimes L_{\alpha_s^{-1}}$.
Computing $\tau (W)$ from  this CW-decomposition of $W$, we obtain
$$\tau(W)=  \eta(\prod_{r=1}^n [b_{\alpha_{2r-1}},
b_{\alpha_{2r}}], 1_L)$$
where $\prod_{r=1}^n [b_{\alpha_{2r-1}},
b_{\alpha_{2r}}]\in L_1$. 
In the case $W=S^2$, all maps $W\to X$ are homotopic to the constant map. The 
invariant $\tau(S^2)$ can be computed from the CW-decomposition of $S^2$
consisting of one vertex, one loop   and two faces (hemispheres).
Fix  an expansion of the canonical element $b_1=\sum_i p_i\otimes q_i\in 
L_1\otimes L_1$.
We have
$$\tau(S^2)=\sum_i \eta(p_i ,1_L)\, \eta(q_i, 1_L)=\eta(\,\sum_i \eta (q_i,1_L) 
\,
p_i,1_L) =\eta (1_L,1_L)=\Dim L_1.$$

  \skipaline \noindent {\bf 7.4.  A lattice 
$(1+1)$-dimensional  HQFT.} We extend the invariant $\tau$ of closed 
$X$-surfaces constructed above to a $(1+1)$-dimensional  HQFT with target 
$X=K(\pi,1)$.
The construction goes in three steps. At the first step we assign  
$K$-modules 
to 
 so-called 
trivialized 1-dimensional  $X$-manifolds. At the second step we extend 
this 
assignement to   a 
preliminary
$(1+1)$-dimensional HQFT defined for  
trivialized $X$-manifolds and $X$-cobordisms with trivialized boundary. 
   At the third step we  get rid of the trivializations. 

\skipaline  {\bf   Step 1.  } We say that  a 
 1-dimensional $X$-manifold $M$ is {\it trivialized} if it is provided with a 
CW-decomposition $T$ (such that the base points of the components of $M$ are 
among 
the vertices) and a characteristic $\pi$-system $g$  on $T$.
  We provide each edge $e$ of $M$ with canonical orientation induced by 
the 
one of $M$.
We associate with the trivialized  1-dimensional $X$-manifold $M=(M,T,g)$ 
the $K$-module
 $A(M) =\otimes_e {L_{g_e}}$
 where $e$ runs over all canonically oriented edges of 
$T$.
It is clear that the module $A(M)$ is projective of finite type.
For disjoint  trivialized 1-dimensional 
$X$-manifolds $M,N$, we have $A(M\coprod N)=A(M) \otimes_K A(N)$.
If $M=\emptyset$ then by definition $M$ is trivialized and   $A(M)=K$.

\skipaline  {\bf   Step 2.  }
Consider   a 2-dimensional  $X$-cobordism
$(W,M_0,M_1)$ as defined in Section 1.1. Assume that the bases 
$M_0,M_1$ of $W$ are trivialized. We define a  $K$-homomorphism 
$\tau(W):A(M_0) \to A(M_1)$ as follows.
Choose a regular CW-decomposition $T$ 
of $W$ extending the given CW-decomposition of $\partial 
W=(-M_0)\cup M_1$. Choose a characteristic
$\pi$-system $g$ on $T$  extending the given $\pi$-system   on 
$\partial 
W$.
Applying the constructions of Section 7.3 to $ g$ 
 we obtain 
 a vector $B_{ g}\in \otimes_{(\Delta,e) }
L(\Delta,e,  g) $ where $(\Delta,e) $ runs over flags in $T$  such that 
$e$ does not lie on   $\partial W$. As  in Section 7.3, we obtain a 
homomorphism $
D_{ g}: \otimes_{(\Delta,e) }
L(\Delta,e, g) \to K$ where $(\Delta,e) $ runs over all flags in $T$.
Contracting these two tensors we obtain
a homomorphism
$\langle g \rangle:  \otimes_{(\Delta,e \subset \partial W)} 
L(\Delta,e,   g) \to K$. We provide each edge $e$ lying on $M_r\, 
(r=0,1)$ with 
the orientation induced by the one of $M_r$.
This orientation of $e$ coincides with the one induced by the only face 
attached 
to   
$e$ 
 if $e\subset M_0$ and is opposite to it if $e\subset M_1$.
Therefore  $L(\Delta,e,  g)=L_{g_e}$ if $e\subset M_0$ 
  and $L(\Delta,e,  g)=L_{(g_e)^{-1}}$ if $e\subset M_1$. We   
identify $L_{(g_e)^{-1}}$ with the dual of
 $ L_{g_e}$ using the inner product $\eta$. In this way we can view  
$\langle g \rangle$
as a  homomorphism
$$ A(M_0)=\bigotimes_{e \subset M_0}
 L_{g_e} \to \bigotimes_{e \subset M_1}  L_{ 
g_e}=A(M_1).$$
We claim that 
$\langle g \rangle $   does not depend on the choice of $g$  and
$T$; this will be checked in Section 8.
 Set  $\tau(W)=\langle g \rangle:A(M_0)  \to A(M_1)$.
By the very definition, $\tau(W) $ depends only on the homotopy class
of the   characteristic map $W\to X$ and satisfies axiom (1.2.5). The next 
lemma 
describes the behavior of $\tau$ under the 
gluing of $X$-cobordisms.

\skipaline \noindent {\bf 7.4.1.  Lemma.  } {\sl Let $M_0,M_1,N$ be 
trivialized 1-dimensional $X$-manifolds (possibly void). If
a 2-dimensional  $X$-cobordism
$(W,M_0,M_1)$ is obtained from two  2-dimensional  $X$-cobordisms
$(W_0,M_0, N)$ and $(W_1,N, M_1)$ by gluing along the identity map $N= N$
then $\tau(W)= \tau(W_1)\circ \tau(W_0):A(M_0) \to A(M_1)$.}

\skipaline {\sl Proof.}    For $r=0,1$, fix a regular CW-decomposition $T_r$ of 
$W_r$ 
extending 
the given 
CW-decomposition of the boundary.
 Gluing $T_0$ and $T_1$ along $N$ we obtain a regular CW-decomposition, $T$, of 
$W$.
 The lemma is a reformulation of the following claim:
 
 \skipaline
 
 $(\ast)$ {\it Let for $r=0,1$, $g_r$ be a  $\pi$-system on $T_r$ extending  
the 
  given $\pi$-systems on $M_r$ and $ N$. Let $g$ be 
the  
unique $\pi$-system on $T$ extending   $g_0$ and $g_1$. 
  Then 
 the homomorphism
 $\langle g \rangle :A(M_0) \to A(M_1)$
 is the composition of 
 $\langle g_0 \rangle :A(M_0 ) \to A(N)$
 and
  $\langle g_1 \rangle: A(N) \to A(M_1) $.}
  
 \skipaline
  
  This claim directly follows from the definition of  
$\langle g\rangle$.

 \skipaline  {\bf   Step 3.  } The constructions above can 
be 
summarized as follows. 
 To any trivialized 1-dimensional $X$-manifold $M$  we assign
a module $A(M) $. 
This module is natural with respect to $X$-homeomorphisms preserving the 
trivializations. To any $X$-cobordism $(W,M_0,M_1)$
between trivialized  1-dimensional $X$-manifolds we assign a
 homomorphism 
$\tau(W):A(M_0)\to A(M_1)$.
This data looks like a HQFT and  satisfies the natural versions of the 
axioms 
(1.2.1) - (1.2.6) and (1.2.8). However, in general $\tau(M\times [0,1])\neq 
\id_{A(M)}$.  
There is a standard procedure which allows to 
pass from such a  
pseudo-HQFT   to a genuine $(1+1)$-dimensional HQFT.
This procedure is described in detail in a   similar setting in [Tu,Section VII.3]
cf. also Section 2.1.
The idea is that if  
$t_1,t_2$ are two trivializations of a  1-dimensional $X$-manifold $M$
then the cylinder $W=M\times [0,1]$ (mapped to $X$ via the composition of 
the projection $M\times [0,1]\to M $ with the 
characteristic map $M\to X$) is an 
$X$-cobordism
between $(M, t_1)$ and $(M, t_2)$.
This gives a   homomorphism $p(t_1,t_2)=\tau(W):A (M, t_1)\to 
A(M, 
t_2)$.
 By Lemma 7.4.1, we have the 
identity
$p(t_1,t_3)=p(t_2,t_3)\, p(t_1,t_2)$.   Taking $t_1=t_2=t_3$ we obtain that 
$p(t_1,t_1)$ is a projection onto a direct summand, $\hat {A}(M,t_1)$,  of 
$A(M,t_1)$. 
Moreover, 
  $p(t_1,t_2)$ maps  
$\hat {A}(M,t_1)$ isomorphically onto $\hat {A}(M,t_2)$. This allows us 
to 
identify the modules $\{\hat {A}(M,t))\}_{t}$ where $t$ runs over all 
trivializations 
of $M$ along these canonical 
isomorphisms and to obtain a module, $\hat {A}(M)$, independent of 
  $t$. Next we observe that for any 2-dimensional $X$-cobordism
$(W,M_0,M_1)$ with trivialized bases, the homomorphism
$\tau(W):A(M_0)\to A(M_1)$ maps $\hat {A} (M_0)\subset A(M_0)$ into 
$\hat {A} (M_1)\subset 
A(M_1)$. This yields a homomorphism $\hat {\tau} (W):\hat {A}(M_0) \to \hat {A} 
(M_1)$ 
independent 
of 
the trivializations of the bases. 
The modules $\hat {A} $ and the homomorphisms $\hat {\tau}  $ form the
  $(1+1)$-dimensional  HQFT  with target 
$X$
associated with   $L$. Note that if $W$ is a closed $X$-surface, then $\hat 
{\tau} (W)=\tau(W)$.

 \skipaline \noindent {\bf 7.5. Remark.} The HQFT  constructed above is  
additive (resp. multiplicative) with respect to direct sums (resp. tensor 
products) of biangular 
$\pi$-algebras.
This HQFT is functorial with respect to   pull-backs of   biangular group-algebras
 along group homomorphisms.
 Note that any  group homomorphism  
 $q:\pi' \to \pi$  induces a map $X'=K(\pi',1)\to K(\pi,1)=X$ 
denoted  $f_q$.
It is easy to deduce from definitions that if
a biangular $\pi'$-algebra  $ L'= q^*(L)$ is the pull-back   
of a 
 biangular $\pi$-algebra $L$ along $q$, then the lattice HQFT  
determined by $L'$
 is obtained from the lattice HQFT   determined by $L$ via 
pulling back along $f_q$, cf.  Section 2.2.
In particular, each closed $X'$-surface $W=(W, 
g:W\to X')$ gives rise to a closed $X$-surface 
 $W_q=(W, f_q\, g:W\to X)$ and 
 $\tau_{L'} (W)=\tau_L (W_q)$. For instance, consider the homomorphism $q:\pi 
\to \{1\}$ where $\{1\}$ is the trivial group. If $L$ is the 1-dimensional 
$\{1\}$-algebra $K$, then $q^*(L)$ is the group ring $K[\pi]$
 and $\tau_{K[\pi]} (W)=\tau_L (W_q)=1\in K$.
 
  The relations with push-fowards of biangular algebras are more subtle, cf. 
Remark 9.7.

\skipaline \centerline {\bf 8.  Lattice models on  skeletons}

\skipaline \noindent {\bf 8.1.  Skeletons of 2-manifolds.} We shall use 
the 
language of graphs.  By a graph we mean a 1-dimensional CW-complex. 
The {\it valency} of a vertex of a graph is the number of edges incident to 
this vertex (counted with multiplicity).

Let $W$ be a compact oriented surface with based (possibly void) boundary.
  A {\it    skeleton}    of 
$W$ is a 
finite graph $\Gamma$ embedded in $W$ such that

(i) $\Gamma$ has no isolated vertices (i.e., vertices of valency $0$);

(ii) each
connected component of $W \backslash \Gamma$ is  either an open 2-disc
or  a half-open 2-disc homeomorphic to $[0,1) \times (0,1)$
and meeting $\partial W$ along the arc $0\times (0,1)$.

(iii)    each point of the set   $\Gamma \cap \partial W$ is a 1-valent vertex 
of $\Gamma$;   the base points of $\partial W$ do not belong to $\Gamma$.

The 1-valent vertices of $\Gamma$ lying on $\partial W$ are called {\it feet} of 
$\Gamma$. An 
open 
subset of  $\Gamma$   consisting of  
a 
foot and the unique open edge  incident to it is called 
a 
{\it leg} of  $\Gamma$.

Conditions (i) - (iii) imply that  each component of $\partial W$ contains at 
least 
one foot 
of $\Gamma$ and  that the part of $\Gamma$ lying in a connected component of 
$W$ is connected. 
 Note that a skeleton may have loops (i.e., edges with coinciding endpoints) 
and multiple edges
(i.e., different edges with the same endpoints). It may also have 2-valent 
vertices and 1-valent 
vertices lying in $\Int(W)$.

There are four basic local moves on  the skeletons   of $W$. 
These moves transform  a skeleton  
of $W$ into another  skeleton  
of $W$
with the same feet.
The first  move, 
called  the {\it contraction} move,
contracts  an edge $e$ into a point provided       $e$ is 
neither a loop nor a leg. 
 The second move, called  the {\it  biangular} move, introduces a small biangle 
(or 
bigon) in the middle of an edge of the skeleton. More precisely, this 
move  
replaces 
a  small  subarc of the edge with 
two parallel embedded subarcs with the same endpoints disjoint from the 
rest of 
the 
skeleton.
  The number of vertices 
decreases by 1 under the contraction move and increases    by 2
 under the   biangular move.   The third and forth moves are the inverses of 
the 
contraction and biangular moves.

 Another useful move on a skeleton $\Gamma$ adds a small loop based at an edge 
of $\Gamma$.
 The loop should be disjoint from $\Gamma$ except at its endpoint and should 
bound a small disc in 
$W\backslash \Gamma$.
 This {\it loop move} is a composition of a biangular move and a contraction 
 along one of the sides of the biangle. Conversely, the biangular move
 is a composition of a loop move and an inverse contraction move.

\skipaline \noindent {\bf 8.1.1.  Lemma.  } {\sl Any two skeletons of $W$ 
having the same number of feet on 
every 
component of 
$ \partial W$
can 
be related by a finite sequence of  basic moves and ambient isotopy of $W$ 
constant on the base points of $\partial W$.}

\skipaline {\sl Proof.}  This lemma is a simple application of the  known 
description of the moves relating generalized spines of $W$.
By a {\it generalized spine} of $W$ we mean a finite  trivalent graph
$G$ embedded in $\Int(W)$ such that all 
connected components of $W \backslash G$ are either open 2-discs
or half-open annuli $S\times [0,1)$ where $S=S\times 0$ runs over 
components 
of 
$\partial W$. There are  three  local moves on the generalized spines.
The first  move, 
called  the $IH$-move,
replaces a piece of the spine looking like the letter $I$
by a piece looking like the letter $H$. 
 The second and third moves are  the biangular move considered above and its 
inverse. 
 It is well known that any two generalized spines of $W$ can be related 
by a finite sequence of such moves and ambient isotopy.

Let us call   a skeleton $\Gamma$ of $W$ {\it trivalent}
if all its vertices lying in $\Int(W)$ are trivalent and at least one of these 
vertices is not incident to a 
leg. It is easy 
to see that any skeleton of $W$  
can be transformed by the basic moves 
into a trivalent skeleton.
Removing from a  trivalent  skeleton of $W$ all its legs 
we 
obtain a generalized spine, $G$, of $W$. We can reconstruct the skeleton by 
adjoining 
to $G$ several legs lying in the annuli components of 
$W\backslash  
G$. Note that different ways to adjoin the legs are related by  
 $IH$-moves. 
Now the results mentioned above  imply that 
  any two trivalent skeletons of $W$ can be related by a sequence of 
$IH$-moves, 
biangular moves and  inverse biangular moves.
  It remains to observe that each $IH$-move is a composition of a contraction 
move and
  an inverse contraction move.

\skipaline \noindent {\bf 8.2.  Skeletons versus CW-decompositions.} Let $W$ be 
a compact 
oriented surface with based (possibly void) boundary.
A  regular CW-decomposition $T$ of $W$ such that $\partial W$ is a subcomplex and the base
points of $\partial  W$ are 
among the vertices 
yields a \lq\lq dual" skeleton $\Gamma_T$ 
of 
$W$ as follows. 
Choose inside every 2-cell and every 1-cell of 
$T$ a point called 
its 
center. Let us connect the center of each 2-cell $\Delta$ to the centers 
of its sides by embedded 
arcs   lying in $\Delta$ 
and disjoint except at their common endpoint.
A  1-cell $e$ of $T$ not lying in $\partial W$ gives rise to a \lq\lq 
dual edge" formed by two arcs joining the center of $e$ with the centers of  
  2-cells adjusted to $e$.
A  1-cell $e$ of $T$   lying in $\partial W$ gives rise to a \lq\lq dual 
edge" which is the arc joining the center of $e$ to the center of the only  
2-cell adjusted to $e$.  The  vertices of $\Gamma_T$ are 
  the  centers of   the 2-cells of $T$  and  the centers of the  1-cells of 
$T$ lying in $\partial W$. The edges of $\Gamma_T$ are the dual 
edges
of the 1-cells of $T$. Note that the legs of $\Gamma_T$ are    the dual 
edges
of the 1-cells of $T$ contained in $\partial W$.

The formula $T\mapsto \Gamma_T$ establishes a bijective correspondence between 
the (isotopy 
classes of) regular CW-decompositions as above and skeletons of $W$.
Here by isotopy we mean ambient isotopy constant on the base points. 
The inverse   construction of    $T=T_\Gamma$ from a   skeleton 
$\Gamma$ of $W$ is as follows. Choose in each component $U$ of $W\backslash 
\Gamma$
 a point called its 
center. We assume that if $U$ meets $\partial W$ along an arc then the center of 
$U$ is chosen on this arc.
Moreover, if  $U$ meets $\partial W$ along an arc containing a base point of $W$ 
then we take this point as the center of $U$. The  centers of the
components of $W\backslash 
\Gamma$ are the vertices of 
$T$. Each   edge $e$ of $\Gamma$   
 gives rise to a  1-cell $e^*$ of $T$ which 
 crosses $e$ transversally in one point  and connects the centers of the
 components of $W\backslash \Gamma$ adjusted to $e$.
Moreover, if $e$ is a  leg   of $\Gamma$ then $e^*\subset 
\partial W$ and $e^*\cap e$ is  the foot of $e$.
This completes the description of the 0-cells  and 1-cells of $T$,  the 
connected components of their complement in $W$ are the 2-cells of $T$.
The 2-cells  of $T$ bijectively correspond to those vertices of 
$\Gamma$ which are not feet. 

Observe that if an edge  $e$ 
   of $\Gamma$  is oriented then 
the dual 1-cell $e^*$ of $T_\Gamma$ can be oriented so that it  
 crosses $e$     from    right to left,
 the  right and  left of $e$ being determined by the orientations of $W$ 
and 
$e$, 
cf. Section 7.3.    Obviously, 
$(e^{-1})^*=(e^*)^{-1}$. This establishes a bijection 
$\Edg(\Gamma)=\Edg(T_\Gamma)$ equivariant with respect to the involution 
$e\mapsto e^{-1}$.

 \skipaline \noindent {\bf 8.3.  $\pi$-systems on skeletons.}  
  Let $\Gamma$ be a  skeleton of a compact oriented surface $W$ with based 
(possibly void) 
boundary. 
A {\it $\pi$-system $\{g_e\}$ } on 
$\Gamma$ 
is  a function which
assigns to any oriented edge $e\in \Edg(\Gamma)$ an element $g_e\in \pi$ such 
that

(i) for any $e\in \Edg(\Gamma)$, we have $g_{e^{-1}}=(g_e)^{-1}$;

(ii) for any vertex $u$ of $\Gamma$ of valency $n\geq 1$,  the oriented edges   
$e_1, e_2, ...,e_n$ of $\Gamma$ with  terminal 
vertex $u$
  and the cyclic order   opposite to 
the 
one induced by the orientation of $W$ in a neighborhood of $u$ satisfy
$g_{e_1} g_{e_2} 
... g_{e_n}=1$.

Two $\pi$-systems $g,g'$ on $\Gamma$ are said to be {\it homotopic }    
if there is
a function $\gamma:\pi_0(W\backslash \Gamma)\to \pi$
such that $\gamma (U)=1$ if $U$ is a component of $W\backslash \Gamma$ 
containing a base 
point of $\partial W$, and
for any   $e\in \Edg(\Gamma)$,   we have $ 
g'_e=\gamma (U) \,g_e \, (\gamma(V))^{-1}$ where $U,V$ are the components of 
$W\backslash 
\Gamma$ lying on the right and   left of $e$, respectively (possibly, 
$U=V$).
If $\gamma$ takes the value $1\in \pi$ on all components of $W\backslash 
\Gamma$
except one component $U$ then we say that $g'$ is obtained from $g$ by a 
{\it homotopy move} at 
$U$. It is clear that two $\pi$-systems on $\Gamma$ are homotopic if and 
only 
if they  
can be related by homotopy moves.

Setting $g_{e^*}= g_e$ we obtain a bijective correspondence between the 
$\pi$-systems on $\Gamma$ and on  the dual CW-decomposition 
$T_\Gamma$. This correspondence preserves the relation of 
homotopy.

Every $\pi$-system $\{g_e\}$ on $\Gamma$ gives rise to a map $\tilde g:W  
\to X$
which sends all vertices of  $T_\Gamma$  
into 
the base point of $K(\pi,1)$ and for each $e\in \Edg(\Gamma)$ sends 
the
oriented 1-cell $e^*$   into a loop representing   $g_e$. Condition 
(ii) above  
ensures that this mapping of the 1-skeleton of $T_\Gamma$ extends to $W$. 
   This establishes a bijective
correspondence between the homotopy classes of   $\pi$-systems on 
$\Gamma$ 
and
the set   of (pointed) homotopy classes of maps $W\to X$.     

Each basic move $\Gamma \mapsto \Gamma'$ on  skeletons lifts to   
$\pi$-systems by  
  transforming any $\pi$-system $g$ on $\Gamma$ 
into a
$\pi$-system $g'$ on $\Gamma'$ coinciding with $g$ on the common part     
of 
$\Gamma$ and 
$\Gamma'$. The system $g'$
 is uniquely determined by $g$ in the case of the contraction move,
 the inverse contraction move and the inverse 
biangular move. 
Under the biangular move, such a $\pi$-system $g'$   exists  and is unique up 
to homotopy moves at the small  biangle component of $W\backslash \Gamma'$ 
created by the move.  For all basic 
moves,
 the mapping $g\mapsto g'$   
establishes a bijection between the homotopy classes
of $\pi$-systems on $\Gamma$ and $\Gamma'$. The classes corresponding to 
each 
other under 
this bijection determine the same homotopy class of maps  $W\to X$.

\skipaline \noindent {\bf 8.4.  State sums on closed $X$-surfaces
via skeletons.}
Fix a biangular $\pi$-algebra $L$. 
Let $W$ be   
 a closed   $X$-surface. Here we
define a numerical invariant  $\tau(W)=\tau_L(W) \in K$ using the 
skeletons of $W$.

 Take    a skeleton 
$\Gamma$ of $W$ and represent the characteristic map 
$W\to X$
by a $\pi$-system $g=\{g_e\}_e$ on $\Gamma$.  Every unoriented edge
of $\Gamma$    gives rise to two oriented edges $e, e^{-1}$ and 
to the vector
$b_{g_e}\in L_{g_e} \otimes L_{(g_e)^{-1}}$.   The tensor product of these 
vectors 
over 
all unoriented edges
of $\Gamma$
is an element,   $B_g\in  \otimes_{e\in \Edg (\Gamma)} 
L_{g_e} $.

Consider a  vertex $u$ of $\Gamma$ of valency $n\geq 1$ and   $n$ oriented 
edges 
$e_1,...,e_n$ of 
$\Gamma$ incident to $u$ and directed towards $u$. We 
  choose the cyclic order $e_1,...,e_n$ so that it is opposite to the one 
induced 
by the orientation of $W$ in a neighborhood of $u$.
 By (8.3.ii), 
$g_{e_1} ... g_{e_n}=1$. The $n$-linear form
$$ L_{g_{e_1}} \otimes ... \otimes L_{g_{e_n}} \to K$$
defined by $(a_1,...,a_n)\mapsto \eta(a_1...a_n,1_L)$ 
with 
$a_r\in L_{g_{e_r}} $ for $r=1,...,n$
is invariant under cyclic permutations
and therefore is uniquely determined by $u$.
The tensor product of these   forms over all vertices of $\Gamma$
yields a homomorphism, $D_g: \otimes_{e\in \Edg (\Gamma)} 
L_{g_e} \to K$. Set 
$\tau(W)=\langle g \rangle= 
D_g(B_g)\in K$. It is clear that for the skeleton dual to a CW-decomposition of 
$W$ 
this definition is equivalent to the one given in Section 7.3. 

We claim that   $\langle g \rangle$ does not depend on the choice of $g$ 
and $\Gamma$. This implies Claim 7.3.1.
It suffices to check that $\langle g \rangle$ is preserved under the basic moves 
and the
homotopy moves.  It follows 
directly from definitions 
and (7.1.d) that   $\langle g \rangle$ is invariant  under the loop move. The 
invariance of $\langle g 
\rangle$  under the 
contraction move follows 
from definitions and Lemma 7.1.3. This implies the 
invariance of $\langle g \rangle$  under  the biangle move which is a 
composition of a loop move and an 
inverse contraction 
move.

 Assume that two $\pi$-systems $g_1, g_2$ on $\Gamma$ are related by a 
homotopy move at a component $U$ of $W\backslash 
\Gamma$. 
We claim that it is possible to relate $g_1$ to $g_2$ by  
basic moves. This will imply that   $\langle g_1 \rangle=\langle g_2 
\rangle$.
To prove our claim, take a small subarc $f$ on an edge of $\Gamma$ adjacent to 
$U$. Let $A,B$ be the endpoints of   
$f$.  Let us connect the points $A,B$ by a small arc,
   $f'$, lying inside $U$ (except at its endpoints $A,B$). The arc $f'$ splits 
the 
disc $U$ into a biangle bounded by $f\cup f'$  and a complementary 
open 2-disc denoted $D$.
It is clear that $\Gamma'=\Gamma \cup f'$ is a  skeleton of
$W$ obtained from $\Gamma$ by the biangular move. This move transforms  
 $g_1,g_2$
into certain $\pi$-systems $g'_1,g'_2$ on  $\Gamma'$ which can be chosen to be  
related by a homotopy move  at the component 
$ D$ of 
$W\backslash 
\Gamma'$.
 Now we transform  $\Gamma'$ as follows.  Let us deform $f'$ inside $U$ keeping 
the endpoint $A$ and gradually 
  pushing the endpoint $B$ along the edges 
adjacent to $U$. We do this until   $B$  traverses the whole boundary of 
$U$ and comes back to the original edge   from the other side
 of $A$.
 At the end of this deformation we obtain a  skeleton $\Gamma''$
 isotopic to $\Gamma'$.
 Note that while   $B$ moves along an edge,  both $\pi$-systems $g'_1,g'_2$ move 
along in the obvious way. When 
$B$ moves across a vertex of $\Gamma'$ adjacent to $U$,
the  skeleton under deformation is transformed via a contraction move and an  
inverse 
contraction move. 
Under these two moves, the $\pi$-systems $g'_1,g'_2$  are transformed in a 
canonical 
way     remaining related by a homotopy move at the image of 
$  D$. At the end of the deformation we obtain $\pi$-systems, 
$g''_1,g''_2$ on $\Gamma''$  related by a homotopy move at the image of 
$  D$ under the deformation. This image is a biangle. Applying the inverse 
biangular move we transform $g''_1,g''_2$ 
into one and the same $\pi$-system on $\Gamma$. This relates 
$g_1$ to $g_2$ by a sequence of basic moves.

  \skipaline \noindent {\bf 8.5.  Construction of a 
$(1+1)$-dimensional  HQFT via skeletons.}
The construction follows along the same lines as   in Section 7.4.

 \skipaline  {\bf   Step 1.  }  We assign  $K$-modules 
to 
the so-called 
split 1-dimensional  $X$-manifolds.   A
connected  non-empty 1-dimensional $X$-manifold $M$ is  {\it split }  if 
it is provided 
with a finite   set of points $x_1,...,x_n$ with $n\geq 1$ and elements
$h_1,...,h_n\in \pi$  such that: the points $x_1,...,x_n$ are 
distinct from 
each other and from the base point;
starting from the base point and moving in the direction given by the 
orientation of $M$ we meet consequtively   $x_1,...,x_n$;
the product $h_1 ...h_n\in \pi$ is the homotopy class of the loop 
represented by the characteristic map $M\to X= K(\pi,1)$. 
We  call $x_1,...,x_n$ {\it split points} and think of each  
$h_r$ as being attached to   $x_r$.
We assign to $M$ the module $A(M) =\otimes_{r=1}^n {L_{h_r}}$.
A non-connected 1-dimensional  $X$-manifold $M$ is {\it split} if its 
connected components are split. We define $A(M)$ as the tensor product of 
the modules $A$ corresponding to the connected components.
If $M=\emptyset$ then by definition $M$ is split and   $A(M)=K$.

Note that each trivialization $(T,g)$ of a 1-dimensional  $X$-manifold $M$ 
in the sense of Section 7.4 gives rise to  a  dual splitting of 
$M$. Namely, as split points we take  centers of the edges of $T$. With  
the center of a (canonically) oriented edge $e$ we associate the 
element $g_e\in \pi$. 

 \skipaline  {\bf  Step 2.  }  Consider   a 2-dimensional  
$X$-cobordism
$(W,M_0,M_1)$ with split bases 
$M_0,M_1$. 
Let $\Gamma$ be a  skeleton  
of $W$   and $g$ a $\pi$-system on $\Gamma$. We say that the pair 
$(\Gamma,g)$  {\it extends the splitting} of $\partial W$  if 

(i) the set of   feet of $\Gamma$ coincides with the set of split 
points of $\partial W$;

(ii) for each split point $x\in M_0$ (resp. $x\in M_1$) the value of $g$ 
on 
the edge of $\Gamma$ incident to $x$ and directed inside (resp. outside) 
$W$ is equal to the element of $\pi$ attached to $x$.

 Choose a  pair $(\Gamma,g)$  extending the splitting of $\partial W$
and such that
  the map $\tilde g:W\to X$ is homotopic to the characteristic map $W\to X$.
Applying the constructions of Section 8.4 to $ g$ 
 we obtain 
 a   homomorphism
$\langle g \rangle:A(M_0) \to A(M_1)$.
This
homomorphism  does not depend on the choice of $g$ and $\Gamma$.
The proof goes exactly as in Section 8.4; the only additional ingredient 
is the fact  that any two homotopic $\pi$-systems on $\Gamma$
 extending the same splitting of $\partial W$  can be obtained from each other 
by homotopy 
moves at those  components of $W\backslash \Gamma$  which are not adjacent to 
$\partial 
W$. This follows from definitions and the fact that each component of 
$\partial W$ has a base point. 

 It is clear that the definition of $\langle g \rangle $ is equivalent to the 
one 
in  Section 7.4 via the passage from trivializations of
 1-dimensional  $X$-manifolds to the dual splittings and the passage from 
CW-decompositions of 2-dimensional  $X$-cobordisms to the dual  skeletons.  
 Therefore, the independence of $\langle g \rangle $   on the choice of 
$\Gamma$ and $g$ implies the similar claim made in Section 7.4 at Step 2.

Set  $\tau(W)=\langle g \rangle:A(M_0)  \to A(M_1)$.
   For 
completeness, we state an analogue 
of Lemma 7.4.1.

\skipaline \noindent {\bf 8.5.1.  Lemma.  } {\sl Let $M_0,M_1,N$ be 
split 1-dimensional $X$-manifolds (possibly void). If
a 2-dimensional  $X$-cobordism
$(W,M_0,M_1)$ is obtained from two  2-dimensional  $X$-cobordisms
$(W_0,M_0, N)$ and $(W_1,N, M_1)$ by gluing along $N$
then $\tau(W)= \tau(W_1) \circ \tau(W_0):A(M_0)  \to A(M_1)$.}

\skipaline {\sl Proof.}    For $r=0,1$, choose a  skeleton $\Gamma_r$ of 
$W_r$ and a $\pi$-system $g_r$ on $\Gamma_r$ such that  
$(\Gamma_r,g_r)$ extends the splitting of $\partial W_r$.
 Gluing $\Gamma_0$ and $\Gamma_1$ along their feet lying on $N$ we obtain 
a  skeleton, $\Gamma$, of $W$.
 Let $g$ be 
the  unique
  $\pi$-system on $\Gamma$ extending   $g_0$ and $g_1$. 
 The lemma follows from the equality
 $\langle g \rangle =\langle g_1\rangle \,\circ \,\langle g_0 \rangle :A(M_0 ) 
\to  
A(M_1) $. 
   This equality directly follows from the definition of  
$\langle g \rangle$.

 \skipaline  {\bf  Step 3.  } We repeat the constructions made in 
Section 7.4 at Step 3  replacing everywhere the words \lq\lq trivialized, 
trivialization" with \lq\lq split, splitting". 
 This gives a $(1+1)$-dimensional HQFT $(\hat A,\hat \tau)$ with target $X$.

The duality between the skeletons and the regular CW-decompositions establishes 
an equivalence between the construction  of this section  
and  the one of Section 7.   Therefore the   
HQFT obtained via    skeletons coincides with the one constructed in 
Section 7. 

\skipaline \noindent {\bf 8.6.  Theorem.  } {\sl Let $L$ be 
a biangular $\pi$-algebra over an algebraically closed field $K$ of characteristic 0. Then the
associated $(1+1)$-dimensional HQFT $(\hat A,\hat \tau)$ is semi-cohomological.}
 \skipaline  
The proof   is based on the following   lemma.

\skipaline \noindent {\bf 8.6.1.  Lemma.  } {\sl Let  $A$ be a finite-dimensional algebra over a
field $K$   and $\eta:A\otimes
A\to K$  be   the bilinear form defined by 
 $\eta(a,b)=\Tr(x\mapsto {ab}x :A\to 
A) $
 where $a,b,x\in A$. If $\eta$  is non-degenerate then $A$ is semisimple.}

   \skipaline {\sl Proof.}  It suffices to prove that the radical $J(A) \subset A$ of $A$ is zero
(see  [Hu] for the relevant definitions).
If $a\in J(A)$ then for any $b\in A$ there is $n\geq 1$ such that $(ab)^n=0$.
The homomorphism $x\mapsto {ab}x :A\to 
A$ is nilpotent and therefore its trace is equal to 0. Hence, $\eta(a,b)=0$ for any $b\in A$.
By assumption, $\eta$ is non-degenerate, so that $a=0$. 

   \skipaline {\sl Proof of Theorem 8.6.}  We first compute the
$\pi$-center $C$ of $L$, cf. Section 7.1.  (Here we   need no assumptions on $K$.)
Let $(  A,   \tau)$ be the $(1+1)$-dimensional \lq\lq HQFT" with target $X=K(\pi,1)$ 
constructed at Steps 1 and 2 in Section 8.5. Recall that this \lq\lq HQFT"  is defined for split
1-dimensional $X$-manifolds and 2-dimensional $X$-cobordisms with split boundary.  
By definition,  $(\hat A, \hat \tau)$
is the $(1+1)$-dimensional HQFT with target $X=K(\pi,1)$  derived from $(  A,   \tau)$ at
Step 3 in Section 8.5.  For  $\alpha \in \pi$, denote by $S^1_\alpha$ a
pointed  circle  endowed with a map to $X$ representing $\alpha$.
   To compute   the module $$C_\alpha=\hat A_{S^1_\alpha}$$
we   first  provide $S^1_\alpha$ with a splitting, cf. Section 8.5. 
 As a splitting of $S^1_\alpha$
we take any non-base  point $s\in S^1$ endowed with $\alpha$. By definition,
$$A_{(S^1_\alpha,s)}=L_\alpha.$$
Now, the $X$-annulus $C_{-+}
(\alpha;1)=S^1_\alpha\times [0,1]$  (cf. Section 4.6)  is an $X$-cobordism between
two copies of the split circle $(S^1_\alpha,s)$. The simplest skeleton of $C_{-+}
(\alpha;1)$ extending the splitting of the boundary consists of the arc $s\times [0,1]$ and an 
arc connecting two distinct points of $s\times [0,1]$ and winding once around the  
annulus. A direct computation from definitions shows that 
the $K$-homomorphism $\tau(C_{-+}
(\alpha;1)):L_\alpha\to L_\alpha$ is   the  homomorphism $\psi_1$ defined in Section 7.1.
Hence $C_\alpha=\psi_1(L_\alpha) \subset L_\alpha$. We leave it to the reader to check that
 multiplication in $C\subset L$ derived from the HQFT $(\hat A, \hat
\tau)$   coincides with multiplication in $L$. 

Now we can prove Theorem 8.6.  By Lemma 8.6.1, the algebra $L_1$ is semisimple. 
Since $K$ is an algebraically closed field,    $L_1$ is a direct sum of 
matrix rings over $K$. For a matrix ring $L_1=\Mat_n(K)$ with $n\geq 1$, the bilinear form
(7.1.b) and the corresponding homomorphism $\psi_1:L_1\to L_1$   can be
computed explicitly. This computation shows that $\psi_1$ is    a projection of    $\Mat_n(K)$
onto its 1-dimensional center. If  $L_1$ is a direct sum of $m$ matrix rings then  
$\psi_1$ is   a projection of    $L_1$ onto its
$m$-dimensional center $K^m$. Hence $C_1=K^m$ and the crossed $\pi$-algebra $C$ is
semisimple. By Theorem 4.3, the $(1+1)$-dimensional HQFT $(\hat A, \hat \tau)$ is
semi-cohomological.

\skipaline \centerline {\bf 9.  Non-degenerate $\pi$-algebras and lattice 
HQFT's}

 \skipaline 

Throughout this section the group $\pi$ is {\it finite} unless explicitly 
stated to the contrary.

\skipaline \noindent {\bf 9.1. Non-degenerate $\pi$-algebras.} 
Let $L=\oplus_{\alpha\in \pi} L_\alpha$ be a
$\pi$-algebra.    
Recall that for $\ell\in L$, we denote 
by 
$\mu_\ell$ the left multiplication by $\ell$ sending any $a\in L$ into 
$\ell 
a\in 
L$. 
  We say that 
$L$ is {\it non-degenerate} if the bilinear form  
$\eta:L \otimes 
L_{}\to K$
defined by $\eta(a,b)= \Tr(\mu_{ab} :L\to L)$ with $a
 , b\in  L$ is non-degenerate. The trace in this formula is well defined since 
 the underlying $K$-module of $L$  is projective of finite type; here we 
use 
the assumption that $\pi$ is finite.
The usual properties of the trace imply that   $\eta$   is 
symmetric and  that
  $(L,\eta)$ is a 
  Frobenius $\pi$-algebra.

For instance, a biangular $\pi$-algebra $L$ is non-degenerate provided
($\pi$ is finite and) the order $\vert \pi \vert$ of $\pi$ is invertible 
in $K$. 
Note that   the inner product $\eta(a,b)= \Tr(\mu_{ab} )$
is equal to $\vert \pi \vert$ times the inner product considered in 
Section 7.1.

The 
  direct sums and tensor products   of 
non-degenerate  $\pi$-algebras are non-degenerate.
The   push-forward   along a  
homomorphism of 
finite groups $q:\pi'\to \pi$ transforms
a non-degenerate $\pi'$-algebra   into a non-degenerate 
$\pi$-algebra. 
 If $q$ is surjective and   $\vert \Ker \,q\vert $ is 
invertible in $K$, then the pull-back along $q$ of a non-degenerate $  
\pi$-algebra is a 
non-degenerate $\pi'$-algebra. 
 
We shall show   that each  non-degenerate $\pi$-algebra  $L$  gives rise to a
  lattice $(1+1)$-dimensional HQFT with target $K(\pi,1)$. 
 The underlying crossed
 $\pi$-algebra,  $C=\oplus_{\alpha\in \pi} C_\alpha$, of this HQFT is  called the
{\it $\pi$-center} of $L$.
(In the case where $L$ is biangular, this definition is equivalent to the one in Section 7.1 up to
rescaling, cf. Remark 9.7.1.) We give here 
an
   algebraic description of   $C $.

Observe first that for every $\alpha\in \pi$,
 the form $\eta:L_\alpha\otimes L_{\alpha^{-1}}\to K$ defined above is non-degenerate 
and 
yields   a 
canonical element $b_\alpha \in
L_\alpha\otimes L_{\alpha^{-1}}$  (cf.  Lemma 7.1.1).
As in Section 7.1, we have an expansion ${ b}_\alpha=\sum_{i} 
 p_i^{\alpha}\otimes  q_i^{\alpha}$.
The expression $\hat { b}_\alpha=\sum_i p_i^{\alpha}  
q_i^{\alpha} 
\in L_1$ does not depend on the choice of this expansion.
We claim that
$$\sum_{\alpha\in \pi} \hat { b}_\alpha = 1_L.\leqno (9.1.a)$$
   To prove it,   observe that for any $\ell \in L_1$,
 $$\eta(\ell,1_L)= \Tr(\mu_\ell :L \to 
L)=\sum_{\alpha\in \pi} \Tr(\mu_\ell\vert_{L_\alpha} :{L_\alpha} \to 
{L_\alpha})
$$
$$=\sum_{\alpha\in \pi} \sum_i\eta(\ell \,p_i^{\alpha}, q_i^{\alpha} )
=\sum_{\alpha\in \pi} \eta(\ell , \sum_i p_i^{\alpha} q_i^{\alpha} )
= \sum_{\alpha\in \pi} \eta(\ell , \hat { b}_\alpha)=\eta(\ell, 
\sum_{\alpha\in \pi}  \hat { b}_\alpha )$$
 where we use      the expression for the trace given in Lemma 7.1.1.
Since the form $\eta: L_1\otimes L_{1}\to K$ is non-degenerate,  
$\sum_{\alpha\in \pi}  \hat { b}_\alpha=1_L$.

 As in  Section 7.1, we   use ${ b}_\alpha=\sum_{i} 
 p_i^{\alpha}\otimes  q_i^{\alpha}$ to define a $K$-homomorphism
 $\psi_\alpha:L\to L$    by $\psi_\alpha(a)= \sum_{i} 
 p_i^{\alpha}\,a\, q_i^{\alpha}$
 where $a\in L$.  This
 homomorphism
 does not depend on the choice of  the  expansion of ${ b}_\alpha$   and sends 
each $ 
 L_\beta  \subset L$ into $L_{ \alpha\beta \alpha^{-1}}$. 

For $\omega\in \pi$, we can push  $(L,\eta)$ back along the conjugation $\alpha\mapsto
\omega \alpha \omega^{-1}: \pi \to \pi$. This gives a Frobenius $\pi$-algebra
$L^{\omega}$ such that $L^{\omega}_\alpha= L_{\omega\alpha\omega^{-1}}$ for all $\alpha\in
\pi$. Let $P$ be the direct sum of the Frobenius $\pi$-algebras
$L^{\omega}$ over all $\omega\in \pi$. Thus,
$$P_\alpha=\bigoplus_{\omega\in \pi} L_{\omega\alpha\omega^{-1}}$$
for all $\alpha\in
\pi$. Given $\beta\in \pi$, consider a homomorphism
$\Psi_\beta:P_\alpha\to P_{\beta\alpha\beta^{-1}}$
defined by the block-matrix 
$[\psi_{\omega'\beta\omega^{-1}}]_{\omega,\omega'\in\pi}$.
Thus, for $a\in L_{\omega\alpha\omega^{-1}}\subset P_\alpha$
we have
$$\Psi_\beta(a)= \bigoplus_{\omega'\in \pi} \psi_{\omega'\beta\omega^{-1}}
(a) \in \bigoplus_{\omega'\in \pi} 
L_{\omega'\beta\alpha\beta^{-1}(\omega')^{-1}}
=P_{\beta\alpha\beta^{-1}}.$$ 
 It turns out that
$\Psi_\beta
\Psi_{\beta'}=\Psi_{\beta \beta'}$ for any $\beta, \beta' \in \pi$
and
 $\Psi_\beta$ is adjoint to $ \Psi_{\beta^{-1}}$ with respect to the inner product in $P$.
 In particular, $\Psi_1$ is an orthogonal projection onto a submodule
 $C=\oplus_{\alpha\in \pi} C_\alpha $ of $ P$ where $C_\alpha= \Psi_1
 (P_\alpha) \subset P_\alpha$. 
  Formula (9.1.a) implies that $1_P\in C$.
   It turns out that $C$ is a subalgebra of $P$.
 Moreover, the inner product in $P$ restricted to $C$      and the
homomorphism  $\pi\to \Aut C$ sending each 
$\beta\in 
\pi$
 to $\Psi_\beta\vert_C:C\to C$ make $C$ a crossed $\pi$-algebra.  
 The crossed $\pi$-algebra $C$ is   the $\pi$-center of
$L$.
All 
these
 claims follow directly by applying the  definitions of Section 4 to the 
HQFT
 associated with $L$ (cf. Section  9.6). 

\skipaline \noindent {\bf 9.1.1. Remark.}  It is curious to note (we shall not use it)
that a  non-degenerate $\pi$-algebra  $L$ is separable. Recall that a  unital
associative $K$-algebra $A$ is separable if there is $b\in A\otimes A$ such that (i) the 
algebra multiplication $ A\otimes A\to A$ sends $b$ into $1$ and  (ii) for any $a\in A$,
$(a\otimes 1) b= b(1\otimes a)$ (see [Pi, Chapter 10]). Such an element $b$ is called  a  
separating idempotent of $A$. We claim that   $b=\sum_{\alpha\in \pi}
{ b}_\alpha \in L\otimes L$ is a separating idempotent of $L$. Condition (i) follows from
(9.1.a), Condition  (ii) follows from the fact that for any 
$\alpha, \beta\in \pi, a\in L_\alpha$,  $$ (a\otimes 1) \,b_\beta=b_{\alpha\beta}
(1\otimes a).\leqno (9.1.b)$$
Summing over all $\beta$, we obtain $(a\otimes 1) b= b(1\otimes a)$.
To prove (9.1.b), it suffices to show that for any $x\in L_{\beta^{-1} \alpha^{-1}}$,
$$\sum_i \eta(a p_i^{\beta},x) \,q_i^{\beta} = \sum_j \eta(  p_j^{\alpha\beta},x) 
\,q_j^{\alpha\beta}  a.$$
Using  
 $\eta(a p_i^{\beta},x)=\eta(xa, p_i^{\beta})$, we easily observe that both
sides are equal to  $xa$.

\skipaline \noindent {\bf 9.2.  $\pi$-systems on CW-complexes re-examined.}
 Let $T$ be a  regular finite  CW-complex with underlying topological space 
$\tilde T$ and a distinguished 
set of base vertices.  
 Let $F_T$ be the group formed by the maps $ \Vert (T)\to \pi$ 
taking value $1\in \pi$ on    all   base vertices of $T$. 
Clearly, $F_T=\pi^{n_T}$
where  $n_T$ is the number of vertices of $T$ distinct from the base 
vertices.
Each homotopy class $G$  of $\pi$-systems  on $T$ 
is an 
orbit  of the action of   $F_T$ on the set of $\pi$-systems on 
$T$, see 
Section 7.2.
This implies that $\card (G)$ is an integer divisor of $\vert 
F_T\vert=\vert 
\pi\vert^{n_T}$.
 More precisely,  $\vert 
\pi\vert^{n_T}/\card (G) =
\vert \Stab_g \vert$ where $\Stab_g \subset F_T$ is the stabilizer of an 
element 
$g\in G$.
If   $\tilde T$ is connected and the set of   base vertices is non-void 
then 
$\Stab_g=1$ so that $\card (G)=\vert \pi\vert^{n_T}$. 
If   $\tilde T$ is connected and the set of   base vertices is empty then 
the  group $\Stab_g$
can be computed  in homotopy terms. Consider 
 the map $\tilde g:\tilde T\to X=K(\pi,1)$ determined by $g\in G$ and
the 
induced homomorphism $\tilde g_{\#}:\pi_1(\tilde T)\to \pi$. Then $  \Stab_g  $ 
is isomorphic to the subgroup 
  of $\pi$ consisting of the elements of $\pi$ 
commuting with all elements of $\tilde g_{\#}(\pi_1(\tilde T))\subset \pi $.

 In Section 9.4 we shall need a   notion of  enriched $\pi$-systems on 
$T$.
By an {\it enriched $\pi$-system} on $T$ we mean a pair $(\omega, g)$
where $\omega$ is an arbitrary map from the set of base vertices of $T$ to 
$\pi$ 
and 
$g $ is a $\pi$-system on $T$. Given an enriched $\pi$-system 
$(\omega,\pi)$,  we  define a $\pi$-system $\omega g$ on $T$ as follows:    
extend
$\omega$ to all vertices of $T$ by assigning $1\in \pi$ to all non-base 
vertices 
and set
$(\omega g)_e=\omega (i_e)\, g_e\, (\omega (t_e))^{-1}$ for any  
 $e\in \Edg(T)$.  The induced map $\tilde {\omega g}:\tilde T \to X$ is 
obtained 
from 
$\tilde { g}:\tilde T \to X$ by pushing the image of each  base vertex 
$z\in T$ 
along 
a 
loop representing $\omega(z)\in \pi$.  
The homomorphisms $\pi_1(T,z)\to \pi$ induced by  
$\tilde g$ and $\tilde {\omega g}$ are    conjugated by 
$\omega(z)\in 
\pi$.
Hence the $\pi$-systems $g$ and $\omega g$ are  not homotopic  
 unless $\omega=1$.

If $T$ is a CW-decomposition of an $X$-manifold or an $X$-cobordism
then an enriched $\pi$-system $(\omega,g)$ on $T$ is {\it characteristic} 
if 
the $\pi$-system $g$ is characteristic.

\skipaline \noindent {\bf 9.3.  State sums on closed $X$-surfaces.} Fix a 
non-degenerate $\pi$-algebra $L$.  Let $W$ be   
 a closed   oriented surface 
 endowed with a  map  $W\to X=K(\pi,1)$.   Here we
define a state sum invariant  $\tau(W)=\tau_L(W) \in K$. Choose a  regular 
CW-decomposition 
$T$ of $W$. 
Let $G$ be the homotopy class of
characteristic $\pi$-systems on $T$. 
Repeating the definitions of  Section 7.3 for   $g\in G$ we obtain  a state sum  
  $\langle g \rangle\in K$. In   general this state sum  depends on the 
choice 
of 
$g$ in 
$G$. 

\skipaline \noindent {\bf 9.3.1. Claim.  } {\sl 
 The  sum 
$\sum_{g\in 
G} 
\langle g \rangle$ does not depend on 
the choice of $T$.}

\skipaline

   We 
outline a     proof   in Section 9.5. Set $\langle W \rangle_L=\sum_{g\in G} 
\,\langle g \rangle\in K$.
  Claim 9.3.1   implies 
that 
$\langle W \rangle_L$ is a well defined invariant of the closed $X$-surface $W$.
To include this state sum invariant into an HQFT, we need to renormalize it
(cf.   the proof of Lemma 9.4.1).
Set $$\tau(W)\,=\, \frac {\vert \pi \vert^{n_T}}{\card (G)} \,\langle W 
\rangle_L=\, \, \frac {\vert \pi 
\vert^{n_T}}{\card 
(G)} \,\sum_{g\in G} 
\,\langle g \rangle\in K \leqno (9.3.a)$$ 
where   $n_T$ is the number of vertices of $T$. According to the results 
of 
Section 9.2, the  
 number  ${\vert \pi \vert^{n_T}}/{\card (G)}$ is
  an integer. This number is   multiplicative with respect 
to 
disjoint 
union of closed $X$-surfaces and  can be computed in 
homotopy 
terms, see Section 9.2. In particular,  this
 number    is independent 
of $T$ and $L$. 
Therefore
$\tau(W) $ is a well defined invariant of  $W$.
By the very definition, $\tau(W) $ depends only on the homotopy class
of the   characteristic map $W\to X$ and is multiplicative with respect 
to 
disjoint union 
of closed $X$-surfaces.

It is useful to rewrite the definition of $\tau(W)$ in terms of the group 
$F_T=\pi^{n_T}$
formed by the maps $\Vert (T)\to \pi$ and its action on the set of 
$\pi$-structures on $T$. Namely,
for any characteristic $\pi$-system $g\in G$ on $T$, we have
$$\tau(W)=\sum_{\gamma\in F_T} \langle \gamma g\rangle.$$

We can give an explicit formula for $\tau(W)$
when $W$ is a  closed connected oriented surface of genus 
$n\geq 1$. Using the same CW-decomposition of $W$ and the same notation as in  
Section 7.3, we   obtain 
$$\tau(W)=\sum_{\beta\in \pi}\, \,\eta\,(\,\prod_{r=1}^n  [b_{\beta 
\alpha_{2r-1}\beta^{-1}},
b_{\beta \alpha_{2r} \beta^{-1}}], 1_L)$$
where $\prod_{r=1}^n [b_{\beta \alpha_{2r-1}\beta^{-1}},
b_{\beta \alpha_{2r} \beta^{-1}}]  \in L_1$. 
Similarly,  
$$\tau(S^2)=\vert \pi \vert\, \eta(1_L,1_L)=\vert \pi \vert\, \Dim L.$$

\skipaline \noindent {\bf 9.4.  Lattice construction of a 
$(1+1)$-dimensional  HQFT.} We extend the invariant   of closed 
$X$-surfaces constructed in Section 9.3  to a $(1+1)$-dimensional  HQFT with 
target 
$X=K(\pi,1)$.
The construction follows the same lines as in Section 7.4. 

\skipaline  {\bf  Step 1.  } 
Let $M$ be a 
 1-dimensional $X$-manifold with distinguished CW-decomposition $T$
 (such that the base points of the components of $M$ 
are 
among 
the vertices).
  We provide each edge $e$ of $M$ with canonical orientation induced by 
the 
one of $M$.
 For a  $\pi$-system  $g$ on $T$,   we define a 
module
 $A(M,g)=\otimes_e {L_{g_e}}$
 where $e$ runs over all canonically oriented edges of $T$.
 Set 
 $$A(M)= \bigoplus_{(\omega,g)} A(M,\omega g) \leqno (9.4.a)$$
 where $(\omega,g)$ runs over all  characteristic enriched $\pi$-systems on 
$T$
 (so that $\tilde g:M\to X$ is homotopic to the  characteristic map 
$M\to 
X$). Since the number of 
such systems is finite, the module $A(M)$ is projective of finite type.
It is clear that for  disjoint 1-dimensional 
$X$-manifolds $M,N$ endowed with CW-decompositions, we have $A(M\coprod N)=A(M) 
\otimes A(N)$.
If $M=\emptyset$ then by definition $M$ has a unique enriched $\pi$-system 
and  
$A(M)=K$.

The module $A(M,\omega g)$ can be explicitly computed as follows.
Assume first that   $M=S^1$.
 Starting from the base point, $z\in M$, and moving along $M$ in the 
direction 
determined by the orientation of $M$ we meet consequtively   the 
oriented edges, 
$e_1,...,e_n$, of $T$ where
  $n\geq 1$ is the number of edges of   $T$. Then
  $$A(M,\omega g)=  L_{\omega(z) g_{e_1}} \otimes (\otimes_{r=2}^{n-1} L_{  
g_{e_r}} )
\otimes  L_{  g_{e_n} (\omega(z))^{-1}}.$$
 If $M$ has several components then restricting $(\omega,g)$   we obtain 
enriched $\pi$-systems on these components and $A(M,\omega g)$ 
is  
the tensor product of the corresponding modules.

\skipaline  {\bf   Step 2.  }  Consider   a 2-dimensional  
$X$-cobordism
$(W,M_0,M_1)$. 
Assume that at least one of the bases $M_0,M_1$ is non-void and fix a 
CW-decomposition of 
the bases. 
Using    (9.4.a), 
we define a  $K$-homomorphism  $\tau(W):A(M_0) \to A(M_1)$  by a  
  block-matrix  of  
homomorphisms 
$\tau(W; (\omega_0,g_0), (\omega_1,g_1)): A(M_0, \omega_0 g_0) \to A(M_1, 
\omega_1 
g_1)$
where $(\omega_r,g_r)$ runs over  characteristic  enriched $\pi$-systems on 
the 
given CW-decomposition 
of $M_r$,
 for 
$r=0,1$. Fix   a pair
$(\omega_0,g_0), (\omega_1,g_1)$.
Choose a CW-decomposition 
$T$ 
of $W$ extending the given CW-decomposition of $\partial 
W=(-M_0)\cup M_1$.
 Denote by
$G$  the set of characteristic
$\pi$-systems $g$ on $T$  extending   $g_0\cup g_1$.
The set $G$ is non-void and finite. 
It is clear that for any $g\in G$ the pair $(\omega=\omega_0 \cup \omega_1,g)$
is an enriched $\pi$-system on $T$. Consider the  $\pi$-system $\omega g$ on 
$T$ extending the $\pi$-system $(\omega_0 g_0) \cup (\omega_1 g_1)$ on the 
boundary.
Repeating the definitions of  Section 7.4 for  
$\omega g$   we obtain  a homomorphism 
  $\langle \omega g \rangle: A(M_0, \omega_0 g_0) \to A(M_1, \omega_1 g_1) $.
  Set
$$\tau(W;  (\omega_0,g_0), (\omega_1,g_1))=  \sum_{g\in G} 
\langle \omega g \rangle : A(M_0, \omega_0 g_0) \to A(M_1, \omega_1 g_1) .$$
We  claim that this
 homomorphism
 is independent of the 
choice of 
$T$; for a proof, see Section 9.5.

By the very definition, $\tau(W) $ depends only on the homotopy class
of the   characteristic map $W\to X$ and satisfies (1.2.5).

\skipaline \noindent {\bf 9.4.1.  Lemma.  } {\sl Let $M_0,M_1,N$ be 
  1-dimensional $X$-manifolds (possibly void) with fixed CW-decompositions. If
a 2-dimensional  $X$-cobordism
$(W,M_0,M_1)$ is obtained from two  2-dimensional  $X$-cobordisms
$(W_0,M_0, N)$ and $(W_1,N, M_1)$ by gluing along $
N$
then $\tau(W)= \tau(W_1) \circ \tau(W_0):A(M_0) \to A(M_1)$.}

\skipaline {\sl Proof.}   Because of the  multiplicativity of $\tau$ with 
respect to disjoint union, it suffices to consider the case where $W$ is 
connected and $N\neq \emptyset$. For $r=0,1$, fix a regular CW-decomposition 
$T_r$ of 
$W_r$ extending the 
given 
CW-decomposition of the boundary.
 Gluing $T_0$ and $T_1$ along $N$ we obtain a regular CW-decomposition, $T$, of 
$W$.

 Assume first that at least one of the 
manifolds $M_0,M_1$ is non-void.
Fix a  characteristic enriched $\pi$-system $(\omega_r,g_r)$  on the  given 
CW-decomposition of $M_r$
for
$r=0,1$. We should prove that 
$$\tau(W;  (\omega_0,g_0), (\omega_1,g_1))=  \sum_{y} 
\tau(W_1;  y, (\omega_1,g_1)) \,\circ \, \tau(W_0;  (\omega_0,g_0), y) 
\leqno(9.4.b)$$
where $y$ runs over all characteristic enriched $\pi$-systems on $N$.

Fix a characteristic $\pi$-system  $h$ on $N$. 
Choose a characteristic $\pi$-system $H_r$ on $T_r$ extending $g_r\cup h$. 
It is clear that there is a unique $\pi$-system, $H$, on $T$ such that 
$H\vert_{T_r}=H_r$ for $r=0,1$. The corresponding map $\tilde H:W\to X$
is homotopic to the  map obtained by gluing the characteristic maps
$W_0\to X, W_1\to X$ along $N$.

Denote by
$F$ be the group formed by the 
maps $ \Vert (T) \to \pi$ taking the value $1\in \pi$ on all vertices of 
$T$ lying in $\partial W$.
 We introduce four subgroups 
$F_0,F_1,F_N,F_b$ of $F$ as follows.
For $r=0,1$, the group $F_r$   consists of   $\gamma\in F$ such that 
$\gamma(u) 
=1\in \pi$   for all vertices $u$ of $T$  except possibly those  lying in 
$W_r\backslash \partial W_r$. 
The group $F_N$   consists of   $\gamma\in F$ such that $\gamma(u)= 1\in 
\pi$ 
for all vertices $u$ of $T$  except possibly those lying in $N$ and 
 distinct from the   base vertices of $N$.
 The group $F_b$   consists of   $\gamma\in F$ such that $\gamma(u)= 1\in 
\pi$ 
for any vertex $u$ of $T$   distinct from the base vertices of $N$.
 
 The subgroups $F_0,F_1,F_N,F_b$ of $F$ commute with each other.
Every   $\gamma\in F$ splits uniquely as the product 
$\gamma=\gamma_0\gamma_1\gamma_N\gamma_b$ with $\gamma_0\in F_0,\gamma_1 \in 
F_1, 
\gamma_N\in F_N, \gamma_b\in F_b$.   

 By definition,
  $$\tau(W; (\omega_0,g_0), (\omega_1,g_1))=  \sum_{g\in G} \,
  \langle \omega g \rangle $$
 where $\omega=\omega_0 \cup \omega_1$ and $G$ is the set of
characteristic $\pi$-systems   on $T$ extending $g_0\cup g_1$.
 Thus,  every $g\in G$ is homotopic to 
$H$
 and coincides with $H$ on $\partial W$.
The definition of homotopy for $\pi$-systems 
 (see Section 7.2) and     the assumption that each component 
of 
$\partial W$ contains a base point imply that $G$ is the orbit of   $H$ under 
the 
action of $F$ on the set of $\pi$-systems   on $T$.
 Since $\partial W\neq \emptyset$,  this action is fixed point free and 
therefore
$$\tau(W; (\omega_0,g_0), (\omega_1,g_1))=\sum_{g\in G} \, \langle  \omega g 
\rangle=\sum_{\gamma\in F}\,  \langle  
\omega 
\gamma 
H \rangle $$
 $$=\sum_{\gamma_0\in F_0}
 \sum_{\gamma_1\in F_1} \sum_{\gamma_N\in F_N} \sum_{\gamma_b\in F_b} 
 \, \langle  {\omega \gamma_0\gamma_1\gamma_N\gamma_b H}  \rangle
 :A(M_0, 
\omega_0 g_0) \to A(M_1, \omega_1 g_1).$$
The claim $(\ast)$ in the proof of Lemma 7.4.1 directly implies that 
 the homomorphism
 $\langle  {\omega \gamma_0\gamma_1\gamma_N\gamma_b H}  \rangle $
 is the composition of 
 the homomorphisms
 $$\langle   {\omega_0 \gamma_0 \gamma_N\gamma_b H_0} \rangle:A(M_0, 
\omega_0 
g_0) \to A(N, \gamma_N\gamma_b h)$$
 and
  $$\langle  {\omega_1 \gamma_1 \gamma_N\gamma_b H_1}  \rangle: 
  A(N, \gamma_N\gamma_b h) \to A(M_1, \omega_1 g_1).$$  
  Summing over  ${\gamma_0\in F_0}
,{\gamma_1\in F_1}$ we obtain that
$$\tau(W; (\omega_0,g_0), (\omega_1,g_1))= $$
$$= \sum_{\gamma_N\in F_N} \sum_{\gamma_b\in F_b} \,\tau(W_1; (\gamma_b, 
\gamma_N   
h), (\omega_1,g_1)) \circ
\tau(W_0;  (\omega_0,g_0), (\gamma_b ,\gamma_N  h)).$$
Finally note that when $\gamma_N$ runs over $ F_N$
and $\gamma_b$ runs over $F_b$  the enriched $\pi$-system $y= (\gamma_b, 
\gamma_N   
h)$
runs over all characteristic enriched $\pi$-systems on $N$. This gives 
(9.4.b). 

 In the case $M_0=M_1=\emptyset$ the same argument  gives
 $$\tau(W)=\sum_{\gamma\in F} \langle \gamma H \rangle=\tau(W_1)\circ 
\tau(W_0).$$

 \skipaline  {\bf   Step 3.  } We apply the   constructions given at 
Step 3 
in Section 7.4 to the  preliminary 
HQFT defined above. It suffices to replace everywhere the word \lq \lq 
trivialized" with \lq\lq provided with a CW-decomposition".   This gives a  
  $(1+1)$-dimensional  HQFT $(\hat A,\hat \tau)$ 
associated with  the non-degenerate $\pi$-algebra  $L$. Note that for any 
closed $X$-surface $W$ we have  $\hat 
{\tau} (W)=\tau(W)$. 

\skipaline  {\bf 9.5. Proof of Claim 9.3.1.  }
The proof goes along the same lines as the proof of Claim 7.3.1 given in 
Section 8. 
One  translates the state sum in terms of the   skeletons 
and proves the invariance under the local moves discussed in Section 
8.1.
 The invariance under the contraction move  
directly 
follows 
 from definitions and Lemma 7.1.3 which applies  to   $L$ without any changes.  
Consider  the loop move $\Gamma\mapsto \Gamma'$. Every $\pi$-system $g$ on 
$\Gamma$ lifts to a finite family, 
$Y(g)$, of $\pi$-systems on $\Gamma'$. (In fact, $\card (Y(g))=\vert 
\pi\vert$). It
follows from definitions and   (9.1.a) that $\langle g\rangle 
=\sum_{g'\in Y(g)} \langle g'\rangle$. This implies the invariance of the sum $
\sum_{g\in  G} \langle g\rangle$ under the loop move. 
   In contrast to the constructions of Section 7,  we do not need   to prove 
the invariance under homotopy moves.

The case of $X$-cobordisms  is   similar.

\skipaline \noindent {\bf 9.6.  Theorem.  } {\sl Let $L$ be 
a non-degenerate $\pi$-algebra over an algebraically closed field $K$ of characteristic 0. Then
the associated $(1+1)$-dimensional HQFT $(\hat A,\hat \tau)$ is semi-cohomological.}
 \skipaline

   \skipaline {\sl Proof.}  Let $\alpha\in \pi$ and   
$b_\alpha= \sum_i p_i^{\alpha}\otimes q_i^{\alpha} \in L_\alpha\otimes L_{\alpha^{-1}}$
be  the canonical element associated with the restriction of the inner product  $\eta$ (defined
in Section 9.1) onto  $L_\alpha\otimes L_{\alpha^{-1}}$.  Recall the homomorphism
$\psi_\alpha:L\to L$ defined  by $\psi_\alpha(a)= \sum_{i} 
 p_i^{\alpha}\,a\, q_i^{\alpha}$
 where $a\in L$.
 Observe that for any $b\in L_1$ and  $x\in L_{\alpha^{-1}}$,  
$$ \sum_i \eta(bp_i^{\alpha},x) \, q_i^{\alpha}=
\sum_i \eta(xb, p_i^{\alpha} ) \, q_i^{\alpha}=xb=
\sum_i \eta(x, p_i^{\alpha}) \, q_i^{\alpha}b=
\sum_i \eta(p_i^{\alpha},x) \, q_i^{\alpha}b.$$
Therefore  
$$\sum_i bp_i^{\alpha}\otimes q_i^{\alpha}=\sum_i p_i^{\alpha}\otimes q_i^{\alpha}b.$$
This     implies that for any $a,b\in L_1$,
$$b \psi_\alpha (a)= \sum_i bp_i^{\alpha} a q_i^{\alpha}= \sum_i p_i^{\alpha}a q_i^{\alpha}b=
\psi_\alpha (a) b.$$
In other words,   $\psi_\alpha(L_1)$   lies in the center, $Z(L_1)$,  of $L_1$
for all $\alpha\in \pi$.

Recall the module $P_\alpha$ introduced in Section 9.1. Each   $a\in P_\alpha$
splits uniquely as a  sum $a=\sum_{\omega\in \pi} a_\omega$ with 
$ a_\omega\in L_{\omega\alpha\omega^{-1}}$.  We view $\{a_\omega\}_\omega$ as the
coordinates of $a$. 
Computations  similar to the ones in the proof of Theorem 8.6 give the algebraic description of
the   $\pi$-center $C$ of  $L$ formulated in Section 9.1.  
In particular, $C_\alpha  \subset P_\alpha$ is the image of the  projection 
$\Psi_1:P_\alpha  \to \ P_\alpha$ given in the coordinates by the block-matrix 
$[\psi_{\omega' \omega^{-1}}]_{\omega, \omega'\in \pi}$. 
An element $a\in P_\alpha$ lies in $C_\alpha$ if
and only if for all $\omega'\in \pi$,
$$a_{\omega'}=\sum_{\omega\in \pi} \psi_{\omega' \omega^{-1}} (a_\omega).$$
Taking $\alpha=1$ and applying the result of the preceding paragraph we obtain that  
$$C_1\subset Z(L_1)^{\vert \pi \vert} \subset P_1=L_1^{\vert \pi \vert}.$$

Using an appropriate skeleton of a disc with two holes, we can compute multiplication in $C$
  as follows. If $a\in C_\alpha, b\in C_\beta$, then the $\nu$-th coordinate
of $ab\in C_{\alpha\beta}$  corresponding to $\nu \in \pi$ is computed by  $$(ab)_\nu=
 \sum_{\omega, \mu\in \pi}  \psi_{\nu \omega^{-1}} (a_\omega)
\psi_{\nu \mu^{-1}} (b_\mu)=
(\sum_{\omega \in \pi}  \psi_{\nu \omega^{-1}} (a_\omega)) (\sum_{  \mu\in \pi}
\psi_{\nu \mu^{-1}} (b_\mu)) = a_\nu b_\nu.$$
In particular, $C_1$ is a subalgebra of $P_1=L_1^{\vert \pi \vert}$. By the results above,
$C_1$ is a subalgebra of $Z(L_1)^{\vert \pi \vert}$.
 
The argument given in the proof of Lemma 8.6.1 shows that the algebra $L_1$ is semisimple.
  Since $K$ is an
algebraically closed field,    $L_1$ is a direct sum of  matrix rings over $K$ and its center  
$Z(L_1)$ is a direct sum of several copies of $K$. 
Hence the algebra $C_1$   is a direct sum of several copies of $K$  and the
crossed $\pi$-algebra $C$ is semisimple. By Theorem 4.3, the   HQFT $(\hat
A, \hat \tau)$ is semi-cohomological.

\skipaline  {\bf 9.7. Remarks.  }  1.  Let $L$ be a biangular $\pi$-algebra. If 
 $\vert \pi \vert$ is invertible in $K$ then $L$ is non-degenerate and we have 
two  $(1+1)$-dimensional HQFT's associated to $L$: the HQFT $(A_1,\tau_1)$ 
defined  in Section 7 and the HQFT $(A_2,\tau_2)$ defined  in this section.
These HQFT's are  equivalent up to rescaling.
  In particular, for a closed $X$-surface $W$ we have
  $\tau_2(W)=\vert \pi\vert^{\chi(W)} \tau_{1}(W)$. The factor $\vert 
\pi\vert^{\chi(W)}$ comes
from the normalization factor in (9.3.a) and from the fact that the inner 
products considered in Sections 7.1 and 9.1 differ by a factor of
$\vert \pi \vert$. It can be shown that $A_1(M)=A_2(M)$ for any 1-dimensional 
$X$-manifold $M$.
Thus the construction of an HQFT given in this section 
generalizes the construction of  Section 7   in the case of finite $\pi$
with $\vert \pi \vert$  invertible in $K$.

2. The HQFT  constructed in this section  is  
additive (resp. multiplicative) with respect to direct sums (resp. tensor 
products) of non-degenerate 
$\pi$-algebras.
The invariant $\tau$ of closed $X$-surfaces defined in Section 9.3 is compatible 
with     
pull-backs and push-forwards of   non-degenerate group-algebras
as follows. Consider
 a group homomorphism $q:\pi' \to \pi$.
 Let $f_q:X'=K(\pi',1)\to K(\pi,1)=X$ be the map induced by $q$.
  Let $L=q_*(L')$ be the non-degenerate $ \pi$-algebra obtained as the push-forward along $q$
  of a non-degenerate $\pi'$-algebra  $  L'$.   Consider a closed   $X$-surface 
$(W, 
g:W\to X)$. It follows from definitions that $
\langle W,g\rangle_L=\sum_{g'} \langle W, g'\rangle_{L'}$
where $g'$ runs over the homotopy classes of maps $W\to X'$
such that $f_q g'$ is homotopic to $g$. For connected $W$, this can be rewritten 
as follows:
$$\tau_L(W,g)=\sum_{g'} \frac {\vert \Stab_{ g} \vert } {\vert \Stab_{g'} \vert 
} \,\tau_{L'} (W,g')$$
where $\vert \Stab_{g} \vert$ is the number of elements of $\pi$ commuting with 
all elements of the image of 
the homomorphism $\pi_1(W) \to \pi$ induced  by $g$.

 To consider the pull-backs, assume that $q(\pi')=\pi$  and    $\vert \Ker 
\,q\vert 
$ is 
invertible in $K$.
Let $q^*(L)$ be the non-degenerate $\pi'$-algebra obtained as the pull-back 
along $q$
  of a non-degenerate $  
\pi$-algebra $L$. 
For  each closed $X'$-surface $W=(W, 
g':W\to X')$, we have a closed $X$-surface 
 $W_q=(W, f_q\, g':W\to X)$. It follows from definitions that
 $\tau_{q^*(L)} (W)= \vert\Ker \,q \vert^{\chi(W)} \tau_L (W_q)$. For instance, 
consider the trivial homomorphism $q:\pi 
\to \{1\}$. If $L$ is the 1-dimensional 
$\{1\}$-algebra $K$, then $q^*(L)$ is the group ring $K[\pi]$
 and $\tau_{K[\pi]} (W)=\vert\pi \vert^{\chi(W)}\, \tau_L (W_q)=\vert\pi 
\vert^{\chi(W)}$.

3.  It would be interesting to generalize 
the  constructions of this section to  $\pi$-algebras over an infinite group 
$\pi$.
 A part of these constructions  can be carried over to so-called finite  
$\pi$-algebras. A $\pi$-algebra $L$ is {\it finite} if $L_\alpha=0$ for all but 
finitely many $\alpha\in \pi$.
The definition of non-degeneracy applies to finite $\pi$-algebras word for 
word. If $L$ is a non-degenerate finite $\pi$-algebra over a possibly infinite 
group $\pi$ then 
 the state sum $\sum_{g\in G} 
\,\langle g \rangle\in K$ in Section 9.3 is finite  and gives a well defined 
invariant of the closed $X$-surface $W$. However, the normalization factor 
$  {\vert \pi \vert^{n_T}}/{\card (G)}$ appearing in Section 9.3 
may be infinite and the modules appearing in Section 9.4 may be of infinite 
type.
Another  
possible approach is to consider   topological groups and to replace 
state sums with  integrals.

\skipaline \centerline  {\bf 10.  Further examples of group-algebras}

\skipaline
In this section we have collected a few miscellaneous   examples and constructions 
of group-algebras.
 
\skipaline \noindent {\bf    10.1.  Biangular group-algebras from   algebra
automorphisms.  }   Let $A$ be an associative unital algebra over $K$ whose underlying
$K$-module is projective of finite type. Let $\pi$ be a group acting on $A$ by algebra
automorphisms.  Consider  the direct sum $L=\oplus_{\alpha\in \pi} A\alpha$ of $\card
(\pi)$ copies of $A$ numerated by the elements of $\pi$. We provide $L$ with multiplication   by
$(a\alpha) (b\beta)= (a \alpha(b) ) (\alpha \beta)$
where $a,b\in A$ and $\alpha,\beta\in \pi$. It is easy to check that $L$ is an associative
algebra.    We provide $L$ with the structure of a $\pi$-algebra by  
$L_\alpha=A\alpha$ for all $\alpha\in \pi$. This $\pi$-algebra always satisfies condition (i) in
the definition of biangular $\pi$-algebras. It satisfies condition (ii) if and only if 
the bilinear form $\eta:A\otimes A\to K$
defined by 
 $\eta(a,b)=\Tr(\mu_{ab}:A\to 
A)$ is non-degenerate. The form $\eta$ is non-degenerate for instance when $A$ is a direct sum
of matrix  rings $\Mat_n(K)$ such that $n$ is invertible in $K$.
This can be easily deduced from the following simple lemma.

\skipaline \noindent {\bf 10.1.1.  Lemma.  } {\sl Let $P,Q$ be 
 free $K$-modules of finite rank. Then: (i) the pairing 
$\Hom (P,Q) \otimes \Hom (Q,P) \to K$ sending   $(f\in  
\Hom (P,Q), g \in \Hom (Q,P))$ to $\Tr(fg)=\Tr(gf)$ is non-degenerate;
(ii) for any $\ell \in \Hom (Q,Q)$,
the trace of the endomorphism of $\Hom(P,Q)$ sending any $f\in \Hom (P,Q)$ into
$ \ell f\in \Hom (P,Q)$
is equal to $\Tr(\ell)\, \Dim (P)$.}

 \skipaline \noindent {\bf   10.2. More   non-degenerate and biangular
group-algebras.} 
  Let $\{V_s\}_{s\in S}$ be a family of free $K$-modules of finite rank
numerated by elements of a finite set $S$.  
With every left action of $\pi$ on $S$ we associate a $\pi$-algebra
$L$ as follows.  For $\alpha\in \pi$, set
$$L_\alpha= \bigoplus_{s\in S} \Hom (V_s, V_{\alpha (s)}).$$
Each element $a\in L_\alpha$ is determined by its \lq\lq coordinates"
$\{a_s\in   \Hom (V_s, V_{\alpha(s)}) \}_{s\in S}$. The $K$-linear structure 
in $L$ is   coordinate-wise. If $a\in L_\alpha, b\in L_\beta$
then the product $ab\in L_{\alpha\beta}$ is defined in  coordinates by
$(ab)_s=a_{\beta(s)} b_s\in  \Hom (V_s, V_{\alpha \beta (s)})$.
It is obvious that this multiplication is associative and makes 
$L=\oplus_{\alpha\in \pi} L_\alpha$ a $\pi$-algebra. The unit $1_L$ is 
determined by $(1_L)_s= \id_{V_s}$ for all $s\in S$. 
It follows from Lemma 10.1.1 that $L$ is  non-degenerate
if (and only if)  $\pi$ is finite and  either $V_s=0$ for all $s\in 
S$ or
 $ 
\Dim V_s\in K$   is invertible in $K$ for all $s\in S$.
 It follows from Lemma 10.1.1 that $L$ is
  biangular if (and only if) either $V_s=0$ for all $s\in 
S$ or  $\Dim V_s\in K$ does not 
depend on $s$ and is invertible in $K$. This gives examples of non-degenerate $\pi$-algebras
which
are not biangular.

If $L$ is biangular then an explicit 
computation of the homomorphism $\psi_1:L_1\to L_1$ defined in Section 7.1 
shows that it is a projection of the algebra $L_1=\oplus_s \End(V_s)$
onto its center $K^{\card (S)}$.  
 Note also that if the action of $\pi$ on $S$ is free and transitive then
$L=\End (\bigoplus_{s\in S}  V_s)$.

Replacing the direct sum with tensor product we obtain another interesting  
$\pi$-algebra. Namely, for $\alpha\in \pi$, set
$$R_\alpha= \bigotimes_{s\in S} \Hom (V_s, V_{\alpha (s)}).$$
The $K$-module $R_\alpha$ is additively generated by the elements $a=\otimes_{s  
} a_s$ where $a_s\in   \Hom (V_s, V_{\alpha(s)})$ for   $s\in S$.
If $  b=\otimes_{s} b_s\in R_\beta$  with $
b_s\in  \Hom (V_s, V_{\beta(s)})$
then the product $ab\in R_{\alpha\beta}$ is defined   by
$ab=\otimes_{s } (ab)_s$ where $(ab)_s=a_{\beta(s)} b_s\in  \Hom (V_s, 
V_{\alpha \beta (s)})$ for   $s\in S$. This multiplication is associative and 
makes $R=\oplus_{\alpha\in \pi} R_\alpha$ a $\pi$-algebra with unit 
$\otimes_{s } \id_{V_s}$. It follows from Lemma 10.1.1 
that this $\pi$-algebra always satisfies (7.1.i) and satisfies (7.1.ii) if and 
only if is either $V_s=0$ for all $s\in S$ or    $ 
\Dim V_s\in K$   is invertible in $K$ for all $s\in S$.

Assume that $K$ is a field of 
characteristic 0 and $V_s\neq 0$ for all $s\in S$. 
   Then  $R$ is biangular and 
 the homomorphism $\psi_1:R_1\to R_1$ defined in 
Section 7.1 is a  projection of the algebra 
$R_1=\otimes_{s\in S} \End (V_s)$ onto its 1-dimensional center.
Therefore the   lattice HQFT  determined by $R$
 is obtained  by $k^{\rho_0}$-rescaling from the primitive cohomological HQFT determined by an
element of  $H^2(\pi; K^*)$  where
 $k=\eta(1_R,1_R)=\Dim R_1=(\prod_s \Dim V_s)^2$.

\skipaline \noindent {\bf    10.3.  Crossed group-algebras via push-forward.  } 
Consider a group epimorphism $q:\pi'\to \pi$  with finite kernel, $H$,  lying in the center of
$\pi'$.  We can push  forward any crossed $\pi'$-algebra $L'$  
along $q$ to obtain  a Frobenius $\pi$-algebra $L=q_*(L')$, cf. 
Section 3.1.       Recall that
  $L=L'$ as Frobenius algebras. If $H=\Ker\,q$ is contained in the kernel of the  
action  $\varphi$ of $\pi'$ on $L'$,   then 
$\varphi$   induces  an action  of $\pi$ on $L$.  We claim that $L$ is a crossed
$\pi$-algebra.  Axioms (3.1.1) - (3.2.3) for  $L$ directly follows from the corresponding axioms
for $L'$. Let us check (3.2.4) for $L$.
 Let $\alpha, \beta \in \pi$ and   $c\in   L_{\alpha\beta \alpha^{-1}
\beta^{-1}}$. Note that for   $u\in q^{-1} (\alpha), v \in q^{-1}
(\beta)$, the commutator $u v u^{-1} v^{-1} 
$ does not depend on the choice of $u$ and $v$. Denote this
distinguished element of  $q^{-1}(\alpha\beta \alpha^{-1}
\beta^{-1})$ by $w_{\alpha,\beta}$.
To check equality (3.2.a), it suffices to consider the case where
$c\in L'_w$ where $w\in q^{-1}(\alpha\beta \alpha^{-1}
\beta^{-1})$. The homomorphism $c 
\varphi_{\beta}= c\varphi_{v}$ sends each direct summand $L'_u$ of
$L_\alpha$ (with $u\in q^{-1} (\alpha)$) into
$L'_{wvuv^{-1}}$.
Therefore 
$$\Tr (c\,
\varphi_{\beta}:L_{\alpha}\to L_{\alpha})= \sum_{u\in q^{-1} (\alpha), wvuv^{-1}=u} 
\Tr (c\,
\varphi_{v}:   L'_{u}\to    L'_{u})$$ 
$$
=\cases\sum_{u\in q^{-1} (\alpha)}  \Tr (c\, \varphi_{v}:   L'_{u}\to    L'_{u}),~ { 
{if}}\,\,\, w=w_{\alpha,\beta}, \\
  0,~    otherwise.\endcases
 $$
The assumption $\varphi(H)=1$ implies that the trace $ \Tr (c\, \varphi_{v}:   L'_{u}\to
   L'_{u})$ does not depend on the choice of $v$ in  $q^{-1}
(\beta)$. The same argument together with  formula (3.2.a) shows that this trace 
does not depend on the choice of $u$ in  $q^{-1}
(\alpha)$. Hence,
$$\Tr (c\,
\varphi_{\beta}:L_{\alpha}\to L_{\alpha})=  \cases (\card\,  H) \,   \Tr
(c\,
\varphi_{v}:   L'_{u}\to    L'_{u}),~ {  {if}}\,\,\, w=w_{\alpha,\beta}, \\
  0,~    otherwise,\endcases
 $$
where $u,v$ are arbitrary elements of   $  q^{-1} (\alpha),   q^{-1}
(\beta)$.
Similarly, 
the homomorphism $\varphi_{\alpha^{-1}} c= \varphi_{u^{-1}} c$ sends each direct
summand $L'_v$ of
$L_\beta$ (with $v\in q^{-1} (\beta)$) into
$L'_{u^{-1}wvu}$. Therefore 
$$\Tr ( \varphi_{\alpha^{-1}} c:L_{\beta}\to L_{\beta}) =\sum_{v\in q^{-1} (\beta),
v=u^{-1}wvu} 
\Tr (\varphi_{u^{-1}} c:   L'_{v}\to    L'_{v})$$
$$=
\cases (\card\,  H) \,   \Tr
(\varphi_{u^{-1}} c:   L'_{v}\to    L'_{v}),~ {  {if}}\,\,\, w=w_{\alpha,\beta}, \\
  0,~    otherwise.\endcases
$$
Now, axiom (3.2.4) for $L'$ implies   (3.2.4) for $L$.

 This construction can be combined with those discussed in Section 3.   Consider for
concreteness  the structure of a crossed
$ 
\pi'$-algebra in the group
ring $K[\pi']$
determined by the trivial class $0\in H^2(\pi', K^*)$. By the argument above, this
gives a structure of a crossed $\pi$-algebra  in the same ring $L=K[\pi']$. 
Here $\pi$ acts by conjugations and the inner product 
is given by 
$ (u,v)\mapsto \delta^u_{v^{-1}}$
for $u,v\in \pi'$ where $\delta$ is the Kronecker delta.
The crossed $\pi$-algebra  $L$ is semisimple if and only if the group ring 
$L_1=K[H]$ is semisimple. This is  for instance the case if   $K$ is
an algebraically closed field.  If $K=\bold  Q$ and $H$ is   non-trivial then $K[H]$ is not 
 semisimple.   It is curious to note that the inner product in $L$ admits deformations 
in the class of crossed $\pi$-algebras.
 Namely, choose  a non-degenerate symmetric
bilinear form $\mu: K[H]\times K[H]\to K$  on   $L_1=K[H]$ such that the pair
$(K[H],\mu)$ is a Frobenius algebra. (There are many such forms as it is easy to see for
cyclic $H$.) We define the inner product    $\eta_{\mu} : L
\otimes L
\to K
$ by  
   $\eta_{\mu}(L_{\alpha}\otimes L_{\beta})=0$ if $\alpha \beta\neq 1$ and    
 $\eta_{\mu}(a,b)=\mu (ab,1_L)$ for $a\in L_{\alpha}, b\in  L_{\alpha^{-1}}$ and any
$\alpha \in \pi$. Here $ab\in L_1=K[H]$ and $1_L\in K[H]$ 
  is  the unit element of 
$ L$.   It is easy to check that   $L$ with this inner product is  a crossed $\pi$-algebra.

If $K$ has a ring involution $k\mapsto
\overline k:K \to K$ 
 and  $\mu$ satisfies 
$\mu(u^{-1},v^{-1})=\overline {\mu(u,v)}$ for any $u,v \in H$ then  the antilinear involution in
$L=K[  \pi']$ sending each element of $    \pi'$ to its inverse defines a Hermitian
structure on  
$L$.

\skipaline \noindent {\bf    10.4.  Crossed group-algebras from    algebra
automorphisms.  }   Let us consider crossed $\pi$-algebras $L$ such that $L_\alpha=0$ for all
$\alpha\neq 1$.  Axioms (3.1.1) - (3.2.3)  amount to saying  that $(L_1,\eta)$ is a  commutative 
Frobenius algebra  with an action
  of $\pi$ by algebra automorphisms preserving $\eta$.  
Let us call such a pair $(L_1,\pi)$  a  $\pi$-F-algebra.
Axiom  (3.2.4)
means that the action of $\pi$ is traceless in the sense that $\Tr (c\,
\varphi_{\alpha}:L_{1}\to L_{1})=0$,
  for any $\alpha \neq  1, c\in   L_{1}$. 
Thus any traceless $\pi$-F-algebra $(L_1,\pi)$
gives rise to a crossed $\pi$-algebra such that $L_\alpha=0$ for  
$\alpha\neq 1$.  Note
that the tensor product $L_1\otimes
L'_1$ of  two  $\pi$-F-algebras $(L_1,\pi), (L'_1,\pi)$ is a
$\pi$-F-algebra. If  $(L_1,\pi)$ or $(L'_1,\pi)$  is traceless then $(L_1\otimes
L'_1,\pi)$  is   traceless.

$\pi$-F-algebras naturally arise in the study of groups of homeomorphisms.
Consider a closed connected oriented even-dimensional manifold $M$ and
set $L_1=L_1(M)=\oplus_{k \,even} H^k(M;\bold  Q)$. The product in $L_1$ is the  
cup-product and the form $\eta$ is defined by $\eta(a,b)=(a\cup b) ([M])$. 
Clearly, $L_1$ is a Frobenius algebra. Now, any group $\pi$ of orientation preserving
self-homeomorphisms of $M$ acts on $L_1$ via induced homomorphisms.
This action preserves the grading in $L_1$ and therefore  
$\Tr (c\,  {\alpha}_*:L_{1}\to L_{1})=0$   for all $\alpha\in \pi$ and 
 $c\in \oplus_{k\neq 0, k \,even} H^k(M;\bold  Q)\subset L_1$.
The only   requirement arises for $c=1\in H^0(M;\bold  Q)$ and consists
in the identity $\Tr ( {\alpha}_*:L_{1}\to L_{1})=0$ for all $\alpha\in \pi, \alpha\neq
1$. For instance, consider an orientation-reversing involution of the
$2n$-dimensional sphere $j_n:S^{2n}\to S^{2n}$. It is clear that the action of $j_n$ on 
$H^*(S^{2n};\bold  Q)$ is traceless. Therefore for any  orientation-reversing involution
$j$ of a  closed connected oriented even-dimensional manifold $M$, the product
$j_n\times j:S^{2n}\times M\to  S^{2n}\times M$ induces a traceless endomorphism
of $ L_1(S^{2n}\times M)$. In particular, $j_n\times j_m$ is  
a traceless endomorphism
of $ L_1(S^{2n}\times S^{2m})$. This   yields  
examples of  traceless   $\bold  Z/2\bold  Z$-F-algebras and hence of
crossed  $\bold  Z/2\bold  Z$-algebras.

This example can be extended in a slightly different direction. Consider 
 a closed connected oriented even-dimensional manifold $M$ and
a fixed point free group $\pi$ of orientation preserving
self-homeomorphisms of $M$. By the Lefschetz theorem,
the super-trace of the action of $\pi$ in the graded space $\oplus_{k } H^k(M;\bold  Q)$ 
is zero. If $M$ has only even-dimensional cohomology, this gives an example 
of a  traceless $\pi$-F-algebra. In general this suggests to consider
a wider class of crossed super-$\pi$-algebras.  

Deformations of Frobenius algebras form a subject of
a deep theory   based on the 
Witten-Dijkgraaf-Verlinde-Verlinde equation (see for instance [Du]). It would be interesting to
generalize this theory to crossed group-algebras.

\skipaline

\centerline{\bf References}

\skipaline

[At] Atiyah, M., Topological quantum field theories.  Publ.  Math.  IHES 68
(1989), 175-186.

[BP] Bachas, C., Petropoulos, P., Topological Models on the Lattice and a Remark on String
Theory Cloning. Comm.
Math.  Phys.  152 (1993), 191-202.

[Di] Dijkgraaf, R., A Geometrical Approach to Two-Dimensional Conformal Field
Theory, Ph.  D.  Thesis (Utrecht, 1989).

[Du] Dubrovin, B., Geometry of $2$D topological field theories.  Integrable
systems and quantum groups (Montecatini Terme, 1993), 120-348, Lecture Notes in
Math., 1620, Springer, Berlin, 1996.

[FQ] Freed, G., Quinn, F., Chern-Simons theory with finite gauge group.  Comm.
Math.  Phys.  156 (1993), 435-472.

[FHK] Fukuma, M., Hosono, S.,  Kawai, H., Lattice Topological Field Theory in Two Dimensions.
Comm.
Math.  Phys.  161 (1994), 157-175.

[HT] Hatcher, A., Thurston, W., A presentation for the mapping class group of a
closed orientable surface, Topology 19 (1980), 221-237.

[Hu]  Hungerford, T.W., Algebra.  New York, Springer, 1974

[Pi] Pierce, R., Associative algebras. Graduate Texts in Math. 88, Springer Verlag, 1982.
 
 [Tu] Turaev, V., Quantum invariants of knots and 3-manifolds.  Studies in
Mathematics 18, Walter de Gruyter, 588 p., 1994.

[TV] Turaev, V., Viro, O., State sum invariants of 3-manifolds
and quantum $6j$-symbols.   Topology 31 (1992), 865-902.

\skipaline

 {Institut de Recherche MathŽmatique AvancŽe, UniversitŽ
Louis Pasteur-CNRS,
7 rue RenŽ Descartes, 67084 Strasbourg Cedex, France}

turaev$\@$math.u-strasbg.fr

\end